\documentclass[twoside,a4paper,11pt]{amsart}

\usepackage[T1]{fontenc}
\usepackage{mathpazo} 
\usepackage{eulervm}

\usepackage[utf8]{inputenc}

\usepackage{fullpage}
\usepackage{amscd}
\usepackage{amsmath}
\usepackage{amssymb}
\usepackage{amsthm}
\usepackage{color}
\usepackage{graphicx}
\usepackage{hyperref}
\usepackage[noabbrev,capitalize]{cleveref}
\usepackage{mathrsfs}
\usepackage{mathtools}
\usepackage{pict2e}
\usepackage{stackrel}
\usepackage[T1]{fontenc}
\usepackage{tikz}
\usepackage{tikz-cd}
\usetikzlibrary{decorations.pathmorphing,decorations.markings,arrows,calc,shapes.geometric,arrows.meta,positioning}
\usepackage[utf8]{inputenc}
\usepackage{wasysym}
\usepackage[all]{xy}
\usepackage{xcolor}
\usepackage{xy}
\usepackage{epigraph}

\hypersetup{
    colorlinks,
    linkcolor=[rgb]{0.8, 0.5451, 0.3961},
    citecolor=[rgb]{0.8, 0.5451, 0.3961},
    urlcolor=[rgb]{0.8, 0.5451, 0.3961}
}

\usepackage{palatino}
\usepackage{mathpazo} 
\allowdisplaybreaks

\newtheorem{lemma}{Lemma}[section]
\newtheorem{theorem}[lemma]{Theorem}
\newtheorem{Corollary}[lemma]{Corollary}
\newtheorem{Proposition}[lemma]{Proposition}

\newtheorem{theoremintro}{Theorem}

\theoremstyle{definition}
\newtheorem{Definition}[lemma]{Definition}
\newtheorem{Remark}[lemma]{\sc Remark}

\newtheorem{Example}[lemma]{\sc Example}
\newtheorem{Notation}[lemma]{Notation}

\def\colim{\mathop{\mathrm{colim}}}

\newcommand{\kk}{\Bbbk}

\newcommand{\ac}{\scriptstyle \textrm{!`}}

\newcommand{\qi}{\xrightarrow{ \,\smash{\raisebox{-0.65ex}{\ensuremath{\scriptstyle\sim}}}\,}}
\newcommand{\lqi}{\xleftarrow{ \,\smash{\raisebox{-0.65ex}{\ensuremath{\scriptstyle\sim}}}\,}}

\newcommand{\C}{\mathcal{C}}

\newcommand{\I}{\mathcal{I}}
\newcommand{\PP}{\mathcal{P}}

\newcommand{\ucom}{u\mathcal{C}om}

\input cyracc.def
\font\tencyr=wncysc10
\def\cyr{\tencyr\cyracc}
\def\diracComb{\mbox{\cyr SH}}

\author{Victor Roca i Lucio}
\title{The integration theory of curved absolute $\mathcal{L}_\infty$-algebras}
\date{\today}

\address{Victor Roca i Lucio, Ecole Polytechnique Fédérale de Lausanne, EPFL,
CH-1015 Lausanne, Switzerland}
\email{\href{mailto:victor.rocalucio@epfl.ch}{victor.rocalucio@epfl.ch}}



\thanks{The author was partially supported by the project ANR-20-CE40-0016 HighAGT funded by the Agence Nationale pour la Recherche.}

\begin{document}
	
\begin{abstract}
In this article, we introduce the notion of a curved absolute $\mathcal{L}_\infty$-algebra, a structure that behaves like a curved $\mathcal{L}_\infty$-algebra where all infinite sums of operations are well-defined by definition. We develop their integration theory by introducing two new methods in integration theory: the complete bar construction and new transferred model category structures. They allow us to generalize all essential results of this theory quickly and from a conceptual point of view. We provide applications of our theory to rational homotopy theory, and show that curved absolute $\mathcal{L}_\infty$-algebras provide us with rational models for finite type nilpotent spaces without any pointed or connected assumptions. Furthermore, we show that the homology of rational spaces can be recovered as the homology of the complete bar construction. We also construct new smaller models for rational mapping spaces without any hypothesis on the source simplicial set. Another source of applications is deformation theory: on the algebraic side, we show that curved absolute $\mathcal{L}_\infty$-algebras are mandatory in order to encode the deformation complexes of $\infty$-morphisms of (co)-algebras. On the geometrical side, we construct a curved absolute $\mathcal{L}_\infty$-algebra from a derived affine stack and show that it encodes the formal geometry of any finite collection of points living in any finite field extension of the base field.
\end{abstract}

\maketitle

\setcounter{tocdepth}{1}

\tableofcontents

\section*{Introduction}
\textbf{Global picture.} In deformation theory, one studies how to "deform" structures (algebraic, geometric, etc) on an object. Here "deform" can have multiple meanings. Let $X$ be some type of object and $\mathcal{P}$ some type of structure. Ideally, one has a \textit{moduli space}  $\mathrm{Def}_X(\mathcal{P})$ where the points are the $\mathcal{P}$-structures that we can endow $X$ with. If one considers these structures up to some equivalences, and remembers the equivalences between the equivalences and so on and so forth, one has an $\infty$-groupoid (or \textit{moduli stack}). In this context, a first meaning of "deform" is, given a $\mathcal{P}$-structure on $X$, a point $x$ in $\mathrm{Def}_X(\mathcal{P})$, to look at the formal neighborhood of $x$ inside $\mathrm{Def}_X(\mathcal{P})$. This formal neighborhood encodes the \textit{infinitesimal deformations}, structures "infinitely close" to $x$.

\medskip

Over the years, it was noticed by D. Quillen, P. Deligne, V. Drinfel'd and many others that in all practical examples of these deformations, the formal neighborhoods were "encoded" by a differential graded (dg) Lie algebras. Examples include the work of K. Kodaira and D. Spencer in \cite{KodairaSpencer58}, encoding deformations of complex structures on manifolds or the work of M. Gerstenhaber in \cite{Gerstenhaber64}, encoding deformations of associative algebra structures on a vector space. This point of view suggests that for every point $x$ in $\mathrm{Def}_X(\mathcal{P})$, there should be a dg Lie algebra $\mathfrak{g}_{x,X}$ such that infinitesimal deformations correspond to elements in $\mathfrak{g}_{x,X}$ which satisfy the \textit{Maurer--Cartan equation}. This principle can be formulated as:
\[
\text{"Every infinitesimal deformation problem in characteristic zero is encoded by a dg Lie algebra."}
\]
This became recently a theorem by J. Lurie and J. P. Pridham, see \cite{Lurie11,Pridham10}. The method was to formalize what an "infinitesimal deformation" is using the notion of $\kk$-\textit{pointed formal moduli problem} and then show that their homotopy category is equivalent to the homotopy category of dg Lie algebras.

\medskip

Another point of view on the deformation theory of algebraic structures is given by operadic deformation complexes. In this context, given an operad $\mathcal{P}$ and a dg module $M$, one can construct an explicit dg Lie algebra where Maurer--Cartan elements are exactly the $\mathcal{P}$-algebra structures on $M$ and where the equivalences (called $\infty$-isotopies) between these structures are given by the action of a group called the \textit{gauge group}. Again, this amounts to studying the $\infty$-groupoid that is encoded by the operadic deformation problem. 

\medskip

\textit{Integration theory} amounts to constructing a way to recover these $\infty$-groupoids from their dg Lie algebras. Its most basic example can be traced back to Lie theory and the Baker--Campbell--Hausdorff formula which produces a group out of a nilpotent Lie algebra. This allows us to recover the gauge group in operadic deformation complexes. A first general approach to the integration of dg Lie algebras is given by the seminal work of V. Hinich in \cite{Hinich01}, using methods from D. Sullivan \cite{Sullivan77}. A refined version of the integration procedure was constructed by E. Getzler in \cite{Getzler09}, where he extended the integration procedure to nilpotent $\mathcal{L}_\infty$-algebras (homotopy Lie algebras). Inspired by the ideas of U. Buijs, Y. Félix, A. Murillo and D. Tanré in \cite{Buijandco}, generalized using operadic calculus, D. Robert-Nicoud and B. Vallette were able in \cite{robertnicoud2020higher} to give a new characterizations of Getzler's functor and obtained an adjunction 

\[
\begin{tikzcd}[column sep=7pc,row sep=3pc]
\mathsf{sSet} \arrow[r, shift left=1.1ex, "\mathcal{L}"{name=F}]      
&\mathcal{L}_\infty\text{-}\mathsf{alg}^{\mathsf{comp}}~. \arrow[l, shift left=.75ex, "\mathcal{R}"{name=U}]
\arrow[phantom, from=F, to=U, , "\dashv" rotate=-90]
\end{tikzcd}
\]
\vspace{0.1pc}

between simplicial sets and complete $\mathcal{L}_\infty$-algebras. This adjunction lies at the crossroad of three domains: Lie theory, deformation theory, and rational homotopy. Here the right adjoint functor $\mathcal{R}$ produces an \textit{integration functor} from complete $\mathcal{L}_\infty$-algebras to $\infty$-groupoids (Kan complexes). But the horn-fillers of this $\infty$-groupoid are a \textit{structure}, not a \textit{property}. They are given by explicit formulas which extend the classical Baker--Campbell--Hausdorff formula in the specific case of Lie algebras. Finally, the functor $\mathcal{L}$ is shown to produce faithful rational models for pointed connected finite type nilpotent spaces, and greatly simplifies the original approach of D. Quillen \cite{Quillen69} to rational homotopy. Similar rational models using complete dg Lie algebras where also constructed by Buijs--Félix--Murillo--Tanré in \cite{Buijandco}. 

\medskip

\textbf{Motivations.} The first goal of this article is to generalize the integration procedure to \textit{curved} $\mathcal{L}_\infty$-algebras. Curved Lie or curved $\mathcal{L}_\infty$-algebras are a notion more general than their dg counterparts: these objects are endowed with a distinguished element, the \textit{curvature}, which perturbs the relationships satisfied by the bracket and higher operations in a dg Lie or $\mathcal{L}_\infty$-algebra. In particular, the "differential" no longer squares to zero, thus the notion of quasi-isomorphism so useful in homotopical algebra disappears. 

\medskip

The main reason to consider curved $\mathcal{L}_\infty$-algebras is that they are the "non-pointed analogue" of $\mathcal{L}_\infty$-algebras. In a classical $\mathcal{L}_\infty$-algebra, the element $0$ is always a Maurer--Cartan element. This gives a canonical base point in the corresponding $\infty$-groupoid. It is no longer the case with curved $\mathcal{L}_\infty$-algebras. From the point of view of Lie theory, these objects behave like Lie groupoids instead of Lie groups. From the point of view of rational homotopy theory, these objects provide us with rational models for non necessarily pointed and non necessarily connected spaces.

\medskip

Developing the integration theory of curved $\mathcal{L}_\infty$-algebras has also many applications to deformation theory. In his PhD. thesis \cite{Nuiten19}, J. Nuiten showed that, if $A$ is a dg unital commutative algebra, then the homotopy category of $A$-\textit{pointed formal moduli problems} is equivalent to the homotopy category of dg Lie algebroids over $A$. Later, it was showed by D. Calaque, R. Campos and J. Nuiten in \cite{calaque2021lie} that the homotopy category of dg Lie algebroids over $A$ is equivalent to the homotopy category of certain curved $\mathcal{L}_\infty$-algebras over the de Rham algebra of $A$, under some assumptions on $A$. Therefore deformation problems "parametrized by $\mathrm{Spec}(A)$" can be encoded with curved $\mathcal{L}_\infty$-algebras. "Parametrized by $\mathrm{Spec}(A)$" means that instead of choosing a $\kk$-point in a moduli space $X$, one chooses a morphism $f: \mathrm{Spec}(A) \longrightarrow X$, and then looks at the formal neighborhood of $\mathrm{Spec}(A)$ inside of $X$. Another approach to parametrized deformation problems is also given by $\mathcal{L}_\infty$-spaces introduced by K. Costello in \cite{Costello}. These are families of curved $\mathcal{L}_\infty$-algebras parametrized by smooth manifolds. They were constructed in order to treat parametrized deformation problems arising in fundamental physics. On the other side, operadic deformation complexes of \textit{unital} algebraic structures form \textit{curved} Lie algebras as shown in \cite{lucio2022curved}. Furthermore, the space of $\infty$-morphisms between types of unital algebras or the space of $\infty$-morphisms between types of counital coalgebras are encoded, in both cases, by convolution curved $\mathcal{L}_\infty$-algebras. In all the above cases, having a good \textit{integration theory} for curved $\mathcal{L}_\infty$-algebras is of primordial importance. 

\medskip
 
\textbf{Framework and methods.} The second main goal of this article is to introduce new methods in integration theory. This is first done by introducing a new kind of algebraic objects that we call \textit{curved absolute} $\mathcal{L}_\infty$\textit{-algebras}. A curved absolute $\mathcal{L}_\infty$-algebra can be thought as curved $\mathcal{L}_\infty$-algebra where all infinite sums of operations have a well defined image \textit{without supposing any underlying topology}. Infinite sums appear naturally in the theory of $\mathcal{L}_\infty$-algebras with the Maurer--Cartan equation. So far, there have been two standard ways to deal with them: either restrict to nilpotent $\mathcal{L}_\infty$-algebras, or use the somewhat \textit{ad hoc} approach of changing the base category in order to consider objects with an underlying complete topology, so that these infinite sums \textit{converge}. Our approach solves these problems altogether. It also includes the classical examples of nilpotent Lie algebras and nilpotent $\mathcal{L}_\infty$-algebras in the sense of \cite{Getzler09}.

\medskip 

Considering curved absolute $\mathcal{L}_\infty$-algebras, which are encoded with \textit{a curved cooperad}, allows us to use the operadic calculus developed by B. Le Grignou and D. Lejay in \cite{grignoulejay18}. We introduce $u\mathcal{CC}_\infty$-coalgebras, a version of homotopy counital cocommutative coalgebras encoded by a dg operad as a Koszul dual notion of curved absolute $\mathcal{L}_\infty$-algebras. We then construct a \textit{complete bar-cobar} adjunction between these two types of objects. In particular, the \textit{complete bar construction}, given in terms of the non-conilpotent cofree $u\mathcal{CC}_\infty$-coalgebra, proves to be a new key ingredient in the theory. Using similar methods as in \cite{robertnicoud2020higher}, we can now construct a $u\mathcal{CC}_\infty$-coalgebra structure on the cellular chain functor and thus obtain the integration functor composing with the aforementioned adjunction. Notice that the methods of \textit{loc.cit.} where definitely not generalizable on the nose, since in the curved case, one has to consider \textit{counital coalgebras, which are not conilpotent}, and thus, which cannot be encoded by cooperads. This makes the usual operadic calculus fail in this context.

\medskip 

Another key ingredient to our results is the use of \textit{new model structures}. In the curved setting, the notion of quasi-isomorphism no longer makes sense, and one needs to find another notion of weak-equivalences. We solve this difficulty by first endowing the category of $u\mathcal{CC}_\infty$-coalgebras with a model category structure where weak-equivalences are given by quasi-isomorphisms and cofibrations by monomorphisms. Secondly, we transfer it onto the category of curved absolute $\mathcal{L}_\infty$-algebras, thus obtaining a meaningful notion of weak-equivalence in the curved context. This transferred model category structure proves to be very natural, as it allows us to generalize directly all the main results of integration theory to our setting. From a conceptual point of view, this allows us to fully embed integration theory in the context of homotopy theory. The functor that associates to a nilpotent Lie algebra a group using the Baker--Campbell--Hausdorff formula is now naturally the restriction of a right Quillen functor.

\medskip 

The last main new method used in this article is the extensive use that we make of the \textit{duality squares} that we introduced in \cite[Section 2]{lucio2022contra}. With them, one can understand homotopically the linear duality functors between types of coalgebras and types of algebras. This tool proves to be essential in order to obtain results in rational homotopy theory and in deformation theory. 

\medskip

\textbf{Main results.} We construct a cosimplicial $u\mathcal{CC}_\infty$-coalgebra $C_*^c(\Delta^\bullet)$ by endowing the cellular chains on the simplices with a $u\mathcal{CC}_\infty$-coalgebra. By Kan extension, this gives an adjunction

\[
\begin{tikzcd}[column sep=7pc,row sep=3pc]
            \mathsf{sSet} \arrow[r, shift left=1.1ex, "C^c_*(-)"{name=F}] 
            &u\mathcal{CC}_\infty\text{-}\mathsf{coalg}~, \arrow[l, shift left=.75ex, "\overline{\mathcal{R}}"{name=U}]
            \arrow[phantom, from=F, to=U, , "\dashv" rotate=-90]
\end{tikzcd}
\]

where $C^c_*(-)$ is the usual cellular chain functor with a $u\mathcal{CC}_\infty$-coalgebra structure. Since we endowed the category of $u\mathcal{CC}_\infty$-coalgebras with an intrinsic model category structure where weak-equivalences are given by quasi-isomorphisms, it is automatic to check that it forms a Quillen adjunction. By composing the above adjunction with the complete bar-cobar adjunction we get the following triangle. 

\begin{theoremintro}[Theorem \ref{thm: triangle of adjunctions}]
The following triangle of Quillen adjunctions commutes

\[
\begin{tikzcd}[column sep=5pc,row sep=2.5pc]
&\hspace{1pc}u\mathcal{CC}_\infty \textsf{-}\mathsf{coalg} \arrow[dd, shift left=1.1ex, "\widehat{\Omega}_{\iota}"{name=F}] \arrow[ld, shift left=.75ex, "\overline{\mathcal{R}}"{name=C}]\\
\mathsf{sSet}  \arrow[ru, shift left=1.5ex, "C^c_*(-)"{name=A}]  \arrow[rd, shift left=1ex, "\mathcal{L}"{name=B}] \arrow[phantom, from=A, to=C, , "\dashv" rotate=-70]
& \\
&\hspace{3pc}\mathsf{curv}~\mathsf{abs}~\mathcal{L}_\infty\textsf{-}\mathsf{alg}~. \arrow[uu, shift left=.75ex, "\widehat{\mathrm{B}}_{\iota}"{name=U}] \arrow[lu, shift left=.75ex, "\mathcal{R}"{name=D}] \arrow[phantom, from=B, to=D, , "\dashv" rotate=-110] \arrow[phantom, from=F, to=U, , "\dashv" rotate=-180]
\end{tikzcd}
\]
\end{theoremintro}

Here $\widehat{\Omega}_{\iota}$ is a completed cobar construction, and $\widehat{\mathrm{B}}_{\iota}$ is a new complete Bar construction. Having this constructions allows to define the \textit{integration functor} $\mathcal{R}$ as the composite of the functor $\overline{\mathcal{R}}$ with the complete bar construction. One of the main novelties of this is the fact that $\mathcal{R}$ defined in this way is automatically \textit{right Quillen functor} with respect to transferred model category structure on curved absolute $\mathcal{L}_\infty$-algebra. This simple fact directly implies the required properties that a well-behaved integration functor needs to satisfy. 

\begin{theoremintro}[Theorem \ref{thm: propriétés de l'intégration}]
\leavevmode

\begin{enumerate}
\item For any curved absolute $\mathcal{L}_\infty$-algebra $\mathfrak{g}$, the simplicial set $\mathcal{R}(\mathfrak{g})$ is a Kan complex.

\item For any degree-wise epimorphism $f: \mathfrak{g} \twoheadrightarrow \mathfrak{h}$ of curved absolute $\mathcal{L}_\infty$-algebras, the induced map
\[
\mathcal{R}(f): \mathcal{R}(\mathfrak{g}) \twoheadrightarrow \mathcal{R}(\mathfrak{h})
\]
is a fibrations of simplicial sets. 

\item The functor $\mathcal{R}(-)$ preserves weak-equivalences. In particular, it sends any graded quasi-isomorphism between complete curved absolute $\mathcal{L}_\infty$-algebras to a weak homotopy equivalence of simplicial sets.

\end{enumerate}
\end{theoremintro}

The fact that $\mathcal{R}(\mathfrak{g})$ is an $\infty$-groupoid is quintessential to integration theory as developed in \cite{Hinich01} and in \cite{Getzler09}. The second statement is a generalization of one of the main theorems of \cite{Getzler09}. The third point, the \textit{homotopy invariance} of the integration functor, is the generalization of the celebrated Goldman--Milson Theorem of  \cite{goldmanmillson} and its extension by Dolgushev--Rogers to $\mathcal{L}_\infty$-algebras in \cite{dolgushevrogers}. Also, since $\mathcal{R}$ is a right Quillen functor, it automatically commutes with homotopy limits and therefore satisfies \textit{descent} in the sense of \cite{hinichdescent}. Moreover, the adjunction $\mathcal{L} \dashv \mathcal{R}$ is shown to be a non-abelian generalization of the Dold-Kan correspondence.

\medskip

We then develop higher absolute Lie theory. We introduce and characterize gauge equivalences for Maurer--Cartan elements in this setting. We construct higher Baker--Campbell--Hausdorff products for curved absolute $\mathcal{L}_\infty$-algebras and we show that these are given by the same explicit formulae as in \cite{robertnicoud2020higher}. We generalize A. Berglund's Theorem of \cite{Berglund15} to the case of curved absolute $\mathcal{L}_\infty$-algebras with completely new proof based on a simple computation of the $u\mathcal{CC}_\infty$-coalgebra structures on the cellular chains of the spheres. This gives a way to compute the homotopy groups of $\mathcal{R}(\mathfrak{g})$ using the \textit{twisted homology} of $\mathfrak{g}$. It is worth mentioning that Getzler's and Hinich's approach to integration theory \textit{would not have worked} in this context, since the tensor product of a commutative algebra with a curved absolute $\mathcal{L}_\infty$-algebra has no natural curved absolute $\mathcal{L}_\infty$-algebra structure.

\medskip

The third section is devoted to constructing rational homotopy models using the new functor $\mathcal{L}$. The fact that curved $\mathcal{L}_\infty$-algebras are the "non-pointed analogue" of $\mathcal{L}_\infty$-algebras allows to obtain rational Lie models for \textit{non necessarily pointed nor connected} finite type nilpotent spaces using our functor $\mathcal{L}$.

\begin{theoremintro}[Theorem \ref{thm: modèles d'homotopie rationnel type fini}]
Let $X$ be a finite type nilpotent simplicial set. The unit of the adjunction

\[
\eta_X: X \qi \mathcal{R}\mathcal{L}(X) 
\]

is a rational homotopy equivalence.
\end{theoremintro}

We show that if a curved absolute $\mathcal{L}_\infty$-algebra $\mathfrak{g}_X$ is a rational model for a space $X$, then we can recover the \textit{homology} of $X$ by computing the \textit{homology of the complete bar construction of} $\mathfrak{g}_X$. This result can be viewed as Eckmann-Hilton dual to the A. Berglund's theorem mentioned before. One can think of the complete bar construction of a curved absolute $\mathcal{L}_\infty$-algebras as a \textit{higher Chevalley-Eilenberg complex}.

\medskip

We use the theory of B. Le Grignou developed in the forthcoming paper \cite{grignou2022mapping} to get a convolution curved absolute $\mathcal{L}_\infty$-algebra structure on the space of graded morphisms between a $u\mathcal{CC}_\infty$-coalgebra and a curved absolute $\mathcal{L}_\infty$-algebra. 

\begin{theoremintro}[Theorem \ref{thm: vrai thm mapping spaces}]
Let $\mathfrak{g}$ be a curved absolute $\mathcal{L}_\infty$-algebra and let $X$ be a simplicial set. There is a weak-equivalence of Kan complexes

\[
\mathrm{Map}(X, \mathcal{R}(\mathfrak{g})) \qi \mathcal{R}\left(\mathrm{hom}(C^c_*(X),\mathfrak{g})\right)~,
\]
\vspace{0.25pc}

which is natural in $X$ and in $\mathfrak{g}$, where $\mathrm{hom}(C^c_*(X),\mathfrak{g})$ denotes the convolution curved absolute $\mathcal{L}_\infty$-algebra. 
\end{theoremintro}

Notice that, for the first time, there are no assumptions on $X$ nor on $\mathfrak{g}$. Furthermore, we can replace $C^c_*(X)$ by the homology $\mathrm{H}_*(X)$ in a non-functorial way using the homotopy transfer theorem. If $Y$ is a finite type nilpotent simplicial set, the above theorem gives an explicit model for the mapping space of $X$ and the $\mathbb{Q}$-localization of $Y$. Futhermore, this model is constructed using the cellular chains on $X$ (or even its homology), hence it smaller than other rational models constructed before. See for instance \cite{Berglund15, BuijMapping, LazarevMapping}.

\medskip

Finally, we also develop the integration theory of (non-curved) absolute $\mathcal{L}_\infty$-algebras and show that the adjunction $\mathcal{L}_* \dashv \mathcal{R}_*$ thus obtained is exactly a \textit{pointed} version of the one constructed before in a precise sense. Furthermore, we show that the rational models for pointed connected finite type nilpotent spaces constructed using absolute $\mathcal{L}_\infty$-algebras are \textit{local} with respect to quasi-isomorphisms, and this allows us to compare our approach with the others present in the literature. We generalize this to the curved case, constructing an \textit{intrinsic} model category structure on curved absolute $\mathcal{L}_\infty$-algebras, transferred from simplicial sets with the Kan--Quillen model structure; we show that it is again a localization of the previous model structure, and that it still encodes rational finite type nilpotent homotopy types.

\medskip

The last section is devoted to applications in deformation theory, where the new notion of curved absolute $\mathcal{L}_\infty$-algebras proves to be mandatory. We construct convolution curved $\mathcal{L}_\infty$-algebras that encode $\infty$-morphisms between unital types of algebras as their Maurer--Cartan elements. This gives new notions of homotopies between $\infty$-morphisms and opens the way to simplicial enrichments. The very nature of convolution algebras is "absolute"; they always admit infinite sums of operations without the need of an underlying filtration to make them converge. Therefore same methods could be applied to $\infty$-morphisms of coalgebras over dg operads or, more generally, $\infty$-morphisms of gebras over properads as defined by E. Hoffbeck, J. Leray and B. Vallette in \cite{hoffbeck2019properadic}. 

\medskip

Finally, we explore the geometrical properties of curved absolute $\mathcal{L}_\infty$-algebras. Starting from a derived affine stack $A$, that is, a dg unital commutative algebra, we construct an explicit curved absolute $\mathcal{L}_\infty$-algebra model $\mathfrak{g}_A$. We then show that the formal neighborhood of any $\mathbb{L}$-point of $A$ can be recovered from $\mathfrak{g}_A$, where $\mathbb{L}$ is a finite field extension of our base field $\kk$. Furthermore, contrary to the $\kk$-pointed formal moduli approach, these points are not "specified in advance".

\begin{theoremintro}[Theorem \ref{thm: modèle géométrique sur les points Artiniens}]
Let $A$ be a dg $u\mathcal{C}om$-algebra. Let $B$ be a dg Artinian algebra. There is a weak-equivalence of simplicial sets
\[
\mathrm{Spec}(A)(B) \simeq \mathcal{R}(\mathrm{hom}( \mathbb{R}(\mathrm{Res}_\varepsilon B)^\circ, \mathfrak{g}_A))~,
\]
natural in $B$ and in $A$, where $\mathbb{R}(-)^\circ$ is the derived Sweedler dual functor.
\end{theoremintro}

The trick behind the above theorem lies in the definition of a dg Artinian algebra. Since we are working with \textit{curved} absolute $\mathcal{L}_\infty$-algebras, dg Artinian algebras need not be augmented. In fact, we simply define them as dg unital commutative algebras in non-negative degrees with degree-wise finite dimensional and bounded homology. Conversely, any curved absolute $\mathcal{L}_\infty$-algebra defines a deformation functor from dg Artinian algebras to $\infty$-groupoids. We view these results as a first step in the theory of "non-pointed" formal moduli problems.

\subsection*{Acknowledgments}
I would like to thank my former PhD. advisor Bruno Vallette for the numerious discussion we had and for his careful readings of this paper. Many thanks to Brice Le Grignou for sharing his work on Mapping Coalgebras before it is available on Arxiv, and for taking the time to explain their content to me. I'm grateful to Joan Bellier-Millès, Damien Calaque, Ricardo Campos, Geoffroy Horel, Joost Nuiten, and Bertrand Toën for interesting discussions. I would also like to thank the Mittag-Leffler Institut which provided great working conditions in order to finish this paper. This paper was written during my PhD. thesis at the Université Sorbonne Paris Nord, I would like to thank its great mathematical community. Finally, it is my pleasure to thank the referee for a great number of insightful comments and suggestions that greatly helped improving this paper.

\subsection*{Conventions} We adopt the same conventions as in \cite{lucio2022contra}. We work over a field $\kk$ of characteristic $0$. The base category will be the symmetric monoidal category $\left(\mathsf{pdg}\textsf{-}\mathsf{mod}, \otimes, \kk\right)$ of pre-differential graded $\kk$-modules with the tensor product. We work with \textit{homological conventions}, therefore the (pre)-differentials will always be of degree $-1$. We denote $(\mathsf{pdg}~\mathbb{S}\textsf{-}\mathsf{mod}, \circ, \I)$ the monoidal category of graded $\mathbb{S}$-modules endowed with the composition product $\circ$. This framework is explained with more details in \cite{lucio2022curved}.

\newpage

\section{Curved absolute $\mathcal{L}_\infty$-algebras}
In this section, we introduce the notion of curved \textit{absolute} $\mathcal{L}_\infty$-algebras. This is a new type of curved $\mathcal{L}_\infty$-algebras, encoded by a conilpotent curved cooperad. They posses a much richer algebraic structure than usual curved $\mathcal{L}_\infty$-algebras. Infinite sums of operations have a well-defined image by definition. In particular, the Maurer--Cartan equation is always defined. The rest of this article is devoted to the study of this new notion and its applications. Notice that our degree conventions will correspond to \textit{shifted} curved $\mathcal{L}_\infty$-algebras. This \textit{shifted convention} will be implicit from now on. Curved absolute algebras appear as the Koszul dual of a specific model for \textit{non-necessarily conilpotent} $\mathbb{E}_{\infty}$-coalgebras, which we call $u\mathcal{CC}_{\infty}$-coalgebras. We construct a complete bar-cobar adjunction that relates these two types of algebraic objects and show that it is a Quillen equivalence.

\subsection{Curved absolute $\mathcal{L}_\infty$-algebras.}
Let $\ucom$ be the operad encoding unital commutative algebras and let $\epsilon: \ucom \longrightarrow \I$ be the canonical morphism of $\mathbb{S}$-modules given by the identity on $\ucom(1) \cong \kk~.$ We denote $\mathrm{B}^{\mathrm{s.a}}\ucom$ its semi-augmented bar construction with respect to this canonical semi-augmentation. We refer to Appendix \ref{Section: Appendix B} for more details on this particular construction.

\begin{Definition}[Curved absolute $\mathcal{L}_\infty$-algebra]
A \textit{curved absolute} $\mathcal{L}_\infty$\textit{-algebra} $\mathfrak{g}$ amounts to the data $(\mathfrak{g},\gamma_\mathfrak{g},d_\mathfrak{g})$ of a curved $\mathrm{B}^{\mathrm{s.a}}\ucom$-algebra. 
\end{Definition}

\begin{Remark}
Brief recollections on the notion of curved algebras over cooperads are given in \cite[Section 1]{lucio2022contra}. We also refer to \cite{mathez} for a more thorough exposition of the theory developed in \cite{grignoulejay18}.
\end{Remark}

Let us unravel this definition. The data of a curved absolute $\mathcal{L}_\infty$-algebra structure on a pdg module $(\mathfrak{g},d_\mathfrak{g})$ amounts to the data of a morphism of pdg modules
\[
\gamma_\mathfrak{g}: \prod_{n \geq 0} \mathrm{Hom}_{\mathbb{S}_n}\left(\mathrm{B}^{\mathrm{s.a}}\ucom(n), \mathfrak{g}^{\otimes n}\right) \longrightarrow \mathfrak{g}~,
\]
which satisfies the conditions of \cite[Definition 1.12 and 1.22]{lucio2022contra}. This map admits a simpler description.

\begin{lemma}\label{lemma: iso invariants avec les coinvariants}
There is an isomorphism of pdg modules
\[
\prod_{n \geq 0} \mathrm{Hom}_{\mathbb{S}_n}\left(\mathrm{B}^{\mathrm{s.a}}\ucom(n), \mathfrak{g}^{\otimes n}\right) \cong \prod_{n \geq 0} \widehat{\Omega}^{\mathrm{s.a}}\ucom^*(n) ~\widehat{\otimes}_{\mathbb{S}_n}~ \mathfrak{g}^{\otimes n}~,
\]
natural in $\mathfrak{g}$, where $\widehat{\otimes}$ denotes the completed tensor product with respect to the canonical filtration on the complete cobar construction.
\end{lemma}

\begin{proof}
By definition, $\mathrm{B}^{\mathrm{s.a}}\ucom$ is a conilpotent curved cooperad. Let us denote by $\mathscr{R}_\omega \mathrm{B}^{\mathrm{s.a}}$ the $\omega$-term of its coradical filtration. There is an isomorphism of conilpotent curved cooperads
\[
\mathrm{B}^{\mathrm{s.a}}\ucom \cong \colim_{\omega}\mathscr{R}_\omega \mathrm{B}^{\mathrm{s.a}}\ucom~.
\]
Notice that for every $n \geq 0$ and every $\omega \geq 0$, $\mathscr{R}_\omega \mathrm{B}^{\mathrm{s.a}}\ucom(n)$ is degree-wise finite dimensional. Thus we have
\begin{align*}
\mathrm{Hom}_{\mathbb{S}_n}\left(\mathrm{B}^{\mathrm{s.a}}\ucom(n), \mathfrak{g}^{\otimes n}\right)
&\cong  \mathrm{Hom}_{\mathbb{S}_n}\left(\colim_{\omega}\mathscr{R}_\omega \mathrm{B}^{\mathrm{s.a}}\ucom(n), \mathfrak{g}^{\otimes n}\right) \\
&\cong \lim_{\omega} \mathrm{Hom}_{\mathbb{S}_n}\left(\mathscr{R}_\omega \mathrm{B}^{\mathrm{s.a}}\ucom(n), \mathfrak{g}^{\otimes n}\right)~.
\end{align*}
for all $n \geq 0$. By Lemma \ref{lemma: Bucom dual lin}, there is an isomorphism 
\[
\left(\mathscr{R}_\omega \mathrm{B}^{\mathrm{s.a}}\ucom\right)^{*} \cong \widehat{\Omega}^{\mathrm{s.c}}\ucom^*/\mathscr{F}_\omega \widehat{\Omega}^{\mathrm{s.c}}\ucom^*~,
\]
where $\mathscr{F}_\omega \widehat{\Omega}^{\mathrm{s.c}}\ucom^*$ denotes the $\omega$-term of its canonical filtration as an absolute partial operad. Therefore there are isomorphisms
\begin{align*}
\lim_{\omega} \mathrm{Hom}_{\mathbb{S}_n} \left(\mathscr{R}_\omega \mathrm{B}^{\mathrm{s.a}}\ucom(n), \mathfrak{g}^{\otimes n}\right)
&\cong \lim_{\omega} \left(\widehat{\Omega}^{\mathrm{s.c}}\ucom^*/\mathscr{F}_\omega \widehat{\Omega}^{\mathrm{s.c}}\ucom^*(n) \otimes \mathfrak{g}^{\otimes n}\right)^{\mathbb{S}_n} \\
&\cong \left(\widehat{\Omega}^{\mathrm{s.c}}\ucom^*(n) ~\widehat{\otimes}~ \mathfrak{g}^{\otimes n}\right)^{\mathbb{S}_n} \\
&\cong ~\widehat{\Omega}^{\mathrm{s.c}}\ucom^*(n) ~\widehat{\otimes}_{\mathbb{S}_n}~ \mathfrak{g}^{\otimes n}~.
\end{align*}
Notice that the last isomorphism identifies invariants with coinvariants, which is possible because of the characteristic zero assumption. Nevertheless, this identification carries non-trivial coefficients, see Remark \ref{Remark: renormalization}.
\end{proof}

\begin{Notation}
Let $\mathrm{CRT}_n^\omega$ denote the set of \textit{corked rooted trees} of arity $n$ and with $\omega$ vertices. A corked rooted tree is a rooted tree where vertices either have at least two incoming edges or zero incoming edges, which are called \textit{corks}. The arity of a corked rooted tree is the number of non-corked leaves. The unique rooted tree of arity $n$ with one vertex is called the $n$-corolla, denoted by $c_n$. Notice that for each $n$, the set $\mathrm{CRT}_n$ of corked rooted trees of arity $n$ is infinite. Likewise, for each $\omega$, the set of corked rooted trees with $\omega$ vertices is also infinite. Nevertheless, the set $\mathrm{CRT}_n^\omega$ is finite. The only corked rooted tree of weight $0$ is the trivial tree of arity one with zero vertices.
\end{Notation}

\begin{Proposition}
Let $(\mathfrak{g},d_\mathfrak{g})$ be a pdg module with a basis $\left\{ g_b ~|~b \in B \right\}~.$ The pdg module
\[
\prod_{n \geq 0} \widehat{\Omega}^{\mathrm{s.a}}\ucom^*(n) ~\widehat{\otimes}_{\mathbb{S}_n} ~ \mathfrak{g}^{\otimes n}
\]
admits a basis given by double series on $n$ and $\omega$ of corked rooted trees in $\mathrm{CRT}_n^\omega$ labeled by the basis elements of $\mathfrak{g}^{\otimes n}$. These basis elements can be written as

\[
\left\{ \sum_{n\geq 0} \sum_{\omega \geq 0} \sum_{\tau \in \mathrm{CRT}_n^\omega} \sum_{i \in \mathrm{I}_\tau} \lambda_\tau^{(i)} \tau \left(g_{i_1}^{(i)}, \cdots, g_{i_n}^{(i)}\right) \right\}~,
\]

where $\mathrm{I}_\tau$ is a finite set, $\lambda_\tau$ is a scalar in $\kk$ and $(i_1,\cdots,i_n)$ is in $B^n$. Here $\tau(g_{i_1}, \cdots, g_{i_n})$ is given by the rooted tree $\tau$ with input leaves decorated by the elements $\tau(g_{i_1}, \cdots, g_{i_n})$. The degree of $\tau(g_{i_1}, \cdots, g_{i_n})$ is $-\omega -1 + |g_{i_1}| + \cdots + |g_{i_n}|~.$ An example of such decoration is given by

\medskip

\begin{center}
\includegraphics[width=130mm,scale=1.3]{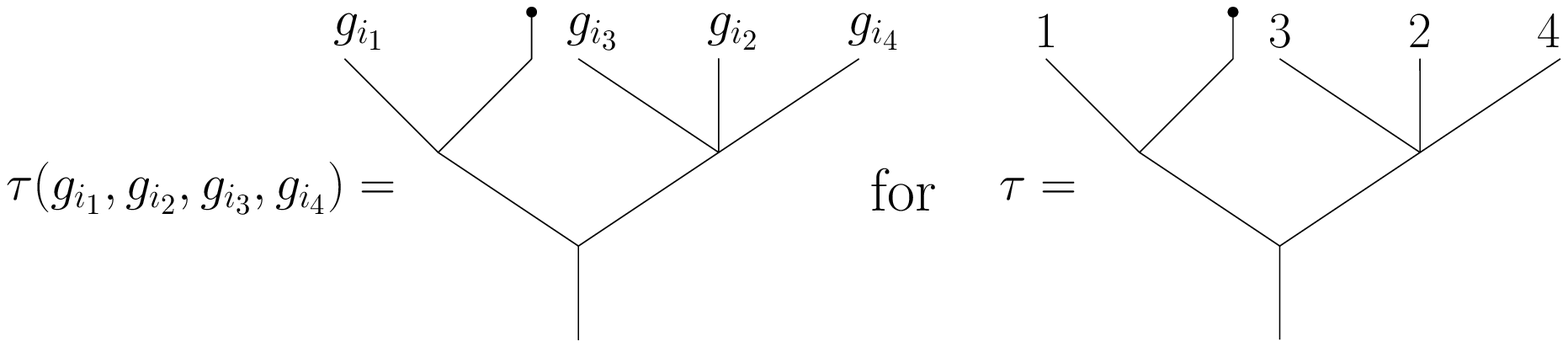}.
\end{center}

The pre differential of the pdg module is given by the sum of two terms: the first term is applies the pre-differential $d_\mathfrak{g}$ on each of the labels $g_{i_j}$ of a corked rooted tree $\tau$, the second splits each vertex into two vertices in all the possible ways including the splitting with a cork. 
\end{Proposition}

\begin{proof}
Corked rooted trees form a basis of the pdg $\mathbb{S}$-module $\widehat{\Omega}^{\mathrm{s.a}}\ucom^*$, hence the results follows by direct inspection.
\end{proof}

\begin{Remark}[Renormalization]\label{Remark: renormalization}
Let $\tau$ be a \textit{rooted tree without corks}. It can be written as $c_m \circ (\tau_1, \cdots, \tau_m)$. We define recursively the following coefficient for rooted trees:

\[
\mathcal{E}(c_n) \coloneqq n! \quad \text{and} \quad \mathcal{E}\left(c_m \circ (\tau_1, \cdots, \tau_m)\right) \coloneqq |\mathrm{G}_m| \prod_{i=1}^m\mathcal{E}(\tau_i )~,
\]
where $G_m$ is the subgroup of permutations in $\mathbb{S}_m$ sending a branch to $\tau_j$ to an isomorphic branch $\tau_{j'}$. Notice that these coefficients appear in Theorem \ref{theom: formulas for the pre differentials of mc}. More generally, when $\tau$ is a corked rooted tree, the coefficient $\mathcal{E}(\tau)$ is given by the symmetry group of the tree. Considering these coefficients, the following series 

\[
\left\{ \sum_{n\geq 0} \sum_{\omega \geq 0} \sum_{\tau \in \mathrm{CRT}_n^\omega} \sum_{i \in \mathrm{I}_\tau} \mathcal{E}(\tau).\lambda_\tau^{(i)} \tau \left(g_{i_1}^{(i)}, \cdots, g_{i_n}^{(i)}\right) \right\}~,
\]

form a basis of 

\[
\prod_{n \geq 0} \mathrm{Hom}_{\mathbb{S}_n}\left(\mathrm{B}^{\mathrm{s.a}}\ucom(n), \mathfrak{g}^{\otimes n}\right) \cong  \prod_{n \geq 0}   \left(\widehat{\Omega}^{\mathrm{s.a}}\ucom^*(n) ~\widehat{\otimes} ~ \mathfrak{g}^{\otimes n} \right)^{\mathbb{S}_n}~.
\]
\end{Remark}

\medskip

\textbf{Structural data.} Thus the data of a curved absolute $\mathcal{L}_\infty$-algebra structure on a pdg module $(\mathfrak{g},d_\mathfrak{g})$ amounts to the data of a morphism of pdg modules
\[
\gamma_\mathfrak{g}: \prod_{n \geq 0} \widehat{\Omega}^{\mathrm{s.a}}\ucom^*(n) ~\widehat{\otimes}_{\mathbb{S}_n} ~ \mathfrak{g}^{\otimes n} \longrightarrow \mathfrak{g}~,
\]
which satisfies the conditions of \cite[Definition 1.12 and 1.22]{lucio2022contra}. In particular, for any infinite sum 
\[
\sum_{n\geq 0} \sum_{\omega \geq 0} \sum_{\tau \in \mathrm{CRT}_n^\omega} \sum_{i \in \mathrm{I}_\tau} \lambda_\tau^{(i)} \tau \left(g_{i_1}^{(i)}, \cdots, g_{i_n}^{(i)}\right)
\]
there is a well-defined image 
\[
\gamma_\mathfrak{g}\left( \sum_{n\geq 0} \sum_{\omega \geq 0} \sum_{\tau \in \mathrm{CRT}_n^\omega} \sum_{i \in \mathrm{I}_\tau} \lambda_\tau^{(i)} \tau \left(g_{i_1}^{(i)}, \cdots, g_{i_n}^{(i)}\right)\right)~,
\]
in the pdg module $\mathfrak{g}$. This \textbf{does not presuppose} an underlying topology on the pdg module $\mathfrak{g}$.

\begin{Notation}
Since elements in $\widehat{\Omega}^{\mathrm{s.c}}\ucom^*/\mathscr{F}_\omega \widehat{\Omega}^{\mathrm{s.c}}\ucom^*(n) \otimes \mathfrak{g}^{\otimes n}$ might not be \textit{simple} tensors, this makes the sums $\sum_{i \in \mathrm{I}_\tau}$ appear in the formulae. From now on, we will omit these sums and work only with simple tensors. Computations for infinite sums of non-simple tensors are, \textit{mutatis mutandis}, the same.
\end{Notation}

\textbf{Pdg condition.} Let us make explicit each of the conditions satisfied by the structural morphism $\gamma_\mathfrak{g}$. It is a morphism of pdg modules if only if, we have that

\begin{small}

\begin{align}\label{pdg condition}
\gamma_\mathfrak{g}\left(\sum_{\substack{n \geq 0 \\ \omega \geq 0}} \sum_{\tau \in \mathrm{CRT}_n^\omega} \lambda_\tau d_{\widehat{\Omega}\ucom}(\tau)(g_{i_1}, \cdots, g_{i_n})\right) + \gamma_\mathfrak{g}\left(\sum_{\substack{n \geq 0 \\ \omega \geq 0}} \sum_{\tau \in \mathrm{CRT}_n^\omega}\sum_{j=0}^n (-1)^j \lambda_\tau\tau(g_{i_1}, \cdots, d_\mathfrak{g}(g_{i_j}), \cdots, g_{i_n})\right) = \nonumber \\
d_\mathfrak{g}\left(\gamma_\mathfrak{g}\left(\sum_{\substack{n \geq 0 \\ \omega \geq 0}} \sum_{\tau \in \mathrm{CRT}_n^\omega} \lambda_\tau\tau(g_{i_1}, \cdots, g_{i_n})\right) \right)~.
\end{align}

\end{small}

\textbf{Associativity condition.} The associativity condition on $\gamma_\mathfrak{g}$ imposed by the diagram of \cite[Definition 1.12]{lucio2022contra} is equivalent to the following equality

\begin{small}

\begin{align}\label{associativity condition}
\gamma_\mathfrak{g} \left(\sum_{\substack{k\geq 0 \\ \omega \geq 0}} \sum_{\tau \in \mathrm{CRT}_k^\omega} \lambda_\tau \tau \left(\gamma_\mathfrak{g}\left(\sum_{\substack{i_1 \geq 0 \\ \omega_1 \geq 0}} \sum_{\tau_1 \in \mathrm{CRT}_{i_1}^{\omega_1}} \lambda_{\tau_1}\tau_1(\bar{g}_{i_1} )\right) ,\cdots, \gamma_\mathfrak{g} \left(\sum_{\substack{i_k \geq 0 \\ \omega_k \geq 0}} \sum_{\tau_k \in \mathrm{CRT}_{i_k}^{\omega_k}} \lambda_{\tau_k} \tau_k(\bar{g}_{i_k} )\right)\right)\right) = \nonumber \\
\gamma_\mathfrak{g} \left(\sum_{n \geq 0} \sum_{k \geq 0} \sum_{i_1 + \cdots + i_k =n} \sum_{\substack{\omega^{\mathrm{tot}} \geq 0  \\ \omega +\omega_1 + \cdots +\omega_k = \omega^{\mathrm{tot}}}} \sum_{\substack{\tau \in \mathrm{CRT}_k^\omega \\ \tau_j \in \mathrm{CRT}_{i_j}^{\omega_j}}} \lambda_\tau(\lambda_{\tau_1}\cdots\lambda_{\tau_k})~ \tau \circ (\tau_1,\cdots, \tau_k)(\bar{g}_{i_1},\cdots,\bar{g}_{i_k}) \right)~, 
\end{align}
\end{small}

where $\bar{g}_{i_j}$ denotes an $i_j$-tuple of elements of $\mathfrak{g}^{\otimes i_j}$ and where $\tau \circ (\tau_1,\cdots, \tau_k)$ denotes the corked rooted tree obtained by grafting $(\tau_1,\cdots, \tau_k)$ onto the leaves of $\tau~.$ If $\gamma_\mathfrak{g}$ satisfies the above conditions, it endows $\mathfrak{g}$ with a pdg $\text{B}^{\mathrm{s.a}}\ucom$-algebra structure.

\begin{Definition}[Elementary operations of a curved absolute $\mathcal{L}_\infty$-algebra]
Let $(\mathfrak{g},\gamma_\mathfrak{g}, d_\mathfrak{g})$ be a curved absolute $\mathcal{L}_\infty$-algebra. The \textit{elementary operations} of $\mathfrak{g}$ are the family of symmetric operations of degree $-1$ given by
\[
\left\{l_n \coloneqq \gamma_\mathfrak{g}(c_n(-,\cdots,-)): \mathfrak{g}^{\odot n} \longrightarrow \mathfrak{g}~\right\}
\]
for all $n \neq 1$. 
\end{Definition}

\textbf{Curved condition.} The condition on $\gamma_\mathfrak{g}$ imposed by the diagram of \cite[Definition 1.22]{lucio2022contra} amounts in this case to the following equation on the elementary operations

\begin{equation}\label{curved condition}
d_\mathfrak{g}^2(g) = l_2(l_0,g)~.
\end{equation}

\textbf{Conclusion.} A curved absolute $\mathcal{L}_\infty$-algebra structure $\gamma_\mathfrak{g}$ on a pdg module $(\mathfrak{g},d_\mathfrak{g})$ amounts to the data of a degree $0$ map 
\[
\gamma_\mathfrak{g}: \prod_{n \geq 0} \widehat{\Omega}^{\mathrm{s.a}}\ucom^*(n) ~\widehat{\otimes}_{\mathbb{S}_n} ~ \mathfrak{g}^{\otimes n} \longrightarrow \mathfrak{g}~,
\]
satisfying conditions \ref{pdg condition}, \ref{associativity condition} and \ref{curved condition}.

\begin{Remark}[Warning]
In general, if $(\mathfrak{g},\gamma_\mathfrak{g}, d_\mathfrak{g})$ is a curved absolute $\mathcal{L}_\infty$-algebra, then 
\[
\gamma_\mathfrak{g}\left(\sum_{n\geq 0} \sum_{\omega \geq 0} \sum_{\tau \in \mathrm{CRT}_n^\omega} \lambda_\tau \tau(g_{i_1}, \cdots, g_{i_n})\right) \neq \sum_{n\geq 0} \sum_{\omega \geq 1} \sum_{\tau \in \mathrm{CRT}_n^\omega} \lambda_\tau \gamma_\mathfrak{g} \left(\tau(g_{i_1}, \cdots, g_{i_n})\right)~,
\]
as the latter expression \textit{is not well-defined} since infinite sums of elements \textit{in} $\mathfrak{g}$ are not well-defined in general.
\end{Remark}

\textbf{Morphisms.} Let $(\mathfrak{g},\gamma_\mathfrak{g}, d_\mathfrak{g})$ and $(\mathfrak{h},\gamma_\mathfrak{h}, d_\mathfrak{h})$ be two curved absolute $\mathcal{L}_\infty$-algebras and let $f: \mathfrak{g} \longrightarrow \mathfrak{h}$ be a morphism of pdg modules. The condition for $f$ to be a morphism of curved absolute $\mathcal{L}_\infty$-algebras can be written as 
\[
f\left(\gamma_\mathfrak{g}\left(\sum_{n\geq 0} \sum_{\omega \geq 0} \sum_{\tau \in \mathrm{CRT}_n^\omega} \lambda_\tau \tau(g_{i_1}, \cdots, g_{i_n})\right)\right) = \gamma_\mathfrak{h}\left(\sum_{n\geq 0} \sum_{\omega \geq 0} \sum_{\tau \in \mathrm{CRT}_n^\omega} \lambda_\tau \tau(f(g_{i_1}), \cdots, f(g_{i_n}))\right)~.
\]

Any curved absolute $\mathcal{L}_\infty$-algebra structure can be restricted along elementary operations to induce a curved $\mathcal{L}_\infty$-algebra structure in the standard sense.

\begin{Proposition}\label{prop: restriction functor}
There is a restriction functor 
\[
\mathrm{Res}: \mathsf{curv}~\mathsf{abs}~\mathcal{L}_\infty\text{-}\mathsf{alg} \longrightarrow \mathsf{curv}~\mathcal{L}_\infty\text{-}\mathsf{alg}~,
\]
from the category of curved absolute $\mathcal{L}_\infty$-algebras to the category of curved $\mathcal{L}_\infty$-algebras which is faithful.
\end{Proposition}

\begin{proof}
Let $(\mathfrak{g},\gamma_\mathfrak{g}, d_\mathfrak{g})$ be a curved absolute $\mathcal{L}_\infty$-algebra, we can restrict the structural map

\[
\begin{tikzcd}
\mathrm{Res}(\gamma_\mathfrak{g}): \displaystyle\bigoplus_{n \geq 0} \widehat{\Omega}^{\mathrm{s.c}}\ucom^*(n) \otimes_{\mathbb{S}_n} \mathfrak{g}^{\otimes n} \arrow[r,"\iota_\mathfrak{g}"]
&\displaystyle \prod_{n \geq 0} \widehat{\Omega}^{\mathrm{s.c}}\ucom^*(n) ~\widehat{\otimes}_{\mathbb{S}_n}~ \mathfrak{g}^{\otimes n} \arrow[r,"\gamma_\mathfrak{g}"]
&\mathfrak{g}
\end{tikzcd}
\]

along the natural inclusion $\iota_\mathfrak{g}$. It endows $\mathfrak{g}$ with a curved $\widehat{\Omega}^{\mathrm{s.c}}\ucom^*$-algebra structure. By Proposition \ref{Prop: Curved algebras over OmegauCom}, this is equivalent to a classical curved $\mathcal{L}_\infty$-algebra structure in the sense of Definition \ref{def: classical curved linfty alg}. Any morphism of curved absolute $\mathcal{L}_\infty$-algebras is in particular a morphism of curved $\mathcal{L}_\infty$-algebras; moreover, the functor $\mathrm{Res}$ is faithful.
\end{proof}

\begin{Proposition}\label{prop: adjunction curved absolute et curved non-absolute}
The restriction functor admits a left adjoint $\mathrm{Abs}$, which called the absolute envelope of a curved $\mathcal{L}_\infty$-algebra. Therefore there is an adjunction 

\[
\begin{tikzcd}[column sep=7pc,row sep=3pc]
\mathsf{curv}~\mathcal{L}_\infty\text{-}\mathsf{alg} \arrow[r, shift left=1.1ex, "\mathrm{Abs}"{name=F}]      
&\mathsf{curv}~\mathsf{abs}~\mathcal{L}_\infty\text{-}\mathsf{alg}~, \arrow[l, shift left=.75ex, "\mathrm{Res}"{name=U}]
\arrow[phantom, from=F, to=U, , "\dashv" rotate=-90]
\end{tikzcd}
\]
\end{Proposition}

\begin{proof}
Both categories are presentable, by \cite[Theorem 1.66]{AdamekRosicky} the functor $\mathrm{Res}$ admits as left-adjoint since it preserves all limits and it is accessible. This last assertion follows from the fact that sufficiently filtered colimits of curved absolute $\mathcal{L}_\infty$-algebras can be computed in the ground category, as a category of algebras over an accessible monad. 
\end{proof}

\begin{Remark}
For more comparison statements between the classical notion of curved $\mathcal{L}_\infty$-algebras as defined in Definition \ref{def: classical curved linfty alg} and this new notion of curved absolute $\mathcal{L}_\infty$-algebra, we refer to \cite[Chapter 3, Appendix A]{mathez}.
\end{Remark}

\textbf{Completeness.} The coradical filtration on the conilpotent curved cooperad $\text{B}^{\mathrm{s.a}}\ucom$ induces a canonical filtration on any curved absolute $\mathcal{L}_\infty$-algebra. See \cite[Definition 1.14]{lucio2022contra}. In this case, this \textit{canonical filtration} on $\mathfrak{g}$ is given by
\[
\mathrm{W}_\omega\mathfrak{g} \coloneqq \mathrm{Im}\left(\gamma_\mathfrak{g} \vert_{\mathscr{F}_\omega} : \prod_{n \geq 0} \mathscr{F}_\omega\widehat{\Omega}^{\mathrm{s.a}}\ucom^*(n) ~ \widehat{\otimes}_{\mathbb{S}_n}~  \mathfrak{g}^{\otimes n} \longrightarrow \mathfrak{g} \right)~,
\]
for $\omega \geq 0$, where $\mathscr{F}_\omega\widehat{\Omega}^{\mathrm{s.a}}\ucom^*(n)$ is spanned by corked rooted trees with a number of vertices greater or equal to $\omega$. Therefore an element $g$ is in $\mathrm{W}_{\omega_0}\mathfrak{g}$ if and only if it can be written as
\[
g = \gamma_\mathfrak{g}\left(\sum_{n\geq 0} \sum_{\omega \geq \omega_0} \sum_{\tau \in \mathrm{CRT}_n^\omega} \lambda_\tau \tau(g_{i_1}, \cdots, g_{i_n})\right)~.
\]
One can check that $\mathfrak{g}/\mathrm{W}_\omega\mathfrak{g}$ has an unique curved absolute $\mathcal{L}_\infty$-algebra structure given by the quotient structural map. 

\begin{Remark}
Condition \ref{pdg condition} implies that $d_\mathfrak{g}(\mathrm{W}_\omega\mathfrak{g}) \subset \mathrm{W}_\omega\mathfrak{g}$ and therefore the pre-differential $d_\mathfrak{g}$ is \textit{continuous} with respect to the canonical filtration.
\end{Remark}

\begin{Definition}[Complete curved absolute $\mathcal{L}_\infty$-algebra] 
Let $(\mathfrak{g},\gamma_\mathfrak{g}, d_\mathfrak{g})$ be a curved absolute $\mathcal{L}_\infty$-algebra. It is \textit{complete} if the canonical epimorphism
\[
\varphi_\mathfrak{g}: \mathfrak{g} \twoheadrightarrow \lim_{\omega} \mathfrak{g}/\mathrm{W}_\omega\mathfrak{g}
\]
is an isomorphism of curved absolute $\mathcal{L}_\infty$-algebras.
\end{Definition}

\begin{Remark}
We have that $\varphi_\mathfrak{g}$ is always an epimorphism by \cite[Proposition 1.16]{lucio2022contra}, this phenomenon is explained in \cite[Remark 1.17]{lucio2022contra}.
\end{Remark}

Let $(\mathfrak{g},\gamma_\mathfrak{g}, d_\mathfrak{g})$ be a complete curved absolute $\mathcal{L}_\infty$-algebra. One can write
\[
\gamma_\mathfrak{g}\left(\sum_{n\geq 0} \sum_{\omega \geq 0} \sum_{\tau \in \mathrm{CRT}_n^\omega} \lambda_\tau \tau(g_{i_1}, \cdots, g_{i_n})\right) =  
\sum_{\omega \geq 0} \gamma_\mathfrak{g}\left(\sum_{n\geq 0} \sum_{\tau \in \mathrm{CRT}_n^\omega} \lambda_\tau \tau(g_{i_1}, \cdots, g_{i_n})\right)~,
\]
using the fact that the canonical filtration on $\mathfrak{g}$ is complete.

\begin{Definition}[Maurer--Cartan element of a curved absolute $\mathcal{L}_\infty$-algebra]
Let $(\mathfrak{g},\gamma_\mathfrak{g}, d_\mathfrak{g})$ be a curved absolute $\mathcal{L}_\infty$-algebra. A \textit{Maurer--Cartan element} $\alpha$ is an element in $\mathfrak{g}$ of degree $0$ which satisfies the following equation
\[
\gamma_\mathfrak{g}\left(\sum_{n \geq 0,~n\neq 1} \frac{c_n(\alpha, \cdots, \alpha)}{n!} \right) + d_\mathfrak{g}(\alpha) = 0~.
\]
The set of Maurer--Cartan elements in $\mathfrak{g}$ is denoted by $\mathcal{MC}(\mathfrak{g})~.$
\end{Definition}

\begin{Remark}
The coefficients $1/n!$ appear because of the isomorphism between invariants and coinvariants in Lemma \ref{lemma: iso invariants avec les coinvariants}. When one considers the renormalization of this equation as in Remark \ref{Remark: renormalization}, these coefficients disappear.
\end{Remark}

\begin{Remark}
In a curved absolute $\mathcal{L}_\infty$-algebra, the element $0$ is not in general a Maurer--Cartan element since:
\[
\gamma_\mathfrak{g}\left(\sum_{n \geq 0,~n\neq 1} \frac{c_n(0, \cdots, 0)}{n!} \right) + d_\mathfrak{g}(0) = \gamma_\mathfrak{g}(c_0) \neq 0~.
\]
Therefore the set $\mathcal{MC}(\mathfrak{g})$ is not canonically pointed, neither non-necessarily non-empty.
\end{Remark}

\begin{Remark}
Notice that if one takes $\alpha$ to be a degree $0$ element of $\mathfrak{g}$, the value of 
\[
\gamma_\mathfrak{g}\left(\sum_{n \geq 0,~n\neq 1} \frac{c_n(\alpha, \cdots, \alpha)}{n!} \right)~.
\]
is not determined, \textit{in general}, by the values of partial sums of the operations $l_n(\alpha,\cdots, \alpha)$. It can happen that these partial sums are all non-trivial and do not converge to zero in the canonical topology, and that $\alpha$ is still a Maurer--Cartan element. For example, one can consider the free (non-curved) absolute $\mathcal{L}_\infty$-algebra on a single generator $x$ in degree $0$, modulo the ideal generated by the element representing the (non-curved) Maurer--Cartan equation. Then $x$ becomes a Maurer--Cartan element, however all partial sums are non-trivial and lie in weight $1$, without converging to zero.
\end{Remark}

\begin{lemma}\label{lemma: splitting of the Maurer--Cartan equation}
Let $(\mathfrak{g},\gamma_\mathfrak{g}, d_\mathfrak{g})$ be a complete curved absolute $\mathcal{L}_\infty$-algebra. Any degree $0$ element $\alpha$ in $\mathrm{W}_1 \mathfrak{g}$ satisfies
\[
\gamma_\mathfrak{g}\left(\sum_{n \geq 0,~n\neq 1} \frac{c_n(\alpha, \cdots, \alpha)}{n!} \right) = \sum_{n \geq 0,~n\neq 1} \frac{l_n(\alpha, \cdots, \alpha)}{n!} ~.
\]
\end{lemma}

\begin{proof}
Since $\alpha$ is in $\mathrm{W}_1 \mathfrak{g}$, it can be written as 
\[
\alpha = \gamma_\mathfrak{g}\left(\sum_{k\geq 0} \sum_{\omega \geq 1} \sum_{\tau \in \mathrm{CRT}_k^\omega} \lambda_\tau \tau(\bar{g}_k)\right)~.
\]
where $\bar{g}_k$ is an $k$-tuple in $\mathfrak{g}$. A straightforward computation using the associativity condition \ref{associativity condition} concludes the proof.
\end{proof}

\textbf{Graded homology groups.} Let $(\mathfrak{g},\gamma_\mathfrak{g}, d_\mathfrak{g})$ be a curved absolute $\mathcal{L}_\infty$-algebra. Recall that Condition \ref{curved condition} says that $d_\mathfrak{g}^2(-) = l_2(l_0,-)$. Thus 
\[
d_\mathfrak{g}^2(\mathrm{W}_\omega \mathfrak{g}) \subseteq \mathrm{W}_{\omega+1} \mathfrak{g}~.
\]
This implies that
\[
\mathrm{gr}_{\omega}(\mathfrak{g}) \coloneqq \mathrm{W}_\omega \mathfrak{g}/\mathrm{W}_{\omega+1} \mathfrak{g}~,
\]
forms a chain complex endowed with the differential induced by $d_\mathfrak{g}$, for all $\omega \geq 0$.

\begin{Remark}
Notice that any morphism $f: \mathfrak{g} \longrightarrow \mathfrak{h}$ of curved absolute $\mathcal{L}_\infty$-algebras is \textit{continuous} with respect to the canonical filtration, i.e: $f(\mathrm{W}_\omega \mathfrak{g}) \subset \mathrm{W}_\omega \mathfrak{h}$. Therefore the morphism of dg modules $\mathrm{gr}_{\omega}(f)$ is always well-defined.
\end{Remark}

\begin{Definition}[Graded quasi-isomorphisms] 
Let $f: \mathfrak{g} \longrightarrow \mathfrak{h}$ be a morphism between two curved absolute $\mathcal{L}_\infty$-algebras. It is a \textit{graded quasi-isomorphism} if 
\[
\begin{tikzcd}
\mathrm{gr}_{\omega}(f): \mathrm{gr}_{\omega}(\mathfrak{g}) \arrow[r]
&\mathrm{gr}_{\omega}(\mathfrak{h})
\end{tikzcd}
\]
is a quasi-isomorphism of dg modules, for all $\omega \geq 0$.
\end{Definition}

\textbf{Absolute $\mathcal{L}_\infty$-algebras.} Let us develop the case of absolute $\mathcal{L}_\infty$-algebras. Algebraically, this case corresponds to the case where the curvature $l_0$ is zero. They are a \textit{pointed} version of their curved counterparts, as they always admit $0$ as a canonical Maurer--Cartan element. 

\begin{Definition}[Absolute $\mathcal{L}_\infty$-algebra]
An \textit{absolute} $\mathcal{L}_\infty$\textit{-algebra} $\mathfrak{g}$ amounts to the data $(\mathfrak{g},\gamma_\mathfrak{g},d_\mathfrak{g})$ of a dg $\mathrm{B}\mathcal{C}om$-algebra. 
\end{Definition}

Analogously to the curved case, the structural data of an absolute $\mathcal{L}_\infty$-algebra admits an explicit description. An absolute $\mathcal{L}_\infty$-algebra structure on a \textit{dg} module $(\mathfrak{g},d_\mathfrak{g})$ amounts to the data of a morphism of \textit{dg} modules
\[
\gamma_\mathfrak{g}: \prod_{n \geq 0} \Omega\mathcal{C}om^*(n) \otimes_{\mathbb{S}_n} \mathfrak{g}^{\otimes n} \longrightarrow \mathfrak{g}~,
\]
which satisfies the analogous conditions to those explained before in the curved setting. They carry an analogous canonical filtration induced by the weight of \textit{rooted trees} (without corks) and an analogous definition of a \textit{complete} absolute $\mathcal{L}_\infty$-algebra can be made.

\begin{Proposition}
There is an adjunction 
\[
\begin{tikzcd}[column sep=7pc,row sep=3pc]
\mathsf{curv}~\mathsf{abs}~\mathcal{L}_\infty\text{-}\mathsf{alg} \arrow[r, shift left=1.1ex, "(-)_*"{name=F}]      
&\mathsf{abs}~\mathcal{L}_\infty\text{-}\mathsf{alg}~, \arrow[l, shift left=.75ex, "\mathrm{U}"{name=U}]
\arrow[phantom, from=F, to=U, , "\dashv" rotate=-90]
\end{tikzcd}
\]
between curved absolute $\mathcal{L}_\infty$-algebras and absolute $\mathcal{L}_\infty$-algebras, where the functor $\mathrm{U}$ is fully faithful. The essential image of $\mathrm{U}$ is given by curved absolute $\mathcal{L}_\infty$-algebras such that the curvature $l_0$ is zero. 
\end{Proposition}

\begin{proof}
There is a canonical morphism of operads $\mathcal{C}om \longrightarrow u\mathcal{C}om$, which induces a morphism of conilpotent curved cooperads $\mathrm{B}\mathcal{C}om \longrightarrow \mathrm{B}^{\mathrm{s.a}}u\mathcal{C}om$ (any dg cooperad is a curved cooperad with zero curvature). In turn, this morphism induces the above adjunction by \cite[Section 7.2]{grignoulejay18}. It is straightforward to check that the functor $\mathrm{U}$ is fully faithful and to identify its essential image. 
\end{proof}

Finally, let us compare them with their classical counterparts, that is, with $\mathcal{L}_\infty$-algebras.

\begin{Proposition}\label{prop: adjonction absolute non-courbe et non-absolute}
There is an adjunction  

\[
\begin{tikzcd}[column sep=7pc,row sep=3pc]
\mathcal{L}_\infty\text{-}\mathsf{alg} \arrow[r, shift left=1.1ex, "\mathrm{Abs}"{name=F}]      
&\mathsf{abs}~\mathcal{L}_\infty\text{-}\mathsf{alg}~, \arrow[l, shift left=.75ex, "\mathrm{Res}"{name=U}]
\arrow[phantom, from=F, to=U, , "\dashv" rotate=-90]
\end{tikzcd}
\]

between $\mathcal{L}_\infty$-algebras and absolute $\mathcal{L}_\infty$-algebras. 
\end{Proposition}

\begin{proof}
The same arguments as used in Proposition \ref{prop: adjunction curved absolute et curved non-absolute}.
\end{proof}

\begin{Example}[Nilpotent $\mathcal{L}_\infty$-algebras]\label{Example: nilpotent stuff}
Let $\mathfrak{g}$ be a nilpotent $\mathcal{L}_\infty$-algebra in the sense of \cite{Getzler09}. These algebras can both be considered as $\mathcal{L}_\infty$-algebras or as absolute $\mathcal{L}_\infty$-algebras. The adjunction of Proposition \ref{prop: adjonction absolute non-courbe et non-absolute} then becomes an equivalence between nilpotent objects on both sides, as, in this case, the absolute envelope functor $\mathrm{Abs}$ is the identity on nilpotent objects, by the same arguments as in \cite[Remark 3.27]{lucio2022contra}. 
\end{Example}

\subsection{On $u\mathcal{CC}_\infty$-coalgebras}
We study the Koszul dual notion of dg coalgebras over the dg operad $\Omega\mathrm{B}^{\mathrm{s.a}}\ucom$. We denote these coalgebras $u\mathcal{CC}_\infty$-coalgebras, since they correspond to counital cocommutative coalgebras relaxed up to homotopy. 

\begin{Definition}[$u\mathcal{CC}_\infty$-coalgebra]
A $u\mathcal{CC}_\infty$\textit{-coalgebra} $C$ is the data $(C,\Delta_C,d_C)$ of a dg $\Omega\mathrm{B}^{\mathrm{s.a}}\ucom$-coalgebra.
\end{Definition}

\begin{lemma}
The data of a $u\mathcal{CC}_\infty$-coalgebra structure $\Delta_C$ on a dg module $(C,d_C)$ is equivalent to the data of elementary decomposition maps
\[
\Big\{ \Delta_\tau: C \longrightarrow C^{\otimes n} \Big\}
\]

of degree $\omega - 1$ for all corked rooted trees $\tau$ in $\mathrm{CRT}_n^\omega$ for all $n \geq 0$ and $\omega \geq 0$. These operations are subject to the following condition: let $E(\tau)$ denote the set of internal edges of $\tau$, then 
\[
\sum_{e \in E(\tau)} \Delta_{\tau_{(1)}} \circ_e \Delta_{\tau_{(2)}} - (-1)^{\omega -1} \sum_{e \in E(\tau)} \Delta_{\tau^{e}} + d_{C^{\otimes n}} \circ \Delta_\tau - (-1)^{\omega -1} \Delta_{\tau} \circ d_C = 0~,
\]
where in the first term the corked rooted tree is split along $e$ into $\tau_{(1)}$ and $\tau_{(2)}$ and where $\tau^{e}$ denotes the corked rooted tree obtained by contracting the internal edge $e$. 
\end{lemma}

\begin{proof}
The data of a $u\mathcal{CC}_\infty$-coalgebra structure $\Delta_C$ on a dg module $(C,d_C)$ is equivalent to the data of a morphism of dg operads
\[
\Delta_C: \Omega\mathrm{B}^{\mathrm{s.a}}\ucom \longrightarrow \mathrm{Coend}_C~,
\]
which in turn is equivalent to the data of a curved twisting morphism
\[
\delta_C: \mathrm{B}^{\mathrm{s.a}}\ucom \longrightarrow \mathrm{Coend}_C~.
\]
Thus the operations $\Delta_\tau$ are given by $\delta_C(\tau)$ and their relationships are given by the Maurer--Cartan equation that $\delta_C$ satisfies.
\end{proof}

\begin{Notation}
We denote by $\mathrm{PCRT}_n^{(\omega, \nu)}$ the set of partitioned corked rooted trees with $\omega$ vertices, $\nu$ blocks, of arity $n$. Such a tree amounts to the data of a corked rooted tree $\tau$ together with the data of a partition of the set of vertices of $\tau$ into $\nu$ different blocks connected blocks.
\end{Notation}

\begin{Remark}
The dg operad $\Omega\mathrm{B}^{\mathrm{s.a}}\ucom$ admits a basis given by partitioned corked rooted trees. The full structure of a $u\mathcal{CC}_\infty$-coalgebra is given by decomposition maps 
\[
\Big\{ \Delta_{\tau}: C \longrightarrow C^{\otimes n} \Big\}
\]

of degree $\omega - \nu$ if $\tau$ is in $\mathrm{PCRT}_n^{(\omega, \nu)}$. Nevertheless, any such decomposition map is obtained as the composition of the elementary decomposition maps given by the sub-corked rooted trees contained inside each partition.
\end{Remark}

Let $(C,d_C)$ be a dg module endowed with a family of decomposition maps $\{\Delta_\tau\}$ for any partitioned corked rooted tree $\tau$. This data allows us to construct a morphism of dg modules 

\[
\begin{tikzcd}[column sep=4pc,row sep=-0.5pc]
C \arrow[r,"\Delta_{C}"]
&\displaystyle \prod_{n \geq 0} \mathrm{Hom}_{\mathbb{S}_n}\left(\Omega\mathrm{B}^{\mathrm{s.a}}\ucom(n), C^{\otimes n}\right)\\
c \arrow[r,mapsto]
&\Big[ \mathrm{ev}_c: \tau \mapsto \Delta_{\tau}(c) \Big]~.
\end{tikzcd}
\]

which satisfies the axioms of \cite[Definition 1.17]{lucio2022contra}. 

\subsection{Model category structures} Since the dg operad encoding $u\mathcal{CC}_\infty$-coalgebras is cofibrant, the category of $u\mathcal{CC}_\infty$-coalgebras admits a canonical model category structure, left-transferred from the category of dg modules. See \cite[Section 1]{lucio2022contra} for an overview of the general results used in this subsection.

\begin{Proposition}\label{prop: model structure on uCC coalgebras}
There is a model category structure on the category of $u\mathcal{CC}_\infty$-coalgebras left-transferred along the cofree-forgetful adjunction
\[
\begin{tikzcd}[column sep=7pc,row sep=3pc]
\mathsf{dg}~\mathsf{mod} \arrow[r, shift left=1.1ex, "\mathscr{C}(u\mathcal{CC}_\infty)(-)"{name=F}]      
&u\mathcal{CC}_\infty\textsf{-}\mathsf{coalg}~, \arrow[l, shift left=.75ex, "\mathrm{U}"{name=U}]
\arrow[phantom, from=F, to=U, , "\dashv" rotate=90]
\end{tikzcd}
\]
where 
\begin{enumerate}

\item the class of weak-equivalences is given by quasi-isomorphisms,

\item the class of cofibrations is given by degree-wise monomorphisms,

\item the class of fibrations is given by right lifting property with respect to acyclic cofibrations.
\end{enumerate}
\end{Proposition}

\begin{proof}
Consequence of the results of \cite[Section 8]{grignoulejay18} applied to this particular case, since $\Omega\mathrm{B}^{\mathrm{s.a}}\ucom$ is cofibrant in the model category of dg operads.
\end{proof}

Furthermore, using the canonical curved twisting morphism 
\[
\iota: \mathrm{B}^{\mathrm{s.a}}\ucom \longrightarrow \Omega\mathrm{B}^{\mathrm{s.a}}\ucom~,
\]
we can construct a complete bar-cobar adjunction relating $u\mathcal{CC}_\infty$-coalgebras and curved absolute $\mathcal{L}_\infty$-algebras.

\begin{Proposition}[Complete bar-cobar adjunction]
The curved twisting morphism $\iota$ induces a complete bar-cobar adjunction
\[
\begin{tikzcd}[column sep=5pc,row sep=3pc]
            u\mathcal{CC}_\infty\text{-}\mathsf{coalg} \arrow[r, shift left=1.1ex, "\widehat{\Omega}_{\iota}"{name=F}] & \mathsf{curv}~\mathsf{abs}~\mathcal{L}_\infty\text{-}\mathsf{alg}^{\mathsf{comp}} \arrow[l, shift left=.75ex, "\widehat{\mathrm{B}}_{\iota}"{name=U}]
            \arrow[phantom, from=F, to=U, , "\dashv" rotate=-90]
\end{tikzcd}
\]
between the category of $u\mathcal{CC}_\infty$-coalgebras and the category of complete curved absolute $\mathcal{L}_\infty$-algebras. 
\end{Proposition} 

\begin{proof}
Consequence of the results explained in \cite[Section 9]{grignoulejay18} applied to this particular case.
\end{proof}

Using this adjunction, one can transfer the model category structure on $u\mathcal{CC}_\infty$-coalgebras in order to endow curved absolute $\mathcal{L}_\infty$-algebras with a model category structure.

\begin{theorem}\label{thm: Equivalence the Quillen CC infini et absolues}
There is a model category structure on the category of complete curved absolute $\mathcal{L}_\infty$-algebras right-transferred along the complete bar-cobar adjunction
\[
\begin{tikzcd}[column sep=5pc,row sep=3pc]
            u\mathcal{CC}_\infty\text{-}\mathsf{coalg} \arrow[r, shift left=1.1ex, "\widehat{\Omega}_{\iota}"{name=F}] & \mathsf{curv}~\mathsf{abs}~\mathcal{L}_\infty\text{-}\mathsf{alg}^{\mathsf{comp}} \arrow[l, shift left=.75ex, "\widehat{\mathrm{B}}_{\iota}"{name=U}]
            \arrow[phantom, from=F, to=U, , "\dashv" rotate=-90]
\end{tikzcd}
\]
where 
\begin{enumerate}

\item the class of weak-equivalences is given by morphisms $f$ such that $\widehat{\mathrm{B}}_\iota(f)$ is a quasi-isomorphism,

\item the class of fibrations is given by morphisms $f$ such that $\widehat{\mathrm{B}}_\iota(f)$ is a fibration,

\item  and the class of cofibrations is given by left lifting property with respect to acyclic fibrations.

\medskip
\end{enumerate}
Furthermore, this complete bar-cobar adjunction is a Quillen equivalence. 
\end{theorem}

\begin{proof}
Consequence of the results explained in \cite[Section 11]{grignoulejay18} applied to this particular case.
\end{proof}

One can give a complete characterization of the fibrations in this model category structure as follows.

\begin{Proposition}
A morphism of complete curved absolute $\mathcal{L}_\infty$-algebras is a fibration if and only if it is a degree-wise epimorphism.
\end{Proposition}

\begin{proof}
Follows from \cite[Proposition 10.16]{grignoulejay18}.
\end{proof}

\begin{Remark}
In particular, all complete curved absolute $\mathcal{L}_\infty$-algebras are fibrant in this model structure.
\end{Remark}

There is a particular kind of weak-equivalences between curved absolute $\mathcal{L}_\infty$-algebras which admit an easy description.

\begin{Proposition}
Let $f: \mathfrak{g} \longrightarrow \mathfrak{h}$ be a graded quasi-isomorphism between two complete curved absolute $\mathcal{L}_\infty$-algebras. Then 
\[
\widehat{\mathrm{B}}_\iota(f): \widehat{\mathrm{B}}_\iota \mathfrak{g}  \qi \widehat{\mathrm{B}}_\iota \mathfrak{h} 
\]
is a quasi-isomorphism of $u\mathcal{CC}_\infty$-coalgebras, and $f$ is a weak-equivalence of complete curved absolute $\mathcal{L}_\infty$-algebras.
\end{Proposition}

\begin{proof}
Follows from \cite[Theorem 10.25]{grignoulejay18}.
\end{proof}

\textbf{The non-unital case.} In a similar way, absolute $\mathcal{L}_\infty$-algebras are Koszul dual to non-unital homotopy cocommutative coalgebras, which are encoded by the dg operad $\Omega\mathrm{B}\mathcal{C}om$. We will denote this operad by $\mathcal{CC}_\infty$ from now on. There is a bar-cobar adjunction 
\[
\begin{tikzcd}[column sep=5pc,row sep=3pc]
            \mathcal{CC}_\infty\text{-}\mathsf{coalg} \arrow[r, shift left=1.1ex, "\widehat{\Omega}_{\iota}^\flat"{name=F}] & \mathsf{abs}~\mathcal{L}_\infty\text{-}\mathsf{alg}^{\mathsf{comp}} \arrow[l, shift left=.75ex, "\widehat{\mathrm{B}}_{\iota}^\flat"{name=U}]
            \arrow[phantom, from=F, to=U, , "\dashv" rotate=-90]
\end{tikzcd}
\]
between $\mathcal{CC}_\infty$-coalgebras and complete absolute $\mathcal{L}_\infty$-algebras, which we denote by $\widehat{\Omega}_{\iota}^\flat \dashv \widehat{\mathrm{B}}_{\iota}^\flat$ in order to avoid confusion with the previous one. \textit{Mutatis mutandis}, the same arguments apply to this case: the category of $\mathcal{CC}_\infty$-coalgebras admits a left-transferred from dg modules, and this model structure can be transferred to the category of complete absolute $\mathcal{L}_\infty$-algebras via the complete bar-cobar adjunction. Then, the complete bar-cobar adjunction becomes an Quillen equivalence. Let us compare their homotopy theory with that of curved absolute $\mathcal{L}_\infty$-algebras and that of classical $\mathcal{L}_\infty$-algebras.

\begin{Proposition}\label{Prop: fully faithful inclusion of non-curved into curved}
The adjunction 
\[
\begin{tikzcd}[column sep=7pc,row sep=3pc]
\mathsf{curv}~\mathsf{abs}~\mathcal{L}_\infty\text{-}\mathsf{alg}^{\mathsf{comp}} \arrow[r, shift left=1.1ex, "(-)_*"{name=F}]      
&\mathsf{abs}~\mathcal{L}_\infty\text{-}\mathsf{alg}^{\mathsf{comp}}~, \arrow[l, shift left=.75ex, "\mathrm{U}"{name=U}]
\arrow[phantom, from=F, to=U, , "\dashv" rotate=-90]
\end{tikzcd}
\]
between complete curved absolute $\mathcal{L}_\infty$-algebras and complete absolute $\mathcal{L}_\infty$-algebras is Quillen adjunction. Furthermore, the functor $\mathrm{U}$ is homotopically fully faithful. 
\end{Proposition}

\begin{proof}
The adjunction is a Quillen adjunction since it fits in the following commutative square of Quillen adjunctions
\[
\begin{tikzcd}[column sep=5pc,row sep=5pc]
\mathcal{CC}_\infty\text{-}\mathsf{coalg} \arrow[r,"\widehat{\Omega}_\iota^\flat"{name=B},shift left=1.1ex] \arrow[d,"(-) \oplus \mathscr{C}(u\mathcal{CC}_\infty)(0)"{name=SD},shift left=1.1ex ]
&\mathsf{abs}~\mathcal{L}_\infty\text{-}\mathsf{alg}^{\mathsf{comp}} \arrow[d,"\mathrm{U}"{name=LDC},shift left=1.1ex ] \arrow[l,"\widehat{\mathrm{B}}_\iota^\flat"{name=C},,shift left=1.1ex]  \\
u\mathcal{CC}_\infty\text{-}\mathsf{coalg} \arrow[r,"\widehat{\Omega}_\iota "{name=CC},shift left=1.1ex]  \arrow[u,"\mathrm{Res}"{name=LD},shift left=1.1ex ]
&\mathsf{curv}~\mathsf{abs}~\mathcal{L}_\infty\text{-}\mathsf{alg}^{\mathsf{comp}}~, \arrow[l,"\widehat{\mathrm{B}}_\iota"{name=CB},shift left=1.1ex] \arrow[u,"(-)_*"{name=TD},shift left=1.1ex] \arrow[phantom, from=SD, to=LD, , "\dashv" rotate=0] \arrow[phantom, from=C, to=B, , "\dashv" rotate=-90]\arrow[phantom, from=TD, to=LDC, , "\dashv" rotate=0] \arrow[phantom, from=CC, to=CB, , "\dashv" rotate=-90]
\end{tikzcd}
\] 
where $\mathscr{C}(u\mathcal{CC}_\infty)(0)$ is the terminal $u\mathcal{CC}_\infty$-coalgebra. Let us show that the functor $\mathrm{U}$ is homotopically fully faithful. It amounts to showing the following: given an absolute $\mathcal{L}_\infty$-algebra $\mathfrak{g}$, the canonical map 
\[
\widehat{\mathrm{B}}_{\iota}^\flat(\mathfrak{g}) \oplus \mathscr{C}(u\mathcal{CC}_\infty)(0) \longrightarrow \widehat{\mathrm{B}}_{\iota}(\mathfrak{g})~,
\]
is a quasi-isomorphism of $\mathcal{CC}_\infty$-coalgebras. For any dg module $V$, there is an isomorphism of $\mathcal{CC}_\infty$-coalgebras
\[
\mathscr{C}(u\mathcal{CC}_\infty)(V) \cong \mathscr{C}(\mathcal{CC}_\infty)(V) \oplus \mathscr{C}(u\mathcal{CC}_\infty)(0)~,
\]
since $\Omega \mathrm{B} \mathcal{C}om$ is the maximal reduces sub-operad of $\Omega\mathrm{B}^{\mathrm{s.a}}\ucom$, see \cite[Proofs of Lemmas 4.5 and 4.7]{operadswithout}. When $\mathfrak{g}$ is an absolute $\mathcal{L}_\infty$-algebra, this isomorphism commutes with the differential induced on the complete bar construction by the absolute $\mathcal{L}_\infty$-algebra structure on $\mathfrak{g}$.
\end{proof}

\begin{Remark}
The terminal $u\mathcal{CC}_\infty$-coalgebra $\mathscr{C}(u\mathcal{CC}_\infty)(0)$ is quasi-isomorphic to $\kk$, see \cite[Lemma 1.24]{mathez}. 
\end{Remark}

Absolute $\mathcal{L}_\infty$-algebra also admit a homotopical comparison with their classical counterparts.

\begin{Proposition}\label{prop: adjonction absolute non-courbe et non-absolute}
The adjunction  

\[
\begin{tikzcd}[column sep=7pc,row sep=3pc]
\mathcal{L}_\infty\text{-}\mathsf{alg} \arrow[r, shift left=1.1ex, "\mathrm{Abs}"{name=F}]      
&\mathsf{abs}~\mathcal{L}_\infty\text{-}\mathsf{alg}^{\mathsf{comp}}~, \arrow[l, shift left=.75ex, "\mathrm{Res}"{name=U}]
\arrow[phantom, from=F, to=U, , "\dashv" rotate=-90]
\end{tikzcd}
\]

between $\mathcal{L}_\infty$-algebras and complete absolute $\mathcal{L}_\infty$-algebras is a Quillen adjunction, where the category of $\mathcal{L}_\infty$-algebras is endowed with a model structure transferred from dg modules. In particular, any weak-equivalence of complete absolute $\mathcal{L}_\infty$-algebras is a quasi-isomorphism.
\end{Proposition}

\begin{proof}
The fact that any transferred weak-equivalence of complete absolute $\mathcal{L}_\infty$-algebras along the complete bar-cobar adjunction is, in particular, a quasi-isomorphism follows from analogue arguments to \cite[Proposition 4.6]{lucio2022contra}. The functor $\mathrm{Res}$ clearly preserves fibrations, and it also preserves weak-equivalences since any transferred weak-equivalence is a quasi-isomorphism. 
\end{proof}

\begin{Remark}[Homotopical completeness]
The above adjunction is not a Quillen equivalence. However, the functor $\mathrm{Abs}$ is homotopically fully faithful precisely on \textit{homotopically complete} $\mathcal{L}_\infty$-algebras. See \cite{HessHarper} for the original reference, and also \cite[Section 4]{lucio2022contra}. 
\end{Remark}

\subsection{Tensor product of $u\mathcal{CC}_\infty$-coalgebras}
The tensor product of two $u\mathcal{CC}_\infty$-coalgebras can naturally be endowed with a $u\mathcal{CC}_\infty$-coalgebra structure. This gives a closed monoidal structure on the category of $u\mathcal{CC}_\infty$-coalgebras which is compatible with its model category structure.

\begin{Proposition}\label{prop: existence of diagonal}
The dg operad $\Omega\mathrm{B}^{\mathrm{s.a}}u\mathcal{C}om$ is a dg Hopf operad, meaning there exists a diagonal morphism of operads
\[
\Delta_{\Omega\mathrm{B}\mathcal{C}om}: \Omega\mathrm{B}u\mathcal{C}om \longrightarrow \Omega\mathrm{B}u\mathcal{C}om \otimes_{H} \Omega\mathrm{B}u\mathcal{C}om
\]
which is coassociative.
\end{Proposition}

\begin{proof}
The operad $u\mathcal{C}om$ can be obtained as the cellular chains on the operad $\mathrm{uCom}$ in the category of sets given by $\{*\}$ in non-negative arities. Any set-theoretical operad is a cocommutative Hopf operad on the nose. By Theorem \ref{thm: Boardman-Vogt and cellular chains}, one has an isomorphism of dg operads

\[
C_*^c(\mathrm{W}(\mathrm{uCom})) \cong \Omega\mathrm{B}^{\mathrm{s.a}}\left(C_*^c(\mathrm{uCom})\right)~,
\]
\vspace{0.1pc}

where $\mathrm{W}$ denotes the Boardmann-Vogt construction of $\mathrm{uCom}$ seen as a discrete topological operad. Given the choice of a cellular approximation of the diagonal on the topological interval $\mathrm{I}$, one has a coassociative Hopf operad structure on $\mathrm{W}(\mathrm{uCom})$ which might not be cocommutative. The cellular chain functor $C_*^c(-)$ is both lax and colax (the Eilenberg-Ziber and Alexander-Whitney maps) in a compatible way by \cite[Chapter 5.3]{AguiarMahajan}. This is also in \cite[Proposition 26]{grignou2022mappingII}, and therefore by \cite[Proposition 25]{grignou2022mappingII}, the cellular chains functor sends Hopf topological operads to dg Hopf operads.
\end{proof}

\begin{Corollary}\label{cor: tensor product of uCC coalgebras}
The category $u\mathcal{CC}_\infty$-coalgebras can be endowed with a monoidal structure given by the tensor product of the underlying dg modules. 
\end{Corollary}

\begin{proof}
Consider two $u\mathcal{CC}_\infty$-coalgebras $(C_1,\Delta_{C_1},d_{C_1})$ and $(C_2,\Delta_{C_2},d_{C_2})$. These $u\mathcal{CC}_\infty$-coalgebra structures amount to two morphisms of dg operads $f_1: \Omega\mathrm{B}\mathcal{C}om \longrightarrow \mathrm{Coend}(C_1)$ and $f_2: \Omega\mathrm{B}\mathcal{C}om \longrightarrow \mathrm{Coend}(C_2)$. Using the Hopf structure on gets 

\[
\begin{tikzcd}[column sep=2pc,row sep=0pc]
\Omega\mathrm{B}\mathcal{C}om \arrow[r,"\Delta_{\Omega\mathrm{B}\mathcal{C}om}"]
&\Omega\mathrm{B}\mathcal{C}om \otimes_{H} \Omega\mathrm{B}\mathcal{C}om \arrow[r, "f_1 ~ \otimes ~ f_2"]
&\mathrm{Coend}(C_1) \otimes_{H} \mathrm{Coend}(C_2) \arrow[r]
&\mathrm{Coend}(C_1 \otimes C_2)~,
\end{tikzcd}
\]
\vspace{0.1pc}

which endows the dg module $(C_1 \otimes C_2, d_{C_1} \otimes d_{C_2})$ with a $u\mathcal{CC}_\infty$-coalgebra structure. The counitality and coassociativity of $\Delta_{\Omega\mathrm{B}\mathcal{C}om}$ ensure that the category of $u\mathcal{CC}_\infty$-coalgebra together with the tensor product forms a monoidal category.
\end{proof}

\begin{Remark}
The dg Hopf operad structure on $\Omega\mathrm{B}^{\mathrm{s.a}}u\mathcal{C}om$ constructed here extends the dg Hopf structure of $u\mathcal{C}om$, therefore the tensor product of two counital cocommutative coalgebras seen as $u\mathcal{CC}_\infty$-coalgebras coincides with the usual structure on the tensor product of two counital cocommutative coalgebras.
\end{Remark}

In fact, the tensor product of two coalgebras over a dg Hopf operad always forms a \textit{biclosed} monoidal category. We specify the general constructions of \cite{grignou2022mappingII} and of the forthcoming paper \cite{grignou2022mapping} in the particular cases of interest for us. The constructions of \textit{loc.cit} hold for a wider variety of cases.

\begin{Proposition}[{\cite[Proposition 16]{grignou2022mappingII}}]
The category of dg $u\mathcal{CC}_\infty$-coalgebras is a biclosed monoidal category, meaning that there exists a left (resp. right) internal hom bifunctor 

\[
\{-,-\}_{L/R}: \left(u\mathcal{CC}_\infty\text{-}\mathsf{coalg}\right)^{\mathsf{op}} \times u\mathcal{CC}_\infty\text{-}\mathsf{coalg} \longrightarrow u\mathcal{CC}_\infty\text{-}\mathsf{coalg}
\]
\vspace{0.1pc}

and, for any triple of $u\mathcal{CC}_\infty$-coalgebras $C,D,E$, there exists isomorphisms

\[
\mathrm{Hom}_{u\mathcal{CC}_\infty\text{-}\mathsf{coalg}}(C \otimes D, E) \cong \mathrm{Hom}_{u\mathcal{CC}_\infty\text{-}\mathsf{coalg}}(C, \{D,E\}_L)~,
\] 
\[
\mathrm{Hom}_{u\mathcal{CC}_\infty\text{-}\mathsf{coalg}}(C \otimes D, E) \cong \mathrm{Hom}_{u\mathcal{CC}_\infty\text{-}\mathsf{coalg}}(D, \{C,E\}_R)~,
\]
\vspace{0.1pc}

which are natural in $C,D$ and $E$.
\end{Proposition}

\begin{Remark}
Let $C_1$ and $C_2$ be two $u\mathcal{CC}_\infty$-coalgebras. Both of the internal homs can be constructed as equalizers
\[
\begin{tikzcd}[column sep=4pc,row sep=4pc]
\mathrm{Eq}\Bigg(\mathscr{C}(u\mathcal{CC}_\infty)(\mathrm{hom}(C_1,C_2)) \arrow[r,"(\Delta_{C_2})_*",shift right=1.1ex,swap]  \arrow[r,"\varrho"{name=SD},shift left=1.1ex ]
&\mathscr{C}(u\mathcal{CC}_\infty)\left(\mathrm{hom}(C_1,\widehat{\mathscr{S}}^c(u\mathcal{CC}_\infty)(C_2) \right)\Bigg)~,
\end{tikzcd}
\]
where $\mathrm{hom}(C_1,C_2)$ denotes the dg module of graded morphisms, where $\varrho$ is a map constructed using the comonad structure of $\mathscr{C}(u\mathcal{CC}_\infty)$ and the Hopf structure of the dg operad $u\mathcal{CC}_\infty$. This kind of construction works for coalgebras over any Hopf dg operad. 
\end{Remark}

The closed monoidal structure on $u\mathcal{CC}_\infty$-coalgebras is compatible with the model structure.

\begin{Proposition}\label{prop: monoidal model category}
The category of $u\mathcal{CC}_\infty$-coalgebras together with their tensor product forms a monoidal model category.
\end{Proposition}

\begin{proof}
This follows directly from the fact that the tensor product is computed in the underlying category of dg modules, and that (trivial cofibrations) of coalgebras are exactly (trivial) cofibrations of dg modules. Therefore it satisfies the axioms of a monoidal model category as stated in \cite{HoveyModel}.
\end{proof}

\begin{Remark}[Symmetric monoidal structures "up to homotopy"]
There are two choices for the Hopf structure in Proposition \ref{prop: existence of diagonal}, which correspond to the two possible cellular approximations of the diagonal of the interval. Algebraically, they correspond to the two coassociative coalgebra structures on $C_*^c(I)$, which are the one where $[0]$ is the counit and the one where $[1]$ is the counit. Each choice gives a different $u\mathcal{CC}_\infty$-coalgebra structure on the tensor product and for each choice there are left and right internal hom bifunctors; both choices are compatible with the underlying model structure. Let us denote $\otimes_0$ and $\otimes_1$ the two monoidal structures given by each of the possible choices. For any $u\mathcal{CC}_\infty$-coalgebras $A$ and $B$, it can be checked that there exists a natural isomorphism
\[
A \otimes_0 B \cong B \otimes_1 A~,
\]
which at the level of chain complexes is simply the natural isomorphism $A \otimes B \cong B \otimes A$. However, neither of the monoidal structures are \textit{symmetric}, since the Hopf operad structures on $\Omega\mathrm{B}^{\mathrm{s.a}}u\mathcal{C}om$ are only coassociative. 

\medskip

Nonetheless, the two cellular approximations of the diagonal of the interval are homotopic. Thus, the two choices for the monoidal structures are naturally weakly equivalent. In particular, it follows that $A \otimes_0 B \simeq B \otimes_0 A$ (and similarly for the product $\otimes_1$). Therefore left and right internal homs are also naturally weakly equivalent. Both monoidal structures on $u\mathcal{CC}_\infty$-coalgebras are "symmetric up to homotopy", although we do not know a good $1$-categorical definition for this notion. Whatever this might mean, it should be interpreted as inducing a symmetric monoidal structure on the underlying $\infty$-category obtained by localizing at quasi-isomorphisms. From now on, we will not specify which choice of cellular approximations of the diagonal of the interval is made, and we will not distinguish between left and right internal hom bifunctors, as they will only appear in a homotopical context. 
\end{Remark}
 
\begin{Definition}[Convolution curved absolute $\mathcal{L}_\infty$-algebra]
Let $C$ be a $u\mathcal{CC}_\infty$-coalgebra and let $\mathfrak{g}$ be a curved absolute $\mathcal{L}_\infty$-algebra. The pdg module of graded morphisms $(\mathrm{hom}(C,\mathfrak{g}),\partial)$ is endowed the following curved absolute $\mathcal{L}_\infty$-algebra structure.

\medskip

The structural map $\gamma_{\mathrm{hom}(C,\mathfrak{g})}$ is given by the following composition
\[
\begin{tikzcd}
\widehat{\mathscr{S}}^c(\mathrm{B}^{\mathsf{s.a}}u\mathcal{C}om)(\mathrm{hom}(C,\mathfrak{g})) \arrow[d, "\mathrm{coev}_C"]     \\
\mathrm{hom}\left(C, \widehat{\mathscr{S}}^c(\mathrm{B}^{\mathsf{s.a}}u\mathcal{C}om)(\mathrm{hom}(C,\mathfrak{g})) \otimes C \right) \arrow[d, "(\Delta_C)_*"] \\
\mathrm{hom} \left(C, \widehat{\mathscr{S}}^c(\mathrm{B}^{\mathsf{s.a}}u\mathcal{C}om)(\mathrm{hom}(C,\mathfrak{g})) \otimes \widehat{\mathscr{S}}^c(\Omega\mathrm{B}^{\mathsf{s.a}}u\mathcal{C}om)(C) \right) \arrow[d, "\xi"] \\
\mathrm{hom}\left(C, \widehat{\mathscr{S}}^c(\mathrm{B}^{\mathsf{s.a}}u\mathcal{C}om \otimes \Omega\mathrm{B}^{\mathsf{s.a}}u\mathcal{C}om)(\mathrm{hom}(C,\mathfrak{g}) \otimes C)\right) \arrow[d,"\widehat{\mathscr{S}}^c(\delta_{\mathrm{B}^{\mathrm{s.a}}u\mathcal{C}om})"] \\
\mathrm{hom}\left(C, \widehat{\mathscr{S}}^c(\mathrm{B}^{\mathsf{s.a}}u\mathcal{C}om)(\mathrm{hom}(C,\mathfrak{g}) \otimes C)\right) \arrow[d," \mathrm{ev}_C "] \\
\mathrm{hom}\left(C, \widehat{\mathscr{S}}^c(\mathrm{B}^{\mathsf{s.a}}u\mathcal{C}om)(\mathfrak{g})\right) \arrow[d,"(\gamma_\mathfrak{g})_* "] \\
\mathrm{hom}(C, \mathfrak{g})~, \\
\end{tikzcd}
\]
where $\mathrm{coev}_C$ and $\mathrm{ev}_C$ are respectively the unit and the counit of the tensor-hom adjunction, and where $\xi$ is the following natural inclusion 
\[
\begin{tikzcd}
\displaystyle \left(\prod_{n \geq 0} \mathrm{Hom}_{\mathbb{S}_n}(M(n), V^{\otimes n})\right) \otimes \left(\prod_{n \geq 0} \mathrm{Hom}_{\mathbb{S}_n}(N(n), W^{\otimes n})\right) \arrow[d,rightarrowtail] \\
\displaystyle \prod_{n \geq 0} \mathrm{Hom}_{\mathbb{S}_n}(M(n) \otimes N(n), (V \otimes W)^{\otimes n})~.
\end{tikzcd}
\]
Finally

\[
\delta_{\mathrm{B}^{\mathsf{s.a}}u\mathcal{C}om}: \mathrm{B}^{\mathsf{s.a}}u\mathcal{C}om \longrightarrow \mathrm{B}^{\mathsf{s.a}}u\mathcal{C}om \otimes \Omega\mathrm{B}^{\mathsf{s.a}}u\mathcal{C}om
\]
\vspace{0.1pc}

is a restriction of the diagonal $\Delta_{\Omega\mathrm{B}\mathcal{C}om}$.
\end{Definition}

\begin{Remark}
Let $C$ be a counital cocommutative coalgebra viewed as a $u\mathcal{CC}_\infty$-coalgebra, and let $\mathfrak{g}$ be a curved absolute $\mathcal{L}_\infty$-algebra. The convolution structure on $\mathrm{hom}(C,\mathfrak{g})$ greatly simplifies. In particular, the elementary operations are given by 
\[
l_n(f_1,\cdots, f_n) = l_n^{\mathfrak{g}} \circ (f_1,\cdots,f_n) \circ \Delta_n \quad \text{for} \quad n \geq 2, 
\]
where $l_n^{\mathfrak{g}}$ is the elementary operation of $\mathfrak{g}$, $\Delta_n$ is the $n-1$-iteration of the coproduct of $C$ and $f_1,\cdots,f_n$ are graded morphisms in $\mathrm{hom}(C,\mathfrak{g})$. The curvature is given by $l_0^{\mathfrak{g}} \circ \epsilon$, where $l_0^{\mathfrak{g}}$ is the curvature of $\mathfrak{g}$ and $\epsilon$ the counit of $C$. Infinite sums of these operations are well-defined since $\mathfrak{g}$ is itself an absolute algebra. 
\end{Remark}

Convolution curved absolute $\mathcal{L}_\infty$-algebras and the internal hom-set of $u\mathcal{CC}_\infty$-coalgebras are compatible in the following sense.

\begin{theorem}[{\cite{grignou2022mapping}}]
Let $C$ be a $u\mathcal{CC}_\infty$-coalgebra and let $\mathfrak{g}$ be a curved absolute $\mathcal{L}_\infty$-algebra. There is an isomorphism of $u\mathcal{CC}_\infty$-coalgebras

\[
\left\{C, \widehat{\mathrm{B}}_\iota(\mathfrak{g}) \right\} \cong \widehat{\mathrm{B}}_\iota \left(\mathrm{hom}(C,\mathfrak{g}) \right)~,
\]
\vspace{0.1pc}

where $\mathrm{hom}(C,\mathfrak{g})$ denotes the convolution curved absolute $\mathcal{L}_\infty$-algebra of $C$ and $\mathfrak{g}$.
\end{theorem}

\section{Higher absolute Lie theory}
In this section, we follow an analogue approach to \cite{robertnicoud2020higher} in order to integrate curved absolute $\mathcal{L}_\infty$-algebras. The integration functor and its left adjoint form a Quillen adjunction between curved absolute $\mathcal{L}_\infty$-algebras and simplicial sets. The rest of this section is devoted to the study of the main properties of this adjunction.

\subsection{Dupont's contraction}\label{subsection: dupont contraction}
In \cite{DUPONT}, J-L Dupont proved that there is a homotopy contraction between the simplicial unital commutative algebra of polynomial differential forms on the geometrical simplex and the simplicial sub-module of Whitney forms. Whitney forms on the simplices are isomorphic to the cellular cochains on the simplices and they are finite dimensional. This allows us to use the homotopy transfer theorem, and obtain a simplicial $u\mathcal{CC}_\infty$-coalgebra structure on the cellular chains of simplices. This coalgebra is not conilpotent, hence we encode it with an operad. This enables us to use the complete bar-cobar construction in order to construct a commuting triangle of Quillen adjunctions. This approach allows us to extend and to refine automatically many of the standard results in the theory of integration. 

\medskip

The simplicial dg unital commutative algebra of piece-wise polynomial differential forms on the standard simplex $\Omega_\bullet$ is given by
\[
\Omega_n \coloneqq \frac{\kk[t_0,\cdots,t_n, dt_0,\cdots, dt_n]}{(t_0 + \cdots + t_n -1, dt_0 + \cdots + dt_n)}~,
\]
with the obvious simplicial structure.

\begin{theorem}[Dupont's contraction]\label{dupontcontraction}
There exists a simplicial contraction
\[
\begin{tikzcd}[column sep=5pc,row sep=3pc]
\Omega_{\bullet} \arrow[r, shift left=1.1ex, "p_{\bullet}"{name=F}] \arrow[loop left]{l}{h_{\bullet}}
& C_c^*(\Delta^{\bullet})~, \arrow[l, shift left=.75ex, "i_{\bullet}"{name=U}]
\end{tikzcd}
\]
where $C_c^*(\Delta^\bullet)$ denotes the simplicial dg module given by the cellular cochains of the standard simplex. 
\end{theorem}

A basis of $C^*_c(\Delta^n)$ is given by $\{\omega_I\}$ for $I = \{ i_0,\cdots,i_k \} \subset [n]$, where $|\omega_I| = -k$. Since we are working in the homological convention, both algebras are concentrated in degrees $\leq 0$, with $dt_i$ being of degree $-1$. 

\medskip

The data of a simplicial contraction amounts to the data of two morphisms of simplicial cochain complexes
\[
i_{\bullet}: \Omega_{\bullet} \longrightarrow C^*_c(\Delta^{\bullet})~, \quad \text{and} \quad p_{\bullet}: C^*_c(\Delta^{\bullet}) \longrightarrow \Omega_{\bullet}~;
\]
and a degree $1$ linear map $h_{\bullet}: \Omega_{\bullet}\longrightarrow \Omega_{\bullet}$, satisfying the following conditions:
\[
p_n i_n = \mathrm{id}_{C^*_c(\Delta^n)}~, \quad i_n p_n - \mathrm{id}_{\Omega_n} = d_n h_n +h_n d_n~, \quad h_n^2 = 0~, \hspace{1pc} p_n h_n = 0~, \hspace{1pc} h_n i_n = 0~,
\]
where $d_n$ stands for the differential of $\Omega_n$. 

\begin{Remark}
For $n \geq 1$, the dg unital commutative algebra $\Omega_n$ is not \textit{naturally} augmented, as we have $1 = t_0 + \cdots + t_n$. In order to keep track of the full structure on $\Omega_n$, one has to see it as a dg $\ucom$-algebra and not as a dg $\mathcal{C}om$-algebra.
\end{Remark}

\begin{lemma}\label{lemma: simplicial CC infinity}
There is a simplicial $u\mathcal{CC}_\infty$-algebra structure on $C^*_c(\Delta^{\bullet})$. 
\end{lemma}

\begin{proof}
The dg operad $\Omega\mathrm{B}^{\mathrm{s.a}}\ucom$ is a cofibrant resolution of $\ucom$. Therefore we can apply the homotopy transfer theorem given in \cite[Theorem 6.5.5]{HirshMilles12} using the Dupont contraction of Theorem \ref{dupontcontraction}. Since this contraction is compatible with the simplicial structure, we obtain a simplicial dg $\Omega\mathrm{B}^{\mathrm{s.a}}\ucom$-algebra.
\end{proof}

Let $C^c_*(\Delta^n)$ denote the cellular chains on the $n$-simplex. 

\begin{lemma}\label{lemma: cosimplicial CC infity}
There is a cosimplicial $u\mathcal{CC}_\infty$-coalgebra structure on $C^c_*(\Delta^\bullet)$.
\end{lemma}

\begin{proof}
The linear dual of a degree-wise finite dimensional algebra over an dg operad is naturally a dg coalgebra over the same operad. Since $C^*_c(\Delta^n)$ is finite dimensional degree-wise for all $n$, its linear dual $C^c_*(\Delta^n)$ carries a canonical a $u\mathcal{CC}_\infty$-coalgebra structure. Applying a contravariant functor to a simplicial object gives a cosimplicial object.
\end{proof}

\begin{Remark}
Notice that $C^c_*(\Delta^\bullet)$ is a non-conilpotent $u\mathcal{CC}_\infty$-coalgebra, therefore it cannot be encoded by a cooperad, like stated in \cite[Section 2]{robertnicoud2020higher}.
\end{Remark}

We use this cosimplicial $u\mathcal{CC}_\infty$-coalgebra to construct an adjunction between simplicial sets and $u\mathcal{CC}_\infty$-coalgebras, using the following seminal result.

\begin{theorem}[{\cite{Kan58}}]\label{thm: Kan seminal result}
Let $\mathsf{C}$ be a locally small cocomplete category. The data of an adjunction
\[
\begin{tikzcd}[column sep=7pc,row sep=3pc]
            \mathsf{sSet} \arrow[r, shift left=1.1ex, "L"{name=F}] &\mathscr{C} \arrow[l, shift left=.75ex, "R"{name=U}]
            \arrow[phantom, from=F, to=U, , "\dashv" rotate=-90]
\end{tikzcd}
\]
is equivalent to the data of a cosimplicial object $F: \Delta \longrightarrow \mathscr{C}$ in $\mathscr{C}$.
\end{theorem}

\begin{proof}
Given an adjunction $L \dashv R$, one can pullback the left adjoint $L$ along the Yoneda embedding $\mathrm{Yo}$
\[
\begin{tikzcd}
	\Delta \arrow[r,"\mathrm{Yo}"] 
	&\mathsf{sSet} \arrow[r,"L"]
	&\mathsf{sSet}~,
\end{tikzcd}
\]
and obtain a cosimplicial object in $\mathscr{C}$. Given a cosimplicial object $F: \Delta \longrightarrow \mathscr{C}$ in $\mathscr{C}$, one can consider its left Kan extension $\mathrm{Lan}_{\mathrm{Yo}}(F)$ since $\mathscr{C}$ is cocomplete. Its right adjoint is given by 
\[
R(-)_\bullet \coloneqq \mathrm{Hom}_{\mathscr{C}}(F(\Delta^{\bullet}),-)~.
\]
\end{proof}

\begin{Proposition}
There is an adjunction 
\[
\begin{tikzcd}[column sep=7pc,row sep=3pc]
            \mathsf{sSet} \arrow[r, shift left=1.1ex, "C^c_*(-)"{name=F}] 
            &u\mathcal{CC}_\infty\text{-}\mathsf{coalg}~, \arrow[l, shift left=.75ex, "\overline{\mathcal{R}}"{name=U}]
            \arrow[phantom, from=F, to=U, , "\dashv" rotate=-90]
\end{tikzcd}
\]
where $C^c_*(-)$ denotes the cellular chain functor endowed with a canonical $u\mathcal{CC}_\infty$-coalgebra structure. 
\end{Proposition}

\begin{proof}
This is a straightforward application of Theorem \ref{thm: Kan seminal result} to the cosimplicial $u\mathcal{CC}_\infty$-coalgebra constructed in Lemma \ref{lemma: cosimplicial CC infity}. It is immediate to check that the left adjoint is given by the cellular chain functor $C^c_*(-)$ as a dg module, endowed with a $u\mathcal{CC}_\infty$-coalgebra structure. This structure is given by the following colimit

\[
C_*^c(X) \cong \colim_{\mathrm{E}(X)} C_*^c(\Delta^\bullet)~.
\]

on the category of elements $E(X)$ of $X$ in the category of $u\mathcal{CC}_\infty$-coalgebras.
\end{proof}

The right adjoint is therefore given, for a $u\mathcal{CC}_\infty$-coalgebra $C$, by the functor
\[
\overline{\mathcal{R}}(C)_\bullet \coloneqq \mathrm{Hom}_{u\mathcal{CC}_\infty\text{-}\mathsf{cog}}(C^c_*(\Delta^\bullet),C)~.
\]

By pushing forward along the complete cobar construction $\widehat{\Omega}_\iota$ the cosimplicial $u\mathcal{CC}_\infty$-coalgebra $C^c_*(\Delta^\bullet)$, we also obtain a cosimplicial complete curved absolute $\mathcal{L}_\infty$-algebra.

\begin{Definition}[Maurer--Cartan cosimplicial algebra]
The \textit{Maurer--Cartan cosimplicial algebra} is given by the cosimplicial complete curved absolute $\mathcal{L}_\infty$-algebra
\[
\mathfrak{mc}^{\bullet} \coloneqq \widehat{\Omega}_{\iota}C^c_*(\Delta^{\bullet})~.
\]
\end{Definition}

\begin{Proposition}
There is an adjunction 
\[
\begin{tikzcd}[column sep=7pc,row sep=3pc]
            \mathsf{sSet} \arrow[r, shift left=1.1ex, "\mathcal{L}"{name=F}] 
            &\mathsf{curv}~\mathsf{abs}~\mathcal{L}_\infty\text{-}\mathsf{alg}^{\mathsf{comp}}~, \arrow[l, shift left=.75ex, "\mathcal{R}"{name=U}]
            \arrow[phantom, from=F, to=U, , "\dashv" rotate=-90]
\end{tikzcd}
\]
between the categories of simplicial sets and of curved absolute $\mathcal{L}_\infty$-algebras.
\end{Proposition}

\begin{proof}
This follows directly from Theorem \ref{thm: Kan seminal result}, applied to the cosimplicial object $\mathfrak{mc}^{\bullet}$.
\end{proof}

\begin{Remark}
The right adjoint is therefore given by, for a curved absolute $\mathcal{L}_\infty$-algebra $\mathfrak{g}$, by  
\[
\mathcal{R}(\mathfrak{g})_\bullet \coloneqq \mathrm{Hom}_{\mathsf{curv}~\mathsf{abs}~\mathcal{L}_\infty\text{-}\mathsf{alg}} \left(\widehat{\Omega}_\iota C^c_*(\Delta^\bullet), \mathfrak{g}\right)~.
\]
\end{Remark}

\begin{theorem}\label{thm: triangle of adjunctions}
The following triangle of adjunctions

\[
\begin{tikzcd}[column sep=5pc,row sep=2.5pc]
&\hspace{1pc}u\mathcal{CC}_\infty \textsf{-}\mathsf{coalg} \arrow[dd, shift left=1.1ex, "\widehat{\Omega}_{\iota}"{name=F}] \arrow[ld, shift left=.75ex, "\overline{\mathcal{R}}"{name=C}]\\
\mathsf{sSet}  \arrow[ru, shift left=1.5ex, "C^c_*(-)"{name=A}]  \arrow[rd, shift left=1ex, "\mathcal{L}"{name=B}] \arrow[phantom, from=A, to=C, , "\dashv" rotate=-70]
& \\
&\hspace{3pc}\mathsf{curv}~\mathsf{abs}~\mathcal{L}_\infty\textsf{-}\mathsf{alg}^{\mathsf{comp}}~, \arrow[uu, shift left=.75ex, "\widehat{\mathrm{B}}_{\iota}"{name=U}] \arrow[lu, shift left=.75ex, "\mathcal{R}"{name=D}] \arrow[phantom, from=B, to=D, , "\dashv" rotate=-110] \arrow[phantom, from=F, to=U, , "\dashv" rotate=-180]
\end{tikzcd}
\]

commutes. Furthermore, all the adjunctions are Quillen adjunctions when we consider the model category structures on $u\mathcal{CC}_\infty$-coalgebras and on complete curved absolute $\mathcal{L}_\infty$-algebras of Theorem \ref{thm: Equivalence the Quillen CC infini et absolues}; and the Kan-Quillen model category structure on simplicial sets.
\end{theorem}

\begin{proof}
Let $\mathfrak{g}$ be a curved absolute $\mathcal{L}_\infty$-algebra. We have that
\[
\mathcal{R}(\mathfrak{g})_\bullet = \mathrm{Hom}_{\mathsf{curv}~\mathsf{abs}~\mathcal{L}_\infty\text{-}\mathsf{alg}}\left(\widehat{\Omega}_\iota C^c_*(\Delta^\bullet), \mathfrak{g} \right) \cong \mathrm{Hom}_{u\mathcal{CC}_\infty\text{-}\mathsf{cog}}\left( C^c_*(\Delta^\bullet),\widehat{\mathrm{B}}_\iota \mathfrak{g} \right) \cong \overline{\mathcal{R}} \left(\widehat{\mathrm{B}}_\iota \mathfrak{g} \right)_\bullet~.
\]

The functor $C^c_*(-)$ sends monomorphisms of simplicial sets to degree-wise monomorphisms of $u\mathcal{CC}_\infty$-coalgebras and it sends weak homotopy equivalences to quasi-isomorphisms. Hence the triangle is made up of Quillen adjunctions.
\end{proof}

\subsection{The integration functor} The functor $\mathcal{R}$ in the above triangle is \textit{the integration functor} we were looking for. Before giving an explicit combinatorial description of it, let us first state some of its fundamental properties.

\begin{theorem}\label{thm: propriétés de l'intégration}\leavevmode

\begin{enumerate}
\item For any complete curved absolute $\mathcal{L}_\infty$-algebra $\mathfrak{g}$, the simplicial set $\mathcal{R}(\mathfrak{g})$ is a Kan complex.

\medskip

\item Let $f: \mathfrak{g} \twoheadrightarrow \mathfrak{h}$ be a degree-wise epimorphism of complete curved absolute $\mathcal{L}_\infty$-algebras. Then 
\[
\mathcal{R}(f): \mathcal{R}(\mathfrak{g}) \twoheadrightarrow \mathcal{R}(\mathfrak{h})
\]
is a fibrations of simplicial sets. 

\medskip

\item The functor $\mathcal{R}$ preserves weak-equivalences. In particular, it sends any graded quasi-isomorphism $f: \mathfrak{g} \qi \mathfrak{h}$ between complete curved absolute $\mathcal{L}_\infty$-algebras to a weak homotopy equivalence of simplicial sets.

\end{enumerate}
\end{theorem}

\begin{proof}
The functor $\mathcal{R}$ is a right Quillen functor. 
\end{proof}

\begin{Remark}
Let us contextualize the results of Theorem \ref{thm: propriétés de l'intégration}.

\begin{enumerate}
\item The first result is quintessential to the classical integration theory of dg Lie algebras or nilpotent $\mathcal{L}_\infty$-algebras as developed in \cite{Hinich01} and in \cite{Getzler09}.

\vspace{0.5pc}

\item The second result generalizes one of the main theorems of \cite{Getzler09}, see Theorem 5.8 in \textit{loc.cit}.

\vspace{0.5pc}

\item The third point, the \textit{homotopy invariance} of the integration functor, is the generalization of Goldman--Milson's invariance theorem proved in \cite{goldmanmillson} and its generalization to $\mathcal{L}_\infty$-algebras given in \cite{dolgushevrogers}. 
\end{enumerate}
\end{Remark}

\begin{Proposition}
The functor $\mathcal{R}$ commutes with homotopy limits.
\end{Proposition}

\begin{proof}
The functor $\mathcal{R}$ is a right Quillen functor. 
\end{proof}

\begin{Remark}
This means that it satisfies \textit{descent} in the sense of \cite{hinichdescent}. 
\end{Remark}

\begin{Proposition}
Let $\mathfrak{g}$ be a complete curved absolute $\mathcal{L}_\infty$-algebra. There is an isomorphism of simplicial sets
\[
\mathcal{R}(\mathfrak{g})_\bullet \cong \lim_{\omega} \mathcal{R}(\mathfrak{g}/\mathrm{W}_\omega \mathfrak{g})_\bullet~.
\]
\end{Proposition}

\begin{proof}
Since $\mathfrak{g}$ is complete, it can be written as 
\[
\mathfrak{g} \cong \lim_{\omega} \mathfrak{g}/\mathrm{W}_\omega \mathfrak{g}~,
\]
where the limit is taken in the category of curved absolute $\mathcal{L}_\infty$-algebras. Since $\mathcal{R}$ is right adjoint, it preserves all limits.
\end{proof}

\begin{Remark}[Comparison with Getzler integration functor]\label{Remark: comparison with Getzler functor}
One way wonder how to compare the integration functor here defined with the integration functor defined by E. Getzler for nilpotent $\mathcal{L}_\infty$-algebras in \cite{Getzler09}. We refer the reader to \cite[Chapter 3, Section 2.7]{mathez}. There we show that our functor is isomorphic to Getzler's original functor on the category of nilpotent $\mathcal{L}_\infty$-algebras (which are particular examples of curved absolute $\mathcal{L}_\infty$-algebras). We also show that it coincides with the integration functor of \cite{robertnicoud2020higher} under some technical hypothesis on the filtrations considered.
\end{Remark}

\textbf{Dold-Kan correspondence.} The adjunction 
\[
\begin{tikzcd}[column sep=5pc,row sep=2.5pc]
\mathsf{sSet}  \arrow[r, shift left=1ex, "\mathcal{L}"{name=B}] 
&\mathsf{curv}~\mathsf{abs}~\mathcal{L}_\infty\textsf{-}\mathsf{alg}^{\mathsf{comp}} \arrow[l, shift left=.75ex, "\mathcal{R}"{name=D}] \arrow[phantom, from=B, to=D, , "\dashv" rotate=-90] 
\end{tikzcd}
\]
is a \textit{generalization of the Dold-Kan correspondence}. Notice that a chain complex $(V,d_V)$ is a particular example of a complete curved absolute $\mathcal{L}_\infty$-algebra where the structural morphism 
\[
\gamma_V: \displaystyle \prod_{n \geq 0} \widehat{\Omega}^{\mathrm{s.c}}\ucom^*(n) ~ \widehat{\otimes}_{\mathbb{S}_n} ~ V^{\otimes n} \longrightarrow V
\]
is the \textit{zero morphism}. Any chain complex $(V,d_V)$ can be endowed with a canonical \textit{abelian} (trivial) curved absolute $\mathcal{L}_\infty$-algebra structure. 

\begin{Proposition}
Let $(V,d_V)$ be a chain complex endowed with the trivial complete curved absolute $\mathcal{L}_\infty$-algebra structure. Then
\[
\mathcal{R}(V) \cong \Gamma(V)~,
\]
where $\Gamma(-)$ is the Dold-Kan functor. 
\end{Proposition}

\begin{proof}
Recall that 
\[
\mathcal{R}(V)_\bullet \cong \mathrm{Hom}_{u\mathcal{CC}_\infty\text{-}\mathsf{cog}}(C^c(\Delta^\bullet),\widehat{\mathrm{B}}_\iota V )~.
\]
Since $\gamma_V = 0$, the complete bar construction $\widehat{\mathrm{B}}_\iota V $ is simply given by the cofree $u\mathcal{CC}_\infty$-coalgebra generated by the dg module $(V,d_V)$. Thus 

\[
\mathcal{R}(V)_\bullet \cong \mathrm{Hom}_{u\mathcal{CC}_\infty\text{-}\mathsf{cog}}(C^c_*(\Delta^\bullet),\mathscr{C}(u\mathcal{CC}_\infty)(V)) \cong \mathrm{Hom}_{\mathsf{dg}\text{-}\mathsf{mod}}(C^c_*(\Delta^\bullet),V) \cong \Gamma(V)~.
\]
\end{proof}

\begin{Remark}
Chain complexes concentrated in positive degrees with the trivial structure correspond via the functor $\mathcal{R}$ to simplicial abelian groups. One can think of general curved absolute $\mathcal{L}_\infty$-algebras concentrated in positive degrees as producing, via $\mathcal{R}$, a up to homotopy notion of a non-necessarily commutative group. This intuition can be made more precise using higher Baker--Campbell-Hausdorff products, see Subsection \ref{subsection: higher BCH}. 
\end{Remark}

\textbf{Derived adjunction.} The previous adjunction computes the \textit{derived functors} of the adjunction
\[
\begin{tikzcd}[column sep=5pc,row sep=2.5pc]
\mathsf{sSet}  \arrow[r, shift left=1.5ex, "C^c_{*}(-)"{name=A}]
&u\mathcal{CC}_\infty \textsf{-}\mathsf{coalg}~~. \arrow[l, shift left=.75ex, "\overline{\mathcal{R}}"{name=C}] \arrow[phantom, from=A, to=C, , "\dashv" rotate=-90]
\end{tikzcd}
\]
Let $C$ be a $u\mathcal{CC}_\infty$-coalgebra. It might not be fibrant, therefore one needs to take a fibrant resolution of $C$ in order to compute the right derived functor $\mathbb{R}\overline{\mathcal{R}}$ at $C$. The Theorem \ref{thm: Equivalence the Quillen CC infini et absolues} provides us with 
\[
\eta_C: C \qi \widehat{\mathrm{B}}_\iota \widehat{\Omega}_\iota C ~,
\]
a functorial fibrant resolution of $C$. Therefore we have that

\[
\mathbb{R}\overline{\mathcal{R}}(C)_\bullet = \mathrm{Hom}_{u\mathcal{CC}_\infty\text{-}\mathsf{cog}}(C^c(\Delta^\bullet),\widehat{\mathrm{B}}_\iota \widehat{\Omega}_\iota  C ) \cong \mathrm{Hom}_{\mathsf{curv}~\mathsf{abs}~\mathcal{L}_\infty\textsf{-}\mathsf{alg}}(\widehat{\Omega}_\iota C^c(\Delta^\bullet) , \widehat{\Omega}_\iota C ) \cong \mathcal{R}(\widehat{\Omega}_\iota  C )~.
\]

\begin{Remark}[$\infty$-morphisms]
Let $C$ and $D$ be two $u\mathcal{CC}_\infty$-coalgebras. The hom-set 

\[
\mathrm{Hom}_{\mathsf{curv}~\mathsf{abs}~\mathcal{L}_\infty\textsf{-}\mathsf{alg}}(\widehat{\Omega}_\iota C , \widehat{\Omega}_\iota  D ) 
\]

is the set of $\infty$-morphisms of $u\mathcal{CC}_\infty$-coalgebras between $C$ and $D$. See \cite[Section 12]{grignoulejay18}.
\end{Remark}

\subsection{Explicit version of the integration functor} Let us now give a combinatorial description of the integration functor. Recall that the integration functor $\mathcal{R}$ is given by 
\[
\mathcal{R}(-)_\bullet \coloneqq \mathrm{Hom}_{\mathsf{curv}~\mathsf{abs}~\mathcal{L}_\infty\textsf{-}\mathsf{alg}}(\widehat{\Omega}_\iota(C^c_* \Delta^\bullet),-) \cong \mathrm{Hom}_{\mathsf{curv}~\mathsf{abs}~\mathcal{L}_\infty\textsf{-}\mathsf{alg}}(\mathfrak{mc}^\bullet,-)~.
\]
Thus it is primordial to describe the complete curved absolute $\mathcal{L}_\infty$-algebra structure on the cosimplicial Maurer--Cartan algebra $\mathfrak{mc}^\bullet$ which represents this functor. Since it is given by the complete cobar construction of the cosimplicial $u\mathcal{CC}_\infty$-coalgebra $C^c_*(\Delta^\bullet)$, let us describe this structure.

\begin{Notation}
Consider $\Omega_n$ with its dg $u\mathcal{C}om$-algebra structure. Let
\[
\mu_k: \ucom(k) \otimes_{\mathbb{S}_k} \Omega_n^{\otimes k} \longrightarrow \Omega_n
\]
denote its dg $\ucom$-algebra structural morphisms for $k \geq 1$, and $u: \kk \longrightarrow \Omega_n$ its unit. 
\end{Notation}

\begin{Proposition}\label{prop: formulas for the HTT}
The transferred $u\mathcal{CC}_\infty$-algebra structure on $C_c^*(\Delta^n)$ via the Dupont contraction of Theorem \ref{dupontcontraction} is determined by morphisms
\[
\Big\{ \mu_\tau: (C_c^*(\Delta^n))^{\otimes m} \longrightarrow C_c^*(\Delta^n) \Big\}
\]
of degree $-\omega + 1$, where $\tau$ is a rooted tree in $\mathrm{RT}_m^\omega$ of arity $m$ and weight $\omega$ with no corks. 

\medskip

For $\tau$ in $\mathrm{RT}_m$, the operation $\mu_\tau$ is given by labeling all the vertices of $\tau$ with an operation $\mu_k$ of $\Omega_n$, where $k$ is the number of inputs of the vertex, by labeling all internal edges with the homotopy $h_n$, all the leaves with the map $i_n$ and the root with the map $p_n$, and composing the labeling operations along the rooted tree $\tau$. Pictorially, it is given by 

\begin{center}
\includegraphics[width=115mm,scale=1.15]{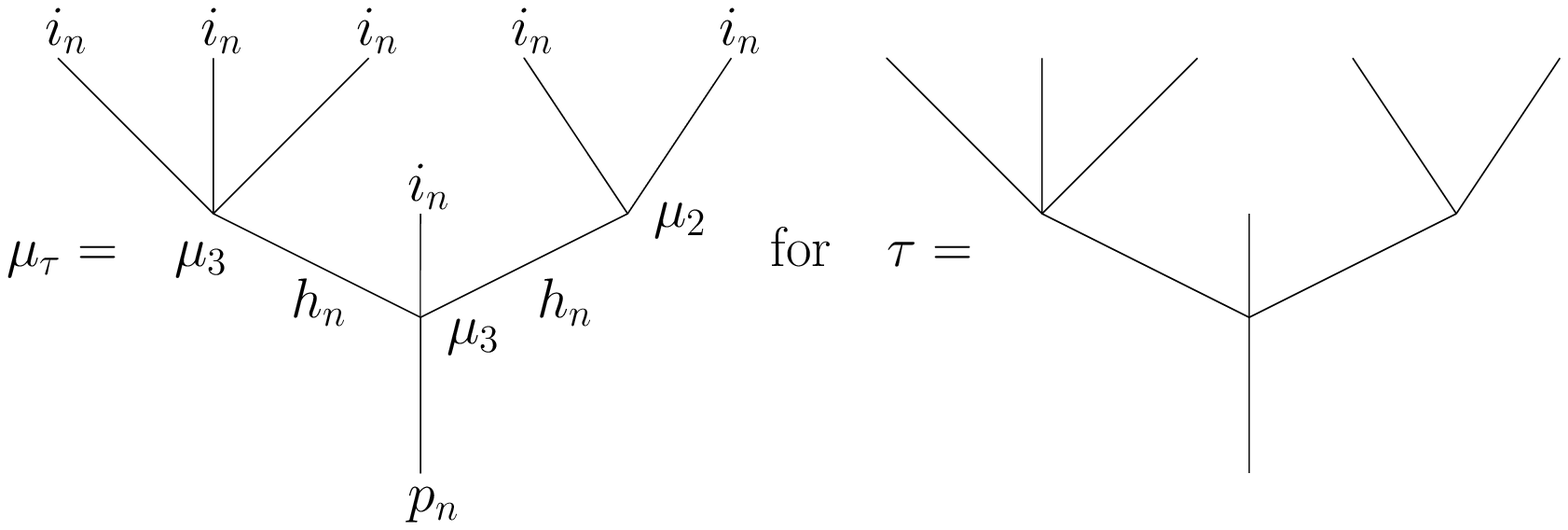}.
\end{center}

The only corked rooted tree acting non-trivially is the single cork morphism given by $p_n \circ u: \kk \longrightarrow C_c^*(\Delta^n) $.
\end{Proposition}

\begin{proof}
The standard transfer formula for the homotopy transfer theorem given in \cite[Theorem 6.5.5]{HirshMilles12} involves operations $\mu_\tau$ where $\tau$ runs over all corked rooted trees in $\mathrm{CRT}$. Two different types of corks can appear: the corks given by $u$ and the corks given by $u \circ h_n$. Pictorially, $\tau$ can be as follows

\begin{center}
\includegraphics[width=75mm,scale=1]{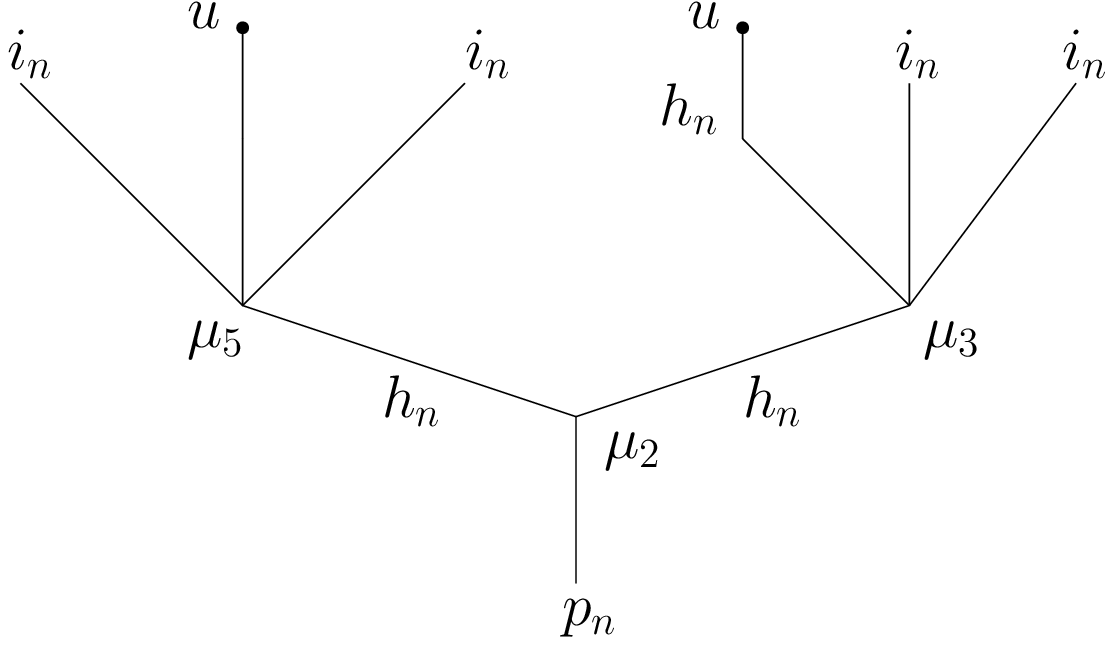},
\end{center}

where the first type of corks is represented on the left branch of the tree and the second type is represented on the right branch of the tree. Since we started with a strict $\ucom$-algebra structure on $\Omega_n$, the composition $\mu_k \circ_i u$ simplifies into $\mu_{k-1}$, hence the first type of corks disappears. Recall that $\Omega_n$ is concentrated in non-positive degrees and that the unit of $\Omega_n$ lies in degree $0$. Hence $h_n \circ u =0$ as $h_n$ raises the degree by one and the second type of corks disappears as well. We are left with operations $\mu_\tau$ where $\tau$ runs over all rooted trees without corks $\mathrm{RT}_m^\omega~,$ with a single operation involving the unit $u$ given by $p_n \circ u$.
\end{proof}

\begin{Remark}
We recover exactly the same structure on $C^*_c(\Delta^n)$ as the one given in \cite[Proposition 2.9]{robertnicoud2020higher} with the addition of the arity $0$ operation $p_n \circ u$.
\end{Remark}

Let us compute the $u\mathcal{CC}_\infty$-coalgebra structure on its linear dual $C^c_*(\Delta^{n})$. The elementary decompositions 
\[
\Big\{ \Delta_\tau: C^c_*(\Delta^{n}) \longrightarrow C^c_*(\Delta^{n})^{\otimes m} \Big\}
\]
which determine the $u\mathcal{CC}_\infty$-coalgebra structure of $C^c_*(\Delta^{n})$ are given by $\Delta_\tau \coloneqq (\mu_\tau)^*$. Let $\{a_I\} = \{\omega_I^*\}$ be the dual basis of $C^c_*(\Delta^{n})$. Let $\tau$ be a rooted tree in $\mathrm{RT}_m^\omega$. For any $m$-tuple $(I_1,\cdots,I_m)$ of non-trivial subsets of $[n]$, suppose we have that
\[
\mu_\tau(\omega_{I_1},\cdots,\omega_{I_m}) = \sum_{I \subseteq [n]} \lambda_I^{\tau,(I_1,\cdots,I_m)}\omega_I~.
\]
where $\{\omega_I\}$ is the basis of $C^*_c(\Delta^{n})$. Then
\[
\Delta_\tau(a_I) = \sum_{\substack{I_1,\cdots,I_m \subseteq [n],I_l \neq \emptyset \\ \lambda_I^{\tau,(I_1,\cdots,I_m)} \neq 0}} \frac{1}{\lambda_I^{\tau,(I_1,\cdots,I_m)}} a_{I_1} \otimes \cdots \otimes a_{I_m}~. 
\]
Notice that $a_{I_1} \otimes \cdots \otimes a_{I_m}$ appears in the decomposition of $a_I$ under $\Delta_\tau$ if and only if $\omega_I$ appears with a non-zero coefficient in $\mu_\tau(\omega_{I_1},\cdots,\omega_{I_m})$. 

\begin{Remark}
Computing the coefficients $\lambda_I^{\tau,(I_1,\cdots,I_m)}$ is a hard task, see \cite[Remark 5.11]{robertnicoud2020higher}.
\end{Remark}

\textbf{Complete cobar construction.} Let $(C,\Delta_C,d_C)$ be a $u\mathcal{CC}_\infty$-coalgebra, where 
\[
\Delta_C: C \longrightarrow \displaystyle \prod_{n \geq 0} \mathrm{Hom}_{\mathbb{S}_n}\left(\Omega\mathrm{B}^{\mathrm{s.a}}\ucom(n), C^{\otimes n}\right)
\]
denotes its structural morphism. The complete cobar construction $\widehat{\Omega}_{\iota} C $ is given by the underlying graded module 
\[
\prod_{n \geq 0} \widehat{\Omega}^{\mathrm{s.c}}\ucom^*(n) ~ \widehat{\otimes}_{\mathbb{S}_n} ~ C^{\otimes n}~.
\]
It is endowed with a pre-differential $d_{\mathrm{cobar}}$ given by the difference of $d_1$ and $d_2$. The first term $d_1$ is given by 
\[
d_1= - \widehat{\mathscr{S}}^c(d_{\widehat{\Omega}^{\mathrm{s.c}}\ucom^*})(\mathrm{id}) + \widehat{\mathscr{S}}^c(\mathrm{id})(\diracComb(\mathrm{id},d_C))~.
\]

The second term $d_2$ is induced by the $u\mathcal{CC}_\infty$-coalgebra structure of $C$, it is the unique derivation extending the map $\varphi$, given by the following composites
\[
\begin{tikzcd}[column sep=4pc,row sep=3pc]
C \arrow[d,"\Delta_C",swap] \arrow[r,"\varphi"]
&\displaystyle \prod_{n \geq 0} \widehat{\Omega}^{\mathrm{s.c}}\ucom^*(n) ~ \widehat{\otimes}_{\mathbb{S}_n} ~ C^{\otimes n} \\
\displaystyle \prod_{n \geq 0} \mathrm{Hom}_{\mathbb{S}_n}\left(\Omega\mathrm{B}^{\mathrm{s.a}}\ucom(n), C^{\otimes n}\right)  \arrow[r,"\widehat{\mathscr{S}}^c(\iota)(\mathrm{id})"] 
&\displaystyle \prod_{n \geq 0} \mathrm{Hom}_{\mathbb{S}_n}\left(\mathrm{B}^{\mathrm{s.a}}\ucom(n), C^{\otimes n}\right)~. \arrow[u,"\cong"]
\end{tikzcd}
\]

\begin{theorem}\label{theom: formulas for the pre differentials of mc}
Let $I = \{i_0,\cdots,i_k\} \subseteq [n]$ with $k > 0$, and $a_I$ be the corresponding basis element in $C^c_*(\Delta^n)$. When $\mathrm{Card}(I) \geq 2$, the image of $a_I$ under the pre-differential $d_{\mathrm{cobar}}$ is given by 
\[
d_{\mathrm{cobar}}(a_I) = \sum_{l=0}^k (-1)^l a_{i_0\cdots\widehat{i_l}\cdots i_k} - \sum_{m \geq 2}\sum_{\substack{\tau \in \mathrm{RT}_m \\ I_1,\cdots,I_m \subseteq [n], I_l \neq \emptyset \\ \lambda_I^{\tau,(I_1,\cdots,I_m)} \neq 0}}\frac{1}{\mathcal{E}(\tau)\lambda_I^{\tau,(I_1,\cdots,I_m)}} \tau(a_{I_1},\cdots, a_{I_m})~,
\]

for $I = \{i_0,\cdots,i_k\}$, where $\mathcal{E}(\tau)$ is the renormalization coefficient of Remark \ref{Remark: renormalization}. When $\mathrm{Card}(I) = 1$, the image of $a_I$ under the pre-differential is given by 
\[
-\sum_{n \geq 0,~ n \neq 1} \frac{1}{n!} c_n(a_0,\cdots,a_0)~,
\]
where $c_n$ denotes the $n$-corolla. 
\end{theorem}

\begin{proof}
By Proposition \ref{prop: formulas for the HTT}, the only arity $0$ operation on $C^c_*(\Delta^n)$ is given by $p_n \circ u$. In order to know in which decompositions the lonely cork is going to appear, it is enough to compute the image of the unit $1$ in  $\Omega_n$ by the morphism $p_n$. Since
\[
p_n(1) = t_0 + \cdots + t_n \in C^*_c(\Delta^n)
\]
for all $n \geq 0$, the cork only appears in the decomposition of the elements of the form $a_{i}$ where $\{i\} \subseteq [n]$. Consequently, all the formulas for $d_{\mathrm{cobar}}(a_I)$ where $\mathrm{Card}(I) \geq 2$ are only indexed by rooted trees without corks. The renormalization coefficient $\mathcal{E}(\tau)$ appears because an identification between invariants and coinvariants was done in the definition of $d_2$.
\end{proof}

\begin{Corollary}\label{cor: same formulas as for RNV}
When $\mathrm{Card}(I) \geq 2$, the formula for $d_{\mathrm{cobar}}(a_I)$ coincides with the formula computed in \cite[Section 2.2]{robertnicoud2020higher} in the non-curved case.
\end{Corollary}

\begin{proof}
Follows directly from the above theorem.
\end{proof}

\begin{Corollary}
Let $(\mathfrak{g},\gamma_\mathfrak{g},d_\mathfrak{g})$ be a curved absolute $\mathcal{L}_\infty$-algebra. There is a bijection
\[
\mathcal{R}(\mathfrak{g})_0 \cong \mathcal{MC}(\mathfrak{g})~,
\]
between the $0$-simplices of $\mathcal{R}(\mathfrak{g})$ and the Maurer--Cartan elements of $\mathfrak{g}$.
\end{Corollary}

\begin{proof}
At $n=0$, the Dupont contraction is trivial since $h_0 = 0$ and $i_0$ and $p_0$ are isomorphism. Thus the resulting $u\mathcal{CC}_\infty$-coalgebra structure on $\kk.a_0$ is simply given by 
\[
\begin{tikzcd}[column sep=4pc,row sep=-0.8pc]
\kk.a_0 \arrow[r]
&\displaystyle \prod_{n \geq 0} \mathrm{Hom}_{\mathbb{S}_n}\left(\Omega\mathrm{B}^{\mathrm{s.a}}\ucom(n), (\kk.a_0)^{\otimes n}\right) \\
a_0 \arrow[r,mapsto]
&\displaystyle \sum_{n \geq 0,~ n \neq 1} \left[\mathrm{ev}_{a_0}: c_n \mapsto a_0 \otimes \cdots \otimes a_0 \right]~.
\end{tikzcd}
\]
Thus 
\[
\mathfrak{mc}^0 \cong \left(\displaystyle \prod_{n \geq 0} \widehat{\Omega}^{\mathrm{s.c}}\ucom^*(n) ~ \widehat{\otimes}_{\mathbb{S}_n} ~ (\kk.a_0)^{\otimes n}, d_{\mathrm{cobar}}(a_0) = - \sum_{n \geq 0,~ n \neq 1} \frac{1}{n!} c_n(a_0,\cdots,a_0) \right)~.
\]

The data of a morphism of curved absolute $\mathcal{L}_\infty$-algebras $f: \mathfrak{mc}^0 \longrightarrow \mathfrak{g}$ is equivalent to the data of $\alpha$ in $\mathfrak{g}$ such that 
\[
d_\mathfrak{g}(\alpha) = - \gamma_\mathfrak{g}\left(\sum_{n \geq 0,~ n \neq 1} \frac{1}{n!} c_n(\alpha,\cdots,\alpha)\right)~,
\]
since $\mathfrak{mc}^0$ is freely generated by $a_0$ as a pdg $\mathrm{B}^{\mathrm{s.a}}\ucom$-algebra, and since the data of a morphism of pdg $\mathrm{B}^{\mathrm{s.a}}\ucom$-algebras is the same as the data of a morphism of curved $\mathrm{B}^{\mathrm{s.a}}\ucom$-algebras.
\end{proof}

\subsection{Gauge actions}
In this subsection, we define intrinsically the gauge action of an element $\lambda$ of degree $1$ on the set of Maurer-Cartan elements. For a complete curved absolute $\mathcal{L}_\infty$-algebra $\mathfrak{g}$, gauge actions correspond to paths in $\mathcal{R}(\mathfrak{g})$, thus we obtain an explicit description of $\pi_0(\mathcal{R}(\mathfrak{g}))$. 

\medskip

Let $\mathfrak{g}$ be a complete curved absolute $\mathcal{L}_\infty$-algebra. We can consider the following continuous function 
\[
\varphi_{\mathcal{MC}}(-) \coloneqq \left[d_\mathfrak{g}(-) + \gamma_\mathfrak{g}\left(\sum_{n \geq 0,~ n \neq 1} \frac{1}{n!} c_n(-,\cdots,-)\right)\right]: \mathfrak{g}_0 \longrightarrow \mathfrak{g}_{-1}
\]
between two complete vector spaces. The Maurer-Cartan set of $\mathfrak{g}$ is given by the $0$ locus in $\mathfrak{g}_{-1}$. Would this function be smooth or algebraic, the resulting set would be a differential or an algebraic variety. One can still formally consider its "tangent space" at a given point $\alpha$ in $\mathcal{MC}(\mathfrak{g})$, which in this case is given by the kernel of $d^\alpha_\mathfrak{g}$, where $d^\alpha_\mathfrak{g}$ is the \textit{twisted differential} by $\alpha$ given by 
\[
d^\alpha_\mathfrak{g}(-) \coloneqq d_\mathfrak{g}(-) + \gamma_\mathfrak{g}\left(\sum_{n \geq 2} \frac{1}{(n-1)!} c_n(\alpha,\cdots,\alpha,-)\right)~.
\]

\begin{lemma}
Let $\mathfrak{g}$ be a curved absolute $\mathcal{L}_\infty$-algebra and let $\alpha$ be a Maurer-Cartan element. Then $d^\alpha_\mathfrak{g}$ squares to zero.
\end{lemma}

\begin{proof}
The formal algebraic computations of \cite[Corollary 5.1]{DSV18} extend \textit{mutatis mutandis} to this case using the associativity Condition \ref{associativity condition}. 
\end{proof}

Therefore an easy way to produce elements in the tangent space $\mathrm{T}_\alpha \mathcal{MC}(\mathfrak{g}) = \mathrm{Ker}(d^\alpha_\mathfrak{g})$ is to take a degree one element $\lambda$ in $\mathfrak{g}$ and consider $d^\alpha_\mathfrak{g}(\lambda)$. Then, given this tangent vector at $\alpha$, one would like to "integrate it" inside $\mathcal{MC}(\mathfrak{g})$. As explained in \cite{DSV18}, these heuristics motivate the following definitions.

\begin{lemma}
Let $(\mathfrak{g},\gamma_\mathfrak{g},d_\mathfrak{g})$ be a complete curved absolute $\mathcal{L}_\infty$-algebra. The graded module  
\[
\kk[[t]] \widehat{\otimes} \mathfrak{g} \coloneqq \lim_\omega \left( \kk[t]/(t^\omega) \otimes \mathfrak{g}/\mathrm{W}_\omega \mathfrak{g} \right) \cong \Bigg\{ \sum_{\omega \geq 0} g_\omega \otimes t^\omega ~~| ~~ g_\omega \in \mathrm{W}_\omega \mathfrak{g} \Bigg\}
\]
has a complete curved absolute $\mathcal{L}_\infty$-algebra structure given by
\[
\gamma_{\kk[[t]] \widehat{\otimes} \mathfrak{g}} \left(\sum_{k \geq 0} \sum_{\omega \geq 0} \sum_{\tau \in \mathrm{CRT}_k^\omega} \lambda_\tau \tau\left(\sum_{i_1 \geq 0}g_{i_1} \otimes t^{i_1}, \cdots, \sum_{i_k \geq 0}g_{i_k} \otimes t^{i_k}\right)\right) \coloneqq 
\]
\[
 \sum_{n \geq 0} \sum_{k \geq 0} \sum_{\omega \geq 0}  \left( \sum_{i_1 + \cdots + i_k = n} \gamma_\mathfrak{g} \left(\sum_{\tau \in \mathrm{CRT}_n^\omega} \lambda_\tau \tau(g_{i_1}, \cdots, g_{i_k}) \right)\right) \otimes t^{n + \omega}~. 
\]
\end{lemma}

\begin{proof}
One directly and easily can check that this structure satisfies the axioms. 
\end{proof}

\begin{Remark}
This structure comes from the convolution complete curved absolute $\mathcal{L}_\infty$-algebra structure structure on 
\[
\mathrm{Hom}_{\mathsf{gr}\textsf{-}\mathsf{mod}}(\kk^c[t],\mathfrak{g})~,
\]
where $\kk^c[t]$ is the conilpotent cofree cocommutative coalgebra cogenerated by $t$ in degree $0$.
\end{Remark}

The algebraic way to integrate a degree one element $\lambda$ in $\mathfrak{g}$ amounts to solving the following ordinary differential equation

\begin{equation}\label{edo jauges}
\frac{\mathrm{d}\gamma(t)}{\mathrm{d}t} = d^{\gamma(t)}_\mathfrak{g}(\lambda)~,
\end{equation}

where $\gamma(t)$ is an element of $\kk[[t]]~ \widehat{\otimes}~ \mathfrak{g}$. Here $d^{\gamma(t)}_\mathfrak{g}$ is again given by 
\[
d^{\gamma(t)}_\mathfrak{g}(\lambda) \coloneqq  d_\mathfrak{g}(\lambda) + \gamma_{\kk[[t]] \widehat{\otimes} \mathfrak{g}}\left(\sum_{n \geq 2}  \frac{1}{(n-1)!} c_k \left(\gamma(t),\cdots,\gamma(t),\lambda\right) \right)~.
\]
where the term $\gamma(t)$ appears $n-1$ times on the right. The solution obviously depends on the initial value $\gamma(0)$. Let us denote it by $\alpha \coloneqq \gamma(0)$.

\begin{lemma}
Let $\mathfrak{g}$ be a curved absolute $\mathcal{L}_\infty$-algebra. Let $\lambda$ be in $\mathfrak{g}_{1}$. For every $\alpha$ in $\mathfrak{g}_0$, the above differential equation has an unique solution, given by: 
\[
\gamma(t) = \sum_{\omega \geq 0} \gamma_\mathfrak{g}\left(\sum_{n \geq 0} \sum_{\tau \in \mathrm{RT}^\omega_n} \frac{|\mathrm{Aut}(\tau)|}{W(\tau)} \tau^\lambda(\alpha, \cdots, \alpha) \right) \otimes t^\omega~,
\]
where $W(\tau)$ is the exponential weight function defined in \cite[Appendix A.2]{robertnicoud2020higher}.
\end{lemma} 

\begin{proof}
We can apply \cite[Proposition A.11]{robertnicoud2020higher} to this context. This sum is here indexed by rooted trees instead of planar trees, and thus rooted trees automorphisms are taken into account, see \cite[Proposition 1.10]{robertnicoud2020higher} for the planar formula.
\end{proof}

\begin{Definition}[Gauge actions]
Let $\mathfrak{g}$ be a complete curved absolute $\mathcal{L}_\infty$-algebra. Let $\lambda$ be in $\mathfrak{g}_{1}$ and let $\alpha$ be a Maurer-Cartan element of $\mathfrak{g}$. The \textit{gauge action} of $\lambda$ on $\alpha$ is given by the element
\[
\lambda \bullet \alpha \coloneqq \gamma(1) = \sum_{\omega \geq 0} \gamma_\mathfrak{g}\left(\sum_{n \geq 0} \sum_{\tau \in \mathrm{RT}^\omega_n} \frac{|\mathrm{Aut}(\tau)|}{W(\tau)} \tau^\lambda(\alpha, \cdots, \alpha) \right)~,
\] 
in $\mathfrak{g}_0$, where $\gamma(t)$ is the unique solution to Equation \ref{edo jauges} having $\alpha$ as initial condition.
\end{Definition}

\begin{Definition}[Gauge equivalence]
Two Maurer-Cartan elements $\alpha$ and $\beta$ of $\mathfrak{g}$ are \textit{gauge equivalent} if there exists $\lambda$ in $\mathfrak{g}_1$ such that $\lambda \bullet \alpha = \beta$. 
\end{Definition}

\begin{Proposition}
The gauge equivalence relation on the set of Maurer-Cartan elements $\mathcal{MC}(\mathfrak{g})$ defines an equivalence relation. 
\end{Proposition}

\begin{proof}
Follows from Proposition \ref{prop: equivalence between gauge and paths} using the fact that paths equivalence in $\mathcal{R}(\mathfrak{g})$ is an equivalence relation on $\mathcal{MC}(\mathfrak{g})$ since it is a Kan complex.
\end{proof}

\begin{Proposition}\label{prop: equivalence between gauge and paths}
Let $\mathfrak{g}$ be a complete curved absolute $\mathcal{L}_\infty$-algebra and let $\alpha$ and $\beta$ be two Maurer-Cartan elements of $\mathfrak{g}$. The data of a gauge equivalence $\lambda$ between $\alpha$ and $\beta$ is equivalent to the data of a morphism of complete curved absolute $\mathcal{L}_\infty$-algebras 
\[
\varphi: \widehat{\Omega}_\iota(C_*^c(\Delta^1)) \longrightarrow \mathfrak{g}
\]
such that $\varphi(a_{01})= \lambda$, $\varphi(a_0) = \alpha$ and $\varphi(a_1) = \beta$. 
\end{Proposition}

\begin{proof}
A morphism of complete curved absolute $\mathcal{L}_\infty$-algebras
\[
\varphi: \widehat{\Omega}_\iota C_*^c(\Delta^1)  \longrightarrow \mathfrak{g}
\]
is equivalent to a twisting morphism 
\[
\zeta_\varphi: C_*^c(\Delta^1) \longrightarrow \mathfrak{g}~.
\]
By Corollary \ref{cor: same formulas as for RNV}, the decomposition of the element $a_{01}$ in $C_*^c(\Delta^1)$ under the coalgebra structure gives the same formulae as in \cite{robertnicoud2020higher}. The data of a twisting morphism $\zeta_\varphi$ is thus equivalent to the data of two Maurer-Cartan elements $\alpha$ and $\beta$ in $\mathfrak{g}$ and to the data of $\lambda$ in $\mathfrak{g}_1$ such that $\lambda \bullet \alpha = \beta$, see \cite[Proposition 4.3]{daniel19} for an analogue statement.
\end{proof}

\begin{theorem}\label{thm: characterisation pi zero et jauges}
Let $\mathfrak{g}$ be a complete curved absolute $\mathcal{L}_\infty$-algebra. There is a bijection 
\[
\pi_0(\mathcal{R}(\mathfrak{g})) \cong \mathcal{MC}(\mathfrak{g})/\sim_{\mathrm{gauge}}~,
\]
where the right-hand side denotes the set of Maurer-Cartan elements up to gauge equivalence. 
\end{theorem}

\begin{proof}
Follows directly from Proposition \ref{prop: equivalence between gauge and paths}.
\end{proof}

\begin{Remark}
Let $\mathfrak{g}$ be a complete curved absolute $\mathcal{L}_\infty$-algebra. Morphisms of the form
\[
\varphi: \widehat{\Omega} C_*^c(\Delta^1) \longrightarrow \mathfrak{g}
\]
are in a one-to-one correspondence with morphisms of $u\mathcal{CC}_\infty$-coalgebras
\[
\varphi^{\sharp}: C_*^c(\Delta^1) \longrightarrow \widehat{\mathrm{B}}_\iota \mathfrak{g}~.
\]
Since $C_*^c(\Delta^1)$ is the interval object in the model category of $u\mathcal{CC}_\infty$-coalgebras, the data of a gauge equivalence between two Maurer-Cartan elements is equivalent to the data of a left homotopy of morphisms of $u\mathcal{CC}_\infty$-coalgebras between the two associated group-like elements $\varphi^{\sharp}|_{a_0}$ and $\varphi^{\sharp}|_{a_1}$.
\end{Remark}

\subsection{Higher Baker--Campbell--Hausdorff formulae}\label{subsection: higher BCH}
In this subsection, we generalize the seminal results about the higher Baker--Campbell--Hausdorff formulae for complete $\mathcal{L}_\infty$-algebras developed in \cite[Section 5]{robertnicoud2020higher} to the case of curved absolute $\mathcal{L}_\infty$-algebras. The higher Baker--Campbell--Hausdorff provide us with \textit{explicit} formulae for the horn-fillers in $\mathcal{R}(\mathfrak{g})$, making it not only an $\infty$-groupoid, but an \textit{algebraic} $\infty$-groupoid with \textit{canonical} horn-fillers.

\medskip

Consider the graded pdg module $U(n)$ which is given by $\kk.u$ in degree $n$ and by $\kk.du$ in degree $n-1$, together with the pre-differential $d(u) \coloneqq du$. Notice that 

\[
\widehat{\mathscr{S}}^c(\mathrm{B}^{\mathrm{s.a}}\ucom)(U(n)) \cong  \prod_{m \geq 0} \widehat{\Omega}^{\mathrm{s.a}}\ucom^*(m) ~\widehat{\otimes}_{\mathbb{S}_m}~ U(n)^{\otimes m}
\]

is the free complete pdg $\mathrm{B}^{\mathrm{s.a}}\ucom$-algebra generated by $U(n)$, but it \textbf{does not} form a complete \textit{curved} $\mathrm{B}^{\mathrm{s.a}}\ucom$-algebra. 

\begin{lemma}\label{lemma: decomposition de L(Delta n)}
There is an isomorphism of complete pdg $\mathrm{B}^{\mathrm{s.a}}\ucom$-algebras 
\[
\mathcal{L}(\Delta^n) \cong \widehat{\mathscr{S}}^c(\mathrm{B}^{\mathrm{s.a}}\ucom)(U(n)) \amalg \mathcal{L}(\Lambda_k^n)
\]
for all $n \geq 2$ and for $1 \leq k \leq n$, where $\Lambda_k^n$ denotes the $k$-horn of dimension $n$. 

\end{lemma}

\begin{proof}
The proof is \textit{mutatis mutandis} the same as that of \cite[Lemma 5.1]{robertnicoud2020higher}, where this time we use the completeness of the underlying canonical filtration instead. The bijections constructed in \textit{loc.cit.} also work in this case, and the formula for the element $a_{\widehat{k}}$ is the same by Theorem \ref{theom: formulas for the pre differentials of mc}. 
\end{proof}

\begin{theorem}\label{thm: canonical horn fillers}
Let $\mathfrak{g}$ be a complete curved absolute $\mathcal{L}_\infty$-algebra. There is a bijection

\[
\mathrm{Hom}_{\mathsf{sSet}}(\Delta^n, \mathcal{R}(\mathfrak{g})) \cong  \mathfrak{g}_n \times \mathrm{Hom}_{\mathsf{sSet}}(\Lambda_k^n,\mathcal{R}(\mathfrak{g}))~,
\]

\vspace{0.5pc}
natural in $\mathfrak{g}$ with respect to morphisms of curved absolute $\mathcal{L}_\infty$-algebras.
\end{theorem}

\begin{proof}
Consider the following isomorphisms

\begin{align*}
\mathrm{Hom}_{\mathsf{sSet}}(\Delta^n, \mathcal{R}(\mathfrak{g})) &\cong \mathrm{Hom}_{\mathsf{curv}~\mathsf{abs}~\mathcal{L}_\infty\text{-}\mathsf{alg}}(\mathcal{L}(\Delta^n), \mathfrak{g}) \\
&\cong \mathrm{Hom}_{\mathsf{pdg}~\mathrm{B}^{\mathrm{s.a}}\ucom\text{-}\mathsf{alg}}(\mathcal{L}(\Delta^n), \mathfrak{g}) \\
&\cong \mathrm{Hom}_{\mathsf{pdg}~\mathrm{B}^{\mathrm{s.a}}\ucom\text{-}\mathsf{alg}}(\widehat{\mathscr{S}}^c(\mathrm{B}^{\mathrm{s.a}}\ucom)(U(n)) \amalg \mathcal{L}(\Lambda_k^n), \mathfrak{g}) \\
&\cong  \mathrm{Hom}_{\mathsf{pdg}~\mathrm{B}^{\mathrm{s.a}}\ucom\text{-}\mathsf{alg}}(\widehat{\mathscr{S}}^c(\mathrm{B}^{\mathrm{s.a}}\ucom)(U(n)), \mathfrak{g}) \times  \mathrm{Hom}_{\mathsf{pdg}~\mathrm{B}^{\mathrm{s.a}}\ucom\text{-}\mathsf{alg}}(\mathcal{L}(\Lambda_k^n), \mathfrak{g}) \\
&\cong \mathfrak{g}_n \times  \mathrm{Hom}_{\mathsf{curv}~\mathsf{abs}~\mathcal{L}_\infty\text{-}\mathsf{alg}}(\mathcal{L}(\Lambda_k^n), \mathfrak{g}) \\
&\cong \mathfrak{g}_n \times \mathrm{Hom}_{\mathsf{sSet}}(\Lambda_k^n,\mathcal{R}(\mathfrak{g}))~,
\end{align*}

where we used several times that curved $\mathrm{B}^{\mathrm{s.a}}\ucom$-algebras are a full subcategory of pdg $\mathrm{B}^{\mathrm{s.a}}\ucom$-algebras, and where we used Lemma \ref{lemma: decomposition de L(Delta n)} in line three.
\end{proof}

\begin{Corollary}
Let $\mathfrak{g}$ be a complete curved absolute $\mathcal{L}_\infty$-algebra. The $\infty$-groupoid $\mathcal{R}(\mathfrak{g})$ has a \textit{canonical structure} of $\infty$-groupoid. 
\end{Corollary}

\begin{proof}
For any lifting problem $\Lambda_k^n \longrightarrow \mathcal{R}(\mathfrak{g})$ there is a canonical lift given by the image of $0$ inside $\mathfrak{g}_n$ under the bijection of Theorem \ref{thm: canonical horn fillers}.
\end{proof}

\begin{Remark}
Thus the integration functor $\mathcal{R}$ lands on the category of \textit{algebraic} $\infty$-groupoids in the sense of \cite{algebraic}. This shifts the point of view in the following way: it not only satisfies a \textit{property} (existence of horn-filler) but instead these horn-fillers a given and thus are considered as a \textit{structure}.
\end{Remark}

An element in $\mathcal{R}(\mathfrak{g})_n$ amounts to the data of a curved twisting morphism
\[
\varphi: C_*^c(\Delta^n) \longrightarrow \mathfrak{g}~.
\]
It can be pictorially represented by geometric $n$-simplex where the $\mathrm{I}$-th face $\mathrm{I} \subseteq [n]$ of the $n$-simplex is labeled by the element $\varphi(a_\mathrm{I})$ in $\mathfrak{g}$. As an example, an element in $\mathcal{R}(\mathfrak{g})_2$ can be represented as 

\begin{center}
\includegraphics[width=60mm,scale=0.7]{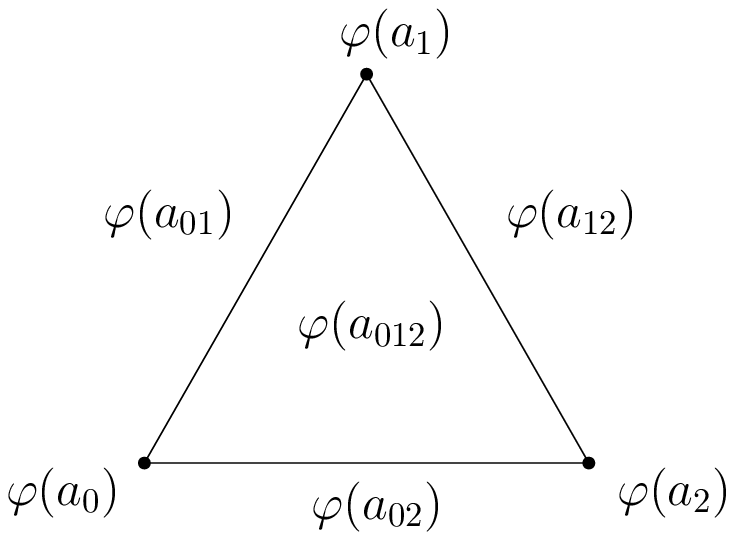},
\end{center}

where $\varphi(a_{i})$ are Maurer-Cartan elements in $\mathfrak{g}_0$ for $0 \leq i \leq 2$, where $\varphi(a_{ij})$ are degree $1$ elements in $\mathfrak{g}$ which induce gauge equivalences and where $\varphi(a_{012})$ is a degree $2$ element satisfying a higher compatibility condition imposed by the curved twisting morphism condition on $\varphi$. Similarly, the $\Lambda_k^n$-horns $\Lambda_k^n \longrightarrow \mathcal{R}(\mathfrak{g})$ admit an analogue pictorial description. As an example, a $\Lambda_1^2$-horn  in $\mathcal{R}(\mathfrak{g})$ can be depicted as

\begin{center}
\includegraphics[width=60mm,scale=0.7]{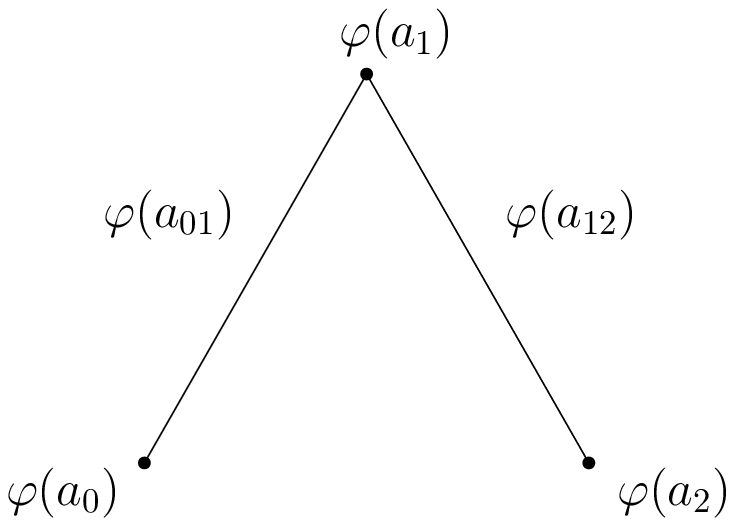},
\end{center}

\begin{Definition}[Higher Baker--Campbell--Hausdorff products]
Let $\mathfrak{g}$ be a complete curved absolute $\mathcal{L}_\infty$-algebra and let $y$ be an element in $\mathfrak{g}$ of degree $n$. The \textit{higher Baker--Campbell--Hausdorff product relative to} $y$ is given by 
\[
\begin{tikzcd}[column sep=4pc,row sep=-0pc]
\Gamma^y: \mathrm{Hom}_{\mathsf{sSet}}(\Lambda_k^n,\mathcal{R}(\mathfrak{g})) \arrow[r]
&\mathfrak{g}_{n-1} \\
x \arrow[r,mapsto]
&\Gamma^y(x) \coloneqq \varphi^{(x;y)}(a_{\widehat{k}})~.
\end{tikzcd}
\]
It sends any $\Lambda_k^n$-horn $x: \Lambda_k^n \longrightarrow \mathcal{R}(\mathfrak{g})$ to the image $\varphi^{(x;y)}(a_{\widehat{k}})$ in $\mathfrak{g}_{n-1}$, where $\varphi^{(x;y)}$ is the curved twisting morphism associated to the element in $\mathcal{R}(\mathfrak{g})_n$ given by the bijection constructed in Theorem \ref{thm: canonical horn fillers}. 
\end{Definition}

These higher Baker--Campbell--Hausdorff products are given by explicit formulae. In our context, we recover the same formulas as the ones obtained in \cite[Proposition 5.10]{robertnicoud2020higher}.

\begin{Proposition}\label{prop: formules higher BCH}
Let $\mathfrak{g}$ be a complete curved absolute $\mathcal{L}_\infty$-algebra. Let $x: \Lambda_k^n \longrightarrow \mathcal{R}(\mathfrak{g})$ be a $\Lambda_k^n$-horn in $\mathcal{R}(\mathfrak{g})$. The higher Baker--Campbell--Hausdorff product relative to an element $y$ in $\mathfrak{g}_n$ is given by the following formula

\[
\Gamma^y(x)= \gamma_\mathfrak{g}\left(\sum_{\begin{subarray}{c}\tau\in\mathrm{PaPRT}\\
\chi\in\mathrm{Lab}^{[n],k}(\tau)\end{subarray}}\ \prod_{\begin{subarray}{c}%
\beta\text{ block of }\tau\\
\lambda^{\beta(\chi)}_{[n]}\neq 0\end{subarray}}\frac{(-1)^{k}}{\lambda^{\beta%
(\chi)}_{[n]}[\beta]!}\, \tau\left(x_{\chi(1)},\ldots,x_{\chi(p)};(-1)^{%
k}d_\mathfrak{g}(y)-\sum_{l\neq k}(-1)^{k+l}x_{\widehat{l}}\right) \right)~.
\]
See \cite[Proposition 5.10]{robertnicoud2020higher} for more details on the coefficients of this formula.
\end{Proposition}

\begin{proof}
The isomorphism of Lemma \ref{lemma: decomposition de L(Delta n)} is given by the same formulae that in \cite{robertnicoud2020higher}, and this isomorphism determines the higher Baker--Campbell-Hausdorff products.
\end{proof}

\begin{Remark}
Using the completeness of $\mathfrak{g}$, this last formula can be rewritten by splitting it along the weight filtration. 
\end{Remark}

These higher Baker--Campbell-Hausdorff products do recover the classical Baker--Campbell-Hausdorff product when one restricts to the case of nilpotent Lie algebras. Remember that any nilpotent Lie algebra is a particular example of a curved absolute $\mathcal{L}_\infty$-algebra, as explained in Example \ref{Example: nilpotent stuff}. In fact, the data of a nilpotent Lie algebra is equivalent to the data of a nilpotent $\mathcal{L}_\infty$-algebra concentrated in degree $1$.

\medskip

Let $\mathfrak{g}$ be a nilpotent Lie algebra viewed in degree $1$, and let $x: \Lambda_1^2 \longrightarrow \mathcal{R}(\mathfrak{g})$ be a $\Lambda_1^2$-horn. Since $\mathfrak{g}$ is concentrated in degree $1$, the only Maurer-Cartan element is $0$ and the data of such a horn amount to the data of two elements $\alpha$ and $\beta$ in $\mathfrak{g}$. 
\vspace{0.4pc}

\begin{center}
\includegraphics[width=45mm,scale=0.5]{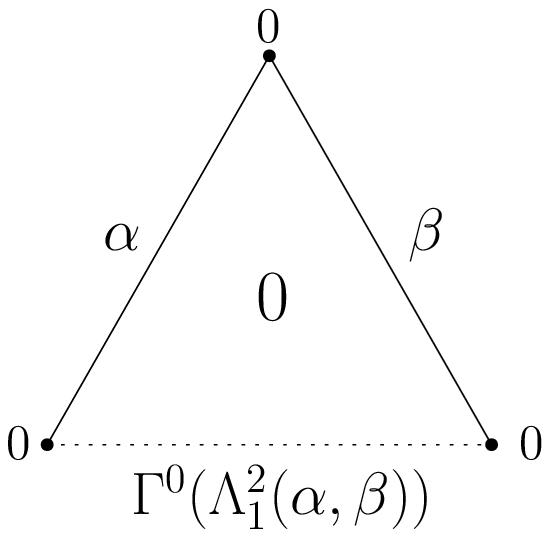},
\end{center}

Since $\mathfrak{g}_2 = \{0\}$, there is only one way to fill this horn. This filling corresponds to the classical Baker--Campbell--Hausdorff formula.

\begin{Proposition}\label{prop: horn 2,1 c'est BCH}
Let $\mathfrak{g}$ be a nilpotent Lie algebra, viewed in degree $1$. Let $\alpha$ and $\beta$ be two degree $1$ elements, then 
\[
\Gamma^{0}(\Lambda_1^2(\alpha,\beta)) = \mathrm{BCH}(\alpha,\beta)~,
\]
where $\Lambda_1^2(\alpha,\beta)$ stands for the horn labeled by $\alpha$ and $\beta$ and where $\mathrm{BCH}(\alpha,\beta)$ stands for the classical Baker--Campbell-Hausdorff formula.
\end{Proposition}

\begin{proof}
By Proposition \ref{prop: formules higher BCH}, we obtain the same formulas as in \cite{robertnicoud2020higher}. By Theorem 5.17 in \textit{loc.cit.}, the formula this higher BCH product coincides with the classical BCH formula. 
\end{proof}

\begin{Remark}
This above result was first proved in R. Bandiera's PhD Thesis, see \cite{bandiera}.
\end{Remark}

\begin{Remark}
It is known that the Baker--Campbell--Hausdorff formula lives in the free complete Lie algebra generated by two elements. This algebra is in fact the \textit{free absolute} Lie algebra, hence the Baker--Campbell--Hausdorff formula naturally belongs to the realm of absolute algebras. 
\end{Remark}

Let $\mathfrak{g}$ be a nilpotent Lie algebra. We denote 
\[
\mathrm{Exp}(\mathfrak{g}) \coloneqq \left(\mathfrak{g}, \mathrm{BCH}, 0 \right)
\]
the group obtained by considering $\mathfrak{g}$ equipped with the BCH formula. 

\begin{Corollary}\label{cor: integration of absolute Lie}
Let $\mathfrak{g}$ be a nilpotent Lie algebra, viewed in degree $1$. There is an isomorphism of simplicial sets 
\[
\mathcal{R}(\mathfrak{g}) \cong \mathcal{N}(\mathrm{B}\mathrm{Exp}(\mathfrak{g}))~,
\]
where $\mathcal{N}$ stands for the nerve of a category.
\end{Corollary}

\begin{proof}
By Theorem \ref{thm: canonical horn fillers}, $n$-simplices in $\mathcal{R}(\mathfrak{g})$ correspond $n$-tuples of elements in $\mathfrak{g}$, and Proposition \ref{prop: horn 2,1 c'est BCH} identifies the composition of paths with the group structure of $\mathrm{Exp}(\mathfrak{g})$.
\end{proof}

\begin{Remark}
The above corollary also follows from the results of \cite{Getzler09}, together with the comparison with Getzler's functor in the case of nilpotent $\mathcal{L}_\infty$-algebras mentioned in Remark \ref{Remark: comparison with Getzler functor}.
\end{Remark}

\begin{Corollary}
Let $\mathfrak{g}$ be a complete absolute Lie algebra, viewed in degree $1$. There is an isomorphism of simplicial sets 
\[
\mathcal{R}(\mathfrak{g}) \cong \mathcal{N}\left(\lim_\omega \mathrm{B}\mathrm{Exp}(\mathfrak{g}/\mathrm{W}_\omega \mathfrak{g})\right)~.
\]
It is therefore the nerve of a pro-$\kk$-unipotent group. 
\end{Corollary}

\begin{proof}
It follows from Corollary \ref{cor: integration of absolute Lie}, using the fact that the nerve functor commutes with limits.
\end{proof}

\subsection{Higher homotopy groups}
In this subsection, we compute the higher homotopy groups of the Kan complex $\mathcal{R}(\mathfrak{g})$. This is done by using the explicit $u\mathcal{CC}_\infty$-coalgebra structures on the cellular chains of the $n$-spheres. This new method allows us to generalize the one of the main results of A. Berglund \cite{Berglund15} to curved absolute $\mathcal{L}_\infty$-algebras. 

\medskip

We consider the model of the $n$-sphere $\mathbb{S}^n$ given by the simplicial set $\Delta^n/\partial \Delta^n$, which has one non-degenerate $0$-simplex $[0]$ and one non-degenerate $n$-simplex $[n]$. This simplicial set is not a Kan complex. Nevertheless, since for any complete curved absolute $\mathcal{L}_\infty$-algebra $\mathfrak{g}$, $\mathcal{R}(\mathfrak{g})$ is a Kan complex, one has

\[
\pi_n(\mathcal{R}(\mathfrak{g}),\alpha) \cong \mathrm{Hom}_{\mathsf{sSets}_*}(\mathbb{S}^n, \mathcal{R}(\mathfrak{g}))/\sim_{\mathrm{hmt}}~,
\]
\vspace{0.1pc}

where the right-hand side is the set of morphisms of simplicial sets which send $[0]$ to $\alpha$ in $\mathcal{R}(\mathfrak{g})_0$ modulo the homotopy relation. Now notice that 

\[
\mathrm{Hom}_{\mathsf{sSets}}(\mathbb{S}^n, \mathcal{R}(\mathfrak{g}))/\sim_{\mathrm{hmt}} ~\cong \mathrm{Hom}_{u\mathcal{CC}_\infty\text{-}\mathsf{coalg}}\left(C_*^c(\mathbb{S}^n), \widehat{\mathrm{B}}_\iota(\mathfrak{g}) \right)/\sim_{\mathrm{hmt}}
\]
\vspace{0.1pc}

since $C_*^c(-) \dashv \overline{\mathcal{R}}$ is a Quillen adjunction and since both $\mathcal{R}(\mathfrak{g})$ and $\widehat{\mathrm{B}}_\iota \mathfrak{g}$ are fibrant objects. Hence by explicitly computing the $u\mathcal{CC}_\infty$-coalgebra structure on $C_*^c(\mathbb{S}^n)$ one can compute the homotopy groups of $\mathcal{R}(\mathfrak{g})$ for any complete curved absolute $\mathcal{L}_\infty$-algebra. 

\begin{lemma}\label{lemma: structure de cogèbre sur la sphere}
Let $n \geq 1$. The dg module $C_*^c(\mathbb{S}^n)$ is given by $\kk.a_0$ in degree $0$ and $\kk.a_{[n]}$ in degree $n$ with zero differential. The $u\mathcal{CC}_\infty$-coalgebra structure on $C_*^c(\mathbb{S}^n)$ is given by the elementary decomposition maps 

\[
\left\{
\begin{tikzcd}[column sep=4pc,row sep=0.5pc]
a_0 \arrow[r,mapsto]
&\Delta_{c_m}(a_0) = a_0 \otimes \cdots \otimes a_0~, \\
a_{[n]} \arrow[r,mapsto]
&\Delta_{c_m}(a_{[n]}) = a_0 \otimes \cdots \otimes a_0 \otimes a_{[n]}~, \\
\end{tikzcd}
\right.
\]

where $c_m$ denotes the $m$-corolla, for all $m \geq 2$, and where $\Delta_\tau$ is the zero morphism for any other corked rooted tree $\tau$ in $\mathrm{CRT}$. 
\end{lemma}

\begin{proof}
Consider the the following pushout in the category of simplicial sets 
\[
\begin{tikzcd}
\partial \Delta^n  \arrow[r,twoheadrightarrow] \arrow[d,rightarrowtail] \arrow[dr, phantom, "\lrcorner", very near end]
&\{*\} \arrow[d] \\
\Delta^n \arrow[r,twoheadrightarrow]
&\mathbb{S}^n~. 
\end{tikzcd}
\]
It induces a pullback in the category of $u\mathcal{CC}_\infty$-algebras 
\[
\begin{tikzcd}
C^*_c(\mathbb{S}^n) \arrow[r] \arrow[d,rightarrowtail] \arrow[dr, phantom, "\ulcorner", very near start]
&C^*_c(\{*\}) \arrow[d] \\
C^*_c(\Delta^n) \arrow[r,twoheadrightarrow]
&C^*_c(\partial \Delta^n)~.
\end{tikzcd}
\]
Therefore there is a monomorphism $C^*_c(\mathbb{S}^n) \rightarrowtail C^*_c(\Delta^n)$ of $u\mathcal{CC}_\infty$-algebras. Using Proposition \ref{prop: formulas for the HTT}, one can see that, for degree reasons, the only operations that survive on $C^*_c(\mathbb{S}^n)$ are the multiplications $\mu_{c_m}$. Furthermore, it is clear that $\mu_{c_m}(a_0,\cdots, a_0, a_{[n]}) = a_{[n]}$. Computing the $u\mathcal{CC}_\infty$-coalgebra structure of the linear dual $C_*^c(\mathbb{S}^n)$ from this is straightforward.
\end{proof}

\begin{Definition}[$\alpha$-homology groups]
Let $\mathfrak{g}$ be a curved absolute $\mathcal{L}_\infty$-algebra and let $\alpha$ be a Maurer--Cartan element. We define the $\alpha$-homology groups of $\mathfrak{g}$ to be
\[
\mathrm{H}^\alpha_*(\mathfrak{g}) \coloneqq \frac{\mathrm{Ker}(d^\alpha_\mathfrak{g})}{\mathrm{Im}(d^\alpha_\mathfrak{g})}~.
\]
The $\alpha$-homology groups are functorial in morphisms of curved absolute $\mathcal{L}_\infty$-algebras: any morphism $f: \mathfrak{g} \longrightarrow \mathfrak{h}$ induces a morphism
\[
f^\alpha:  \mathrm{H}^\alpha_*(\mathfrak{g}) \longrightarrow \mathrm{H}^{f(\alpha)}_*(\mathfrak{h})~,
\]
as $f$ clearly preserves Maurer--Cartan elements and commutes with the twisted differentials.
\end{Definition}

\begin{Remark}\label{Remark: twisting procedure}
There should be a analogue notion of \textit{the twisting procedure} for curved absolute $\mathcal{L}_\infty$-algebras where given a Maurer-Cartan element $\alpha$, one can construct an \textit{absolute} $\mathcal{L}_\infty$-algebra $\mathfrak{g}^\alpha$. If this was the case, then we would have that 
\[
\mathrm{H}^\alpha_*(\mathfrak{g}) \cong \mathrm{H}_*(\mathfrak{g}^\alpha)~.
\]
Furthermore, given a curved absolute $\mathcal{L}_\infty$-algebra $\mathfrak{g}$, its image under the integration functor decomposes as 
\[
\mathcal{R}(\mathfrak{g}) \cong \underset{\alpha \in \pi_0(\mathfrak{R}(\mathfrak{g}))}{\coprod} \mathcal{R}(\mathfrak{g})^{\alpha}~,
\]
into a sum over its connected components. The twisted absolute $\mathcal{L}_\infty$-algebra $\mathfrak{g}^\alpha$ should be a model for the connected component $\mathcal{R}(\mathfrak{g})^{\alpha}$. 

\medskip

We could just set the twisted formulae in the \textit{absolute} setting and try to check by hand that they satisfy Condition \ref{pdg condition}, \ref{associativity condition} and \ref{curved condition}. We would nevertheless lack a conceptual explanation of these twisted algebras in terms of the action of a \textit{gauge group} as explained in \cite{DotsenkoShadrinVallette16}.
\end{Remark}

\begin{theorem}\label{thm: Berglund isomorphism}
Let $\mathfrak{g}$ be a complete curved absolute $\mathcal{L}_\infty$-algebra and let $\alpha$ be a Maurer--Cartan element. There is a bijection
\[
\pi_n(\mathcal{R}(\mathfrak{g}),\alpha) \cong \mathrm{H}^\alpha_n(\mathfrak{g})~,
\]
for all $n \geq 1$, which is natural in $\mathfrak{g}$.
\end{theorem}

\begin{proof}
Recall that the data of a morphism of $u\mathcal{CC}_\infty$-coalgebras $f_{\nu}: C_*^c(\mathbb{S}^n) \longrightarrow \widehat{\mathrm{B}}_\iota(\mathfrak{g})$ is equivalent to the data of a curved twisting morphism $\nu: C_*^c(\mathbb{S}^n) \longrightarrow \mathfrak{g}$ relative to $\iota$. In this case, by Lemma \ref{lemma: structure de cogèbre sur la sphere}, $\nu: C_*^c(\mathbb{S}^n) \longrightarrow \mathfrak{g}$ is a curved twisting morphism if and only if 
\[
\gamma_\mathfrak{g}\left(\sum_{n \geq 0,~ n \neq 1}\frac{1}{n!}c_n(\nu(a_0), \cdots , \nu(a_0)) \right) + d_\mathfrak{g}(\nu(a_0)) = 0~,
\]
that is, $\nu(a_0)$ is a Maurer--Cartan element, and if 

\[
\gamma_\mathfrak{g}\left(\sum_{n \geq 0,~ n \neq 1}\frac{1}{n!}c_n(\nu(a_0), \cdots , \nu(a_{[n]})) \right) + d_\mathfrak{g}(\nu(a_{[n]})) = 0~,
\]
that is, $d^{\nu(a_0)}_\mathfrak{g}(\nu(a_{[n]}) = 0$, where $d^{\nu(a_0)}$ is the twisted differential by $\nu(a_0)$. Therefore we have a bijection
\[
\mathrm{Hom}_{\mathsf{sSets}_*}((\mathbb{S}^n,a_0), (\mathcal{R}(\mathfrak{g}),\alpha)) \cong \mathrm{Z}_n^{\alpha}(\mathfrak{g})~.
\]
Let $f_\nu: C_*^c(\mathbb{S}^n) \longrightarrow \widehat{\mathrm{B}}_\iota(\mathfrak{g})$ and $g_{\rho}: C_*^c(\mathbb{S}^n) \longrightarrow \widehat{\mathrm{B}}_\iota(\mathfrak{g})$ be two morphisms of $u\mathcal{CC}_\infty$-coalgebras which send $a_0$ to $\alpha$. The data of an homotopy between them amounts to the data of a morphism
\[
h: C_*^c(\mathbb{S}^n) \otimes  C_*^c(\Delta^1) \longrightarrow \widehat{\mathrm{B}}_\iota(\mathfrak{g})
\]
such that $h(- \otimes b_0) = f_{\nu}$ and $h(- \otimes b_1) = g_{\rho}$, where $C_*^c(\Delta^1) = \kk.b_0 \oplus \kk.b_1 \oplus \kk.b_{01}$, since $C_*^c(\Delta^1)$ is the interval object in the category of $u\mathcal{CC}_\infty$-coalgebra. The $u\mathcal{CC}_\infty$-coalgebra structure on $C_*^c(\mathbb{S}^n) \otimes  C_*^c(\Delta^1)$ is given by the tensor product structure of Corollary \ref{cor: tensor product of uCC coalgebras}. One can check that the data of such of a morphism $h$ is equivalent to the data of an element $\lambda$ of degree $n+1$ such that 
\[
d^{\alpha}(\lambda) = \nu(a_{[n]}) - \rho(a_{[n]})~.
\]
Therefore there is a bijection 
\[
\mathrm{Hom}_{\mathsf{sSets}_*}((\mathbb{S}^n,a_0), (\mathcal{R}(\mathfrak{g}),\alpha))/\sim_{\mathrm{hmt}} ~\cong \mathrm{H}_n^{\alpha}(\mathfrak{g})~.
\]
\end{proof}

\begin{Corollary}
Let $\alpha$ and $\beta$ be two gauge equivalent Maurer--Cartan elements of $\mathfrak{g}$. There is an isomorphism
\[
\mathrm{H}^\alpha_n(\mathfrak{g}) \cong \mathrm{H}^\beta_n(\mathfrak{g})~,
\]
for all $n \geq 1$. 
\end{Corollary}

\begin{proof}
It follows from Theorem \ref{thm: Berglund isomorphism} and Theorem \ref{thm: characterisation pi zero et jauges}.
\end{proof}

\begin{Remark}[Berglund map]
One can also prove this result by following the same ideas of \cite{Berglund15}, that is, by constructing an explicit map
\[
\begin{tikzcd}[column sep=4pc,row sep=0pc]
\mathscr{B}_n^{\alpha}: \mathrm{H}_n^{\alpha}(\mathfrak{g}) \arrow[r]
&\pi_n(\mathcal{R}(\mathfrak{g}),\alpha) \\
\left[u\right] \arrow[r, mapsto]
& \left[ f_{[u]}: C^c_*(\Delta^n) \longrightarrow \mathfrak{g} \right]~,
\end{tikzcd}
\]
where $f_{[u]}$ is the curved twisting morphism defined by $f_{[u]}(a_{[n]}) = u$ and $f_{[u]}(a_{i}) = \alpha$, for all $0 \leq i \leq n$, and which is zero on any other element $a_I$ of $C^c_*(\Delta^n)$.
\end{Remark}

\begin{Remark}
The group structure on $\pi_n(\mathcal{R}(\mathfrak{g}),\alpha)$ is induced by the pinch map 
\[
\mathrm{pinch}: \mathbb{S}^n \longrightarrow \mathbb{S}^n \vee \mathbb{S}^n~.
\]
which on cellular chains is given by $C^c_*(\mathrm{pinch})(a_{[n]}) = a_{[n]}^{(1)} + a_{[n]}^{(2)}$. Using this, it can be shown that for $n \geq 2$, there is an isomorphism of abelian groups
\[
\pi_n(\mathcal{R}(\mathfrak{g}),\alpha) \cong \mathrm{H}_n^\alpha(\mathfrak{g})~,
\]
where we consider the sum of homology classes on the right-hand side. 
\end{Remark}

\medskip

Let $h: \Lambda_1^2(\lambda_1,\lambda_2) \longrightarrow \mathfrak{g}$ be a $\Lambda_1^2$-horn in $\mathcal{R}(\mathfrak{g})$. Let $h(a_{01}) = \lambda_1$ and $h(a_{12}) = \lambda_2$. These are two degree one elements in $\mathfrak{g}$ which define a gauge equivalences $\lambda_1 \bullet \alpha = \alpha$ and $\lambda_2 \bullet \alpha = \alpha$. The \textit{Baker--Campbell--Hausdorff product relative to $\alpha$} of $\lambda_1$ and $\lambda_2$ is given by 

\[
\mathrm{BCH}^{\alpha}(\lambda_1,\lambda_2) \coloneqq \Gamma^0\left( \Lambda_1^2(\lambda_1,\lambda_2) \right)~.
\]

\begin{Remark}
The Baker--Campbell--Hausdorff product relative to $\alpha$ is not well-defined in general for any two degree one elements of $\mathfrak{g}$.
\end{Remark}

\begin{Corollary}
Let $\alpha$ be a Maurer--Cartan element of a curved absolute $\mathcal{L}_\infty$-algebra $\mathfrak{g}$. The Baker--Campbell--Hausdorff product relative to $\alpha$ defines a group structure with unit $[\alpha]$ on $\mathrm{H}_1^\alpha(\mathfrak{g})$, which is isomorphic to the group $\pi_1(\mathcal{R}(\mathfrak{g}),\alpha)$.
\end{Corollary}

\begin{proof}
This is straightforward from Theorem \ref{thm: Berglund isomorphism}.
\end{proof}

\section{Rational homotopy models}
In this section, we show that curved absolute $\mathcal{L}_\infty$-algebras are models for finite type nilpotent rational spaces without any pointed or connectivity assumptions. We do this by relating the adjunction $\mathcal{L} \dashv \mathrm{R}$ constructed in the previous section to the derived adjunction constructed by Bousfield and Guggenheim in \cite{BousfieldGugenheim} using Sullivan's functor \cite{Sullivan77}. We also construct \textit{small} models for mapping spaces without any hypothesis on the sources using the machinery developed so far.

\subsection{Comparing derived units of adjunctions}\label{Subsection: rational models}
From now on, we assume the base field $\kk$ to be the field of rational numbers $\mathbb{Q}$. Furthermore, we consider the category of simplicial sets endowed with the rational model structure constructed in \cite{Bousfield75}, where cofibrations are given by monomorphisms and where weak-equivalences are given by morphisms $f: X \longrightarrow Y$ such that $\mathrm{H}_*(f,\mathbb{Q}): \mathrm{H}_*(X,\mathbb{Q}) \longrightarrow \mathrm{H}_*(Y,\mathbb{Q})$ is an isomorphism.

\medskip

\textbf{Sullivan's adjunction.} One can use a dual version of Theorem \ref{thm: Kan seminal result} in order to induce the following contravariant adjunction 
\[
\begin{tikzcd}[column sep=5pc,row sep=2.5pc]
\mathsf{sSet} \arrow[r, shift left=1.5ex, "\mathrm{A}_{\mathrm{PL}}(-)"{name=A}]
&\mathsf{dg}~u\mathcal{C}om\textsf{-}\mathsf{alg}^{\mathsf{op}}~, \arrow[l, shift left=.75ex, "\langle - \rangle"{name=C}] \arrow[phantom, from=A, to=C, , "\dashv" rotate=-90]
\end{tikzcd}
\]
where $\mathrm{A}_{\mathrm{PL}}(-)$ is the piece-linear differential forms functor obtained by taking the left Kan extension of $\Omega^\bullet$. Its right adjoint functor, called the geometrical realization functor, is given by 
\[
\langle - \rangle_\bullet \coloneqq \mathrm{Hom}_{\mathsf{dg}~u\mathcal{C}om\textsf{-}\mathsf{alg}}(-, \Omega^\bullet)~.
\] 
This adjunction is in fact a Quillen adjunction when one considers the standard model structure on dg unital commutative algebras, where weak-equivalences are given by quasi-isomorphisms and fibrations by degree-wise epimorphisms.

\medskip

The quasi-isomorphism of dg operads $\varepsilon: \Omega \mathrm{B}^{\mathsf{s.a}}u\mathcal{C}om \qi u\mathcal{C}om$ induces the following Quillen equivalence 
\[
\begin{tikzcd}[column sep=5pc,row sep=2.5pc]
u\mathcal{CC}_\infty \textsf{-}\mathsf{alg} \arrow[r, shift left=1.5ex, "\mathrm{Ind}_\varepsilon"{name=A}]
&\mathsf{dg}~u\mathcal{C}om \textsf{-}\mathsf{alg} ~, \arrow[l, shift left=.75ex, "\mathrm{Res}_\varepsilon "{name=C}] \arrow[phantom, from=A, to=C, , "\dashv" rotate=-90]
\end{tikzcd}
\]
where $\mathrm{Res}_\varepsilon$ is fully faithful. Thus one obtains the following commutative triangle of Quillen adjunctions 
\[
\begin{tikzcd}[column sep=5pc,row sep=2.5pc]
&\hspace{2pc} \left(\mathsf{dg}~u\mathcal{C}om \textsf{-}\mathsf{alg}\right)^{\mathsf{op}} \arrow[dd, shift left=1.1ex, "\mathrm{Res}_\varepsilon^{\mathsf{op}}"{name=F}] \arrow[ld, shift left=.75ex, "\langle - \rangle"{name=C}]\\
\mathsf{sSet}  \arrow[ru, shift left=1.5ex, "\mathrm{A}_{\mathrm{PL}}(-)"{name=A}]  \arrow[rd, shift left=1ex, "\mathrm{C}_{\mathrm{PL}}(-)"{name=B}] \arrow[phantom, from=A, to=C, , "\dashv" rotate=-70]
& \\
&\hspace{1.5pc} \left(u\mathcal{CC}_\infty \textsf{-}\mathsf{alg}\right)^{\mathsf{op}} ~, \arrow[uu, shift left=.75ex, "\mathrm{Ind}_\varepsilon^{\mathsf{op}}"{name=U}] \arrow[lu, shift left=.75ex, "\langle - \rangle_\infty"{name=D}] \arrow[phantom, from=B, to=D, , "\dashv" rotate=-110] \arrow[phantom, from=F, to=U, , "\dashv" rotate=-180]
\end{tikzcd}
\]
by defining $\mathrm{C}_{\mathrm{PL}}(-) \coloneqq \mathrm{Res}_\varepsilon^{\mathsf{op}} \circ \mathrm{C}_{\mathrm{PL}}(-)$ and $\langle - \rangle_\infty \coloneqq \langle - \rangle \circ \mathrm{Ind}_\varepsilon^{\mathsf{op}}$.

\begin{lemma}\label{lemma: nieme equivalence}
Let $X$ be a simplicial set. There is a natural weak-equivalence of simplicial sets
\[
\mathbb{R}\langle \mathrm{C}_{\mathrm{PL}}(X) \rangle_\infty \simeq \mathbb{R}\langle \mathrm{A}_{\mathrm{PL}}(X)\rangle~.
\]
\end{lemma}

\begin{proof}
The adjunction induced by $\varepsilon: \Omega \mathrm{B}^{\mathsf{s.a}}u\mathcal{C}om \qi u\mathcal{C}om$ is a Quillen equivalence.
\end{proof}

\textbf{Cellular cochains functor.}
Consider again the simplicial $u\mathcal{CC}_\infty$-algebra $C_c^*(\Delta^\bullet)$ of Lemma \ref{lemma: simplicial CC infinity} given by the cellular cochains on the standard simplices together with their transferred $u\mathcal{CC}_\infty$-algebra structure. It induces a Quillen adjunction

\[
\begin{tikzcd}[column sep=5pc,row sep=2.5pc]
\mathsf{sSet} \arrow[r, shift left=1.5ex, "C_c^*(-)"{name=A}]
&u\mathcal{CC}_\infty \textsf{-}\mathsf{alg}^{\mathsf{op}}~, \arrow[l, shift left=.75ex, "\mathcal{S}"{name=C}] \arrow[phantom, from=A, to=C, , "\dashv" rotate=-90]
\end{tikzcd}
\]

by applying again Theorem \ref{thm: Kan seminal result}, where $C_c^*(-)$ is the cellular cochain functor endowed with a $u\mathcal{CC}_\infty$-algebra structure, and where the right adjoint functor $\mathcal{S}$ is given by 

\[
\mathcal{S}(-) = \mathrm{Hom}_{u\mathcal{CC}_\infty \textsf{-}\mathsf{alg}}(-,C_c^*(\Delta^\bullet))~.
\]

Our goal is to compare the derived unit of this adjunction with the derived unit of the adjunction obtained by extending Sullivan's functor to the category of $u\mathcal{CC}_\infty$-algebras. In order, to do so, we need to be able to construct cofibrant resolutions. The canonical curved twisting morphism $\iota: \mathrm{B}^{\mathsf{s.a}}u\mathcal{C}om \longrightarrow \Omega \mathrm{B}^{\mathsf{s.a}}u\mathcal{C}om$ induces a Quillen equivalence
\[
\begin{tikzcd}[column sep=5pc,row sep=2.5pc]
\mathsf{curv}~ \mathrm{B}^{\mathsf{s.a}}u\mathcal{C}om \textsf{-}\mathsf{coalg} \arrow[r, shift left=1.5ex, "\Omega_\iota"{name=A}]
&u\mathcal{CC}_\infty \textsf{-}\mathsf{alg} ~, \arrow[l, shift left=.75ex, "\mathrm{B}_\iota"{name=C}] \arrow[phantom, from=A, to=C, , "\dashv" rotate=-90]
\end{tikzcd}
\]
where the model structure considered on the left hand side is the one obtained by transfer along this adjunction, see \cite{grignou2019} for more details. An $\infty$-morphism of $u\mathcal{CC}_\infty$-algebras $f: A \rightsquigarrow B$ is the data of a morphism 
\[
f: \mathrm{B}_\iota A \longrightarrow \mathrm{B}_\iota B
\]
of curved $\mathrm{B}^{\mathsf{s.a}}u\mathcal{C}om$-coalgebras. It is an $\infty$-quasi-isomorphism if the term $f_{\mathrm{id}}: A \longrightarrow B$ is a quasi-isomorphism. This is equivalent to $f$ being a weak-equivalence in the transferred model structure.

\begin{lemma}
Let $X$ be a simplicial set. There is a pair of inverse $\infty$-quasi-isomorphism of $u\mathcal{CC}_\infty$-algebras
\[
(p_\infty)_X: \mathrm{C}_{\mathrm{PL}}(X) \rightsquigarrow C_c^*(X)~,~(i_\infty)_X : C_c^*(X) \rightsquigarrow \mathrm{C}_{\mathrm{PL}}(X)~,
\]
which are natural in $X$.
\end{lemma}

\begin{proof}
Let $X$ be a simplicial set. Since Dupont's contraction is compatible with the simplicial structures, it induces a homotopy retract 

\[
\begin{tikzcd}[column sep=5pc,row sep=3pc]
\mathrm{C}_{\mathrm{PL}}(X) \arrow[r, shift left=1.1ex, "p_X"{name=F}] \arrow[loop left]{l}{h_X}
& C_c^*(X)~. \arrow[l, shift left=.75ex, "i_X"{name=U}]
\end{tikzcd}
\]
One can show, using analogue arguments to the proof of \cite[Proposition 7.12]{robertnicoud2020higher}, that the transferred $u\mathcal{CC}_\infty$-algebra structure from this contraction onto $C_c^*(X)$ is equal to the $u\mathcal{CC}_\infty$-algebra structure obtained by considering the left Kan extension of $C_c^*(\Delta^\bullet)$. Therefore the morphisms of dg modules $p_X$ and $i_X$ can be extended to two (homotopy) inverse $\infty$-quasi-isomorphisms $(p_\infty)_X$ and $(i_\infty)_X$, which are both natural in $X$, see \cite[Remark 7.14]{robertnicoud2020higher}.
\end{proof}

\begin{Proposition}\label{prop: equivalence cochains CC with Apl}
Let $X$ be a simplicial set. There is a natural weak-equivalence of simplicial sets
\[
\mathbb{R}\langle \mathrm{C}_{\mathrm{PL}}(X) \rangle_\infty \simeq \mathbb{R}\mathcal{S}(C_c^*(X))~.
\]
\end{Proposition}

\begin{proof}
The above lemma implies that there is a natural weak-equivalence of curved $\mathrm{B}^{\mathsf{s.a}}u\mathcal{C}om$-coalgebras
\[
(p_\infty)_X: \mathrm{B}_\iota \mathrm{C}_{\mathrm{PL}}(X) \qi \mathrm{B}_\iota C_c^*(X)~.
\]
This in turn implies that there is a natural quasi-isomorphism 
\[
\Omega_\iota((p_\infty)_X) : \Omega_\iota\mathrm{B}_\iota \mathrm{C}_{\mathrm{PL}}(X) \qi \Omega_\iota\mathrm{B}_\iota C_c^*(X)
\]
of $u\mathcal{CC}_\infty$-algebras. Therefore we have 

\begin{align*}
\mathbb{R}\langle \mathrm{C}_{\mathrm{PL}}(X)\rangle_\infty &\simeq \mathrm{Hom}_{u\mathcal{C}om\text{-}\mathsf{alg}}(\mathrm{Ind}_\varepsilon \Omega_\iota\mathrm{B}_\iota \mathrm{C}_{\mathrm{PL}}(X), \Omega_\bullet)\\
& \cong \mathrm{Hom}_{u\mathcal{CC}_\infty\text{-}\mathsf{alg}}(\Omega_\iota\mathrm{B}_\iota \mathrm{C}_{\mathrm{PL}}(X), \mathrm{Res}_\varepsilon \Omega_\bullet) \\
&\simeq \mathrm{Hom}_{u\mathcal{CC}_\infty\text{-}\mathsf{alg}}(\Omega_\iota\mathrm{B}_\iota C_c^*(X), C_c^*(\Delta^\bullet)) \\
&\simeq \mathbb{R}\mathcal{S}(C_c^*(X))~,\\
\end{align*}

where we used the existence of an $\infty$-quasi-isomorphism of $u\mathcal{CC}_\infty$-algebras between $\mathrm{Res}_\varepsilon \Omega_\bullet$ and $C_c^*(\Delta^\bullet)$, given again by the Dupont contraction. Notice that the intermediate equivalences are obtained using the fact that these are hom-spaces between a cofibrant and a fibrant object, thus stable by weak-equivalences.
\end{proof}

\textbf{Finite type nilpotent spaces.} We recall what finite type nilpotent spaces are and we state the main theorem of this subsection.

\begin{Definition}[Finite type simplicial set]
Let $X$ be a simplicial set. It is said to be of \textit{finite type} if the homology groups $\mathrm{H}_n(X,\mathbb{Z})$ are finitely generated for all $n \geq 0$.
\end{Definition}

\begin{Remark}
In particular, $X$ has a finitely many connected components since $\mathrm{H}_0(X,\mathbb{Z})$ is finite dimensional.
\end{Remark}

\begin{Definition}[Nilpotent simplicial set]
Let $X$ be a simplicial set. It is said to be \textit{nilpotent} if for every $0$-simplex $\alpha$ it satisfies the following conditions:

\begin{enumerate}
\item The group $\pi_1(X,\alpha)$ is nilpotent.

\medskip

\item The $\pi_1(X,\alpha)$-module $\pi_n(X,\alpha)$-module is a nilpotent $\pi_1(X,\alpha)$-module.
\end{enumerate}
\end{Definition}

\begin{Proposition}\label{prop: n-ieme equivalence}
Let $X$ be a finite type simplicial set. There is a weak-equivalence of simplicial sets

\[
\mathcal{R}\mathcal{L}(X) \simeq \mathbb{R}\mathcal{S}(C_c^*(X))~,
\]

which is natural on the subcategory of finite type simplicial sets.
\end{Proposition}

\begin{proof}
Recall from the proof above that 
\[
\mathbb{R}\mathcal{S}(C_c^*(X)) \simeq \mathrm{Hom}_{u\mathcal{CC}_\infty}(\Omega_\iota\mathrm{B}_\iota C_c^*(X), C_c^*(\Delta^\bullet))~.
\]

The following square of Quillen adjunctions
\[
\begin{tikzcd}[column sep=5pc,row sep=5pc]
u\mathcal{CC}_\infty\text{-}\mathsf{alg}^{\mathsf{op}} \arrow[r,"\mathrm{B}_\iota^{\mathsf{op}}"{name=B},shift left=1.1ex] \arrow[d,"(-)^\circ "{name=SD},shift left=1.1ex ]
&\mathsf{curv}~\mathrm{B}^{\mathsf{s.a}}u\mathcal{C}om\text{-}\mathsf{coalg}^{\mathsf{op}} \arrow[d,"(-)^*"{name=LDC},shift left=1.1ex ] \arrow[l,"\Omega_\iota^{\mathsf{op}}"{name=C},,shift left=1.1ex]  \\
u\mathcal{CC}_\infty\text{-}\mathsf{coalg} \arrow[r,"\widehat{\Omega}_\iota "{name=CC},shift left=1.1ex]  \arrow[u,"(-)^*"{name=LD},shift left=1.1ex ]
&\mathsf{curv}~\mathrm{B}^{\mathsf{s.a}}u\mathcal{C}om\text{-}\mathsf{alg}^{\mathsf{comp}}~, \arrow[l,"\widehat{\mathrm{B}}_\iota"{name=CB},shift left=1.1ex] \arrow[u,"(-)^\vee"{name=TD},shift left=1.1ex] \arrow[phantom, from=SD, to=LD, , "\dashv" rotate=0] \arrow[phantom, from=C, to=B, , "\dashv" rotate=-90]\arrow[phantom, from=TD, to=LDC, , "\dashv" rotate=0] \arrow[phantom, from=CC, to=CB, , "\dashv" rotate=-90]
\end{tikzcd}
\] 

commutes, by \cite[Theorem 2.22]{lucio2022contra} applied to this situation. Since $X$ is a finite type simplicial set, the homology of $\Omega_\iota\mathrm{B}_\iota C_c^*(X)$ is degree-wise finite dimensional and bounded below. Therefore the generalized Sweedler dual $(-)^\circ$ is homotopically fully-faithful by \cite[Theorem 2.25]{lucio2022contra}, and we get 

\[
\mathrm{Hom}_{u\mathcal{CC}_\infty\text{-}\mathsf{alg}}(\Omega_\iota\mathrm{B}_\iota C_c^*(X), C_c^*(\Delta^\bullet)) \simeq \mathrm{Hom}_{u\mathcal{CC}_\infty\text{-}\mathsf{coalg}}(C_*^c(\Delta^\bullet), \widehat{\mathrm{B}}_\iota \widehat{\Omega}_\iota C_*^c(X) )~,
\]

which concludes the proof.
\end{proof}

\begin{Corollary}\label{cor: equivalences de foncteurs derives Apl et RL}
Let $X$ be a finite type simplicial set.  There is a weak-equivalence of simplicial sets

\[
\mathbb{R}\langle \mathrm{A}_{\mathrm{PL}}(X)\rangle \simeq \mathcal{R}\mathcal{L}(X)~,
\]

which is natural in the subcategory of finite type simplicial sets.
\end{Corollary}

\begin{proof}
This is a direct consequence of Lemma \ref{lemma: nieme equivalence} and Propositions \ref{prop: equivalence cochains CC with Apl} and \ref{prop: n-ieme equivalence}.
\end{proof}

Now we can transfer the known results about Sullivan's rational models along the above equivalence.

\begin{theorem}[{\cite[Theorem C]{MarklLazarev}}]\label{thm: Markl et Lazarev}
Let $X$ be a finite type nilpotent simplicial set. The unit 

\[
\eta_X: X \qi \mathbb{R}\langle \mathrm{A}_{\mathrm{PL}}(X)\rangle
\]

is a rational homotopy equivalence.
\end{theorem}

\begin{Remark}
The original construction of \cite{BousfieldGugenheim} only works for connected finite type nilpotent simplicial set and was subsequently extended to disconnected finite type nilpotent simplicial sets by Markl--Lazarev in \cite{MarklLazarev}.
\end{Remark}

\begin{theorem}\label{thm: modèles d'homotopie rationnel type fini}
Let $X$ be a finite type nilpotent simplicial set. The unit of the adjunction

\[
\eta_X: X \qi \mathcal{R}\mathcal{L}(X) 
\]

is a rational homotopy equivalence.
\end{theorem}

\begin{proof}
Follows directly from Corollary \ref{cor: equivalences de foncteurs derives Apl et RL} and Theorem \ref{thm: Markl et Lazarev}.
\end{proof}

\begin{Corollary}
Let $X$ be a pointed connected finite type simplicial set. The unit of the adjunction 
\[
\eta_X: X \longrightarrow \mathcal{RL}(X)
\]
is weakly equivalent to the $\mathbb{Q}$-completion of Bousfield-Kan. 
\end{Corollary}

\begin{proof}
This is implied by \cite[Theorem 12.2]{BousfieldGugenheim}, using the equivalence $\mathbb{R}\langle \mathrm{A}_{\mathrm{PL}}(X)\rangle \simeq \mathcal{RL}(X)$ of Corollary \ref{cor: equivalences de foncteurs derives Apl et RL}.
\end{proof}

\begin{Remark}
One can characterise the homotopical essential image of the cellular chains functor $C_c^*$ inside $u\mathcal{CC}_\infty$-coalgebras, see \cite[Theorem 3.37]{mathez}. 
\end{Remark}

\begin{Example}\label{Example: two point set and empty set}
Let us mention some examples of models in terms of curved absolute $\mathcal{L}_\infty$-algebras.

\medskip

\begin{enumerate}
\item A model for the constant simplicial set with one point is given by $\mathfrak{mc}^0$. A model for the constant simplicial  set with two points is given by the coproduct $\mathfrak{mc}^0 \amalg \mathfrak{mc}^0$. 

\medskip

\item A model for the empty set is given by $\kk.\varphi$, where $\varphi$ lies in degree $-1$. The only non-trivial structure map is the curvature $l_0: \kk \longrightarrow \kk.\varphi$, which maps $1$ to $\varphi$. In particular, it can be checked that it has no Maurer--Cartan element. 
\end{enumerate}
\end{Example}

\subsection{Minimal models}
In this subsection we show that complete curved absolute $\mathcal{L}_\infty$-algebras always admit a minimal resolutions, which are unique up to isomorphism. Furthermore, as a corollary of our constructions, we show that one can also recover the \textit{homology} of a space $X$ via the homology of the complete bar construction of its models.

\begin{Definition}[Minimal model]
Let $\mathfrak{g}$ be a complete curved absolute $\mathcal{L}_\infty$-algebra. A \textit{minimal model} $(V,\varphi_{d_V},\psi_V)$ amounts to the data of 

\begin{enumerate}
\item A graded module $V$ together with a map

\[
\varphi_{d_V}: V \longrightarrow \prod_{n \geq 0} \widehat{\Omega}^{\mathrm{s.a}}\ucom^*(n)^{(\geq 1)} ~\widehat{\otimes}_{\mathbb{S}_n} ~ V^{\otimes n}~
\]

which lands on elements of weight greater or equal to $1$, such that the induced derivation $d_V$ on the free pdg $\widehat{\Omega}^{\mathrm{s.a}}\ucom^*$-algebra on $V$ satisfies the curvature equation

\[
d_V^2 = l_2 \circ_1 l_0~. 
\]
\medskip

\item A weak-equivalence of complete curved absolute $\mathcal{L}_\infty$-algebras 

\[
\psi_V: \left( \widehat{\Omega}^{\mathrm{s.a}}\ucom^*(n) ~\widehat{\otimes}_{\mathbb{S}_n} ~ V^{\otimes n},d_V \right) \qi \mathfrak{g}~.
\]
\end{enumerate}
\end{Definition}

\begin{Remark}
Given a graded module $V$, the data of the derivation $d_V$ amounts to a $u\mathcal{CC}_\infty$-coalgebra structure on $V$. Thus the data $(V,\varphi_{d_V})$ above amounts to the data of a \textit{minimal} $u\mathcal{CC}_\infty$-coalgebra, that is, a $u\mathcal{CC}_\infty$-coalgebra where the underlying differential is zero.
\end{Remark}

\begin{Proposition}
Let $\mathfrak{g}$ be a complete curved absolute $\mathcal{L}_\infty$-algebra and let $(V,\varphi_{d_V},\psi_V)$ and $(W,\varphi_{d_W},\psi_W)$ be two minimal models of $\mathfrak{g}$. Then there is an isomorphism of graded modules $V \cong W$. 
\end{Proposition}

\begin{proof}
If $(V,\varphi_{d_V},\psi_V)$ and $(W,\varphi_{d_W},\psi_W)$ are two minimal models of $\mathfrak{g}$, then, by definition, there exists a zig-zag of weak-equivalences of curved absolute $\mathcal{L}_\infty$-algebras between them. Since both objects are fibrant-cofibrant, this gives a single weak-equivalence between them

\[
\left( \widehat{\Omega}^{\mathrm{s.a}}\ucom^*(n) ~\widehat{\otimes}_{\mathbb{S}_n} ~ V^{\otimes n},d_V \right) \qi \left( \widehat{\Omega}^{\mathrm{s.a}}\ucom^*(n) ~\widehat{\otimes}_{\mathbb{S}_n} ~ W^{\otimes n},d_W\right)~.
\]

Therefore there is a quasi-isomorphism between the graded modules $V$ and $W$, which implies they are isomorphic.
\end{proof}

\begin{Proposition}\label{prop; generateurs du model minimal, homologie de la Bar}
Let $\mathfrak{g}$ be a complete curved absolute $\mathcal{L}_\infty$-algebra. Then it admits a minimal model, where the graded module of generators is given by
\[
V \cong \mathrm{H}_*\left(\widehat{\mathrm{B}}_\iota \mathfrak{g} \right)~,
\]
that is, the homology of its complete bar construction.
\end{Proposition}

\begin{proof}
Since $\kk$ is a field of characteristic $0$, one can always choose a contraction between $\widehat{\mathrm{B}}_\iota \mathfrak{g}$ and its homology in order to apply the Homotopy Transfer Theorem. The transferred $u\mathcal{CC}_\infty$-coalgebra structure provides us with the derivation by applying the complete cobar construction. The image of the inclusion
\[
i:  \mathrm{H}_*\left(\widehat{\mathrm{B}}_\iota \mathfrak{g} \right) \twoheadrightarrow\widehat{\mathrm{B}}_\iota \mathfrak{g}~,
\]

under the complete cobar construction allows us obtain the required weak-equivalence by pre-composing it with the counit $\epsilon_\mathfrak{g}: \widehat{\Omega}_\iota\widehat{\mathrm{B}}_\iota \mathfrak{g} \qi \mathfrak{g}~.$
\end{proof}

\begin{Remark}
The same minimal model constructions hold for absolute $\mathcal{L}_\infty$-algebras using  the complete bar construction $\widehat{\mathrm{B}}_\iota^\flat$ relative to $\mathcal{CC}_\infty$-coalgebras.
\end{Remark}

\begin{Definition}[Rational curved absolute $\mathcal{L}_\infty$-algebras] 
Let $\mathfrak{g}$ be a complete curved absolute $\mathcal{L}_\infty$-algebra. It is \textit{rational} if there exists a simplicial set $X$ and a zig-zag of weak-equivalences of curved absolute $\mathcal{L}_\infty$-algebras 
\[
\mathcal{L}(X) \lqi \cdot \qi \cdots \lqi \cdot \qi \mathfrak{g}~.
\]
\end{Definition}

\begin{Remark}
If $\mathfrak{g}$ is rational, then the homology of $\widehat{\mathrm{B}}_\iota(\mathfrak{g})$ is also concentrated in positive degrees. Indeed, $\mathcal{L}(X)$ is weakly equivalent to $\mathfrak{g}$ if and only if $C_*^c(X)$ is quasi-isomorphic to the complete bar construction of $\mathfrak{g}$.
\end{Remark}

\begin{Definition}[Rational models] 
Let $\mathfrak{g}$ be a complete curved absolute $\mathcal{L}_\infty$-algebra. It is a \textit{rational model} if it is rational for some simplicial set $X$ and if furthermore $\mathcal{R}(\mathfrak{g})$ is weakly equivalent to $X$. 
\end{Definition}

\begin{Proposition}
Let $\mathfrak{g}$ be a complete curved absolute $\mathcal{L}_\infty$-algebra which is rational for a simplicial set $X$. Then its minimal model is generated by $\mathrm{H}_*(X)$.
\end{Proposition}

\begin{proof}
In this particular case, we have that 
\[
\mathrm{H}_*\left(\widehat{\mathrm{B}}_\iota \mathfrak{g} \right) \cong \mathrm{H}_*(X)~,
\]
therefore by Proposition \ref{prop; generateurs du model minimal, homologie de la Bar}, the graded module $\mathrm{H}_*(X)$ is the generator of the minimal model. 
\end{proof}

\begin{Remark}
Notice that the situation here is Koszul dual to Sullivan's minimal models, where the minimal Sullivan model for $\mathrm{A}_{\mathrm{PL}}(X)$ of a simply connected space $X$ is generated by the linear dual of the homotopy groups $\pi_*(X)^*$. 
\end{Remark}

\textbf{Recovering the homology groups.} Let $\mathfrak{g}$ be curved absolute $\mathcal{L}_\infty$-algebra, let us try to understand the homotopy type of $\mathcal{R}(\mathfrak{g})$. The \textit{homotopy groups} of $\mathcal{R}(\mathfrak{g})$ can be fully described by the $\alpha$-homology groups $\mathrm{H}_*^\alpha(\mathfrak{g})$, where $\alpha$ runs over the set of Maurer--Cartan elements of $\mathfrak{g}$, as it was shown in Theorem \ref{thm: Berglund isomorphism}. In the case where $\mathfrak{g}$ is a rational model, then one can also recover the \textit{homology groups} of $\mathcal{R}(\mathfrak{g})$ using only the complete bar construction with respect to $\mathfrak{g}$. 

\begin{theorem}
Let $\mathfrak{g}$ be a curved absolute $\mathcal{L}_\infty$-algebra that is a rational model. The canonical morphism of $u\mathcal{CC}_\infty$-algebras 
\[
C_*^c(\mathcal{R}(\mathfrak{g})) \qi \widehat{\mathrm{B}}_\iota \mathfrak{g}
\]
is a quasi-isomorphism.
\end{theorem}

\begin{proof}
In this situation, $\mathcal{L}(\mathcal{R}(\mathfrak{g})) \qi \mathfrak{g}$ is a weak-equivalence of curved absolute $\mathcal{L}_\infty$-algebras, which is equivalent to the canonical morphism $C_*^c(\mathcal{R}(\mathfrak{g})) \qi \widehat{\mathrm{B}}_\iota \mathfrak{g} $ being a quasi-isomorphism of $u\mathcal{CC}_\infty$-coalgebras.
\end{proof}

\begin{Remark}
One can think of $\widehat{\mathrm{B}}_\iota \mathfrak{g}$ as a \textit{higher Chevalley-Eilenberg} complex adapted to the setting of curved absolute $\mathcal{L}_\infty$-algebras. 
\end{Remark}

\subsection{Models for mapping spaces}
In this section, we construct explicit rational models for mapping spaces, without any assumption on the source simplicial set. Furthermore, these models are relatively small. They can be constructed using the cellular chains of the source, if one wishes to preserve naturality, or even with the homology of the source, if one does not need naturality. 

\medskip

Recall that, if $X$ and $Y$ are simplicial sets, there is an explicit model for their mapping space given by 

\[
\mathrm{Map}(X,Y)_\bullet \coloneqq \mathrm{Hom}_{\mathsf{sSet}}(X \times \Delta^\bullet, Y)~,
\]
\vspace{0.25pc}

which forms a Kan complex when $Y$ is so.

\begin{lemma}\label{prop: les chaines sont lax monoidales en infini morphismes}
Let $X$ and $Y$ be two simplicial sets. There is an $\infty$-quasi-isomorphism 
\[
\psi_{X,Y}: C^c_*(X \times Y) \rightsquigarrow  C^c_*(X) \otimes C_*^c(Y)~,
\]
of $u\mathcal{CC}_\infty$-coalgebras which is natural in $X$ and $Y$.
\end{lemma}

\begin{proof}
For any simplicial sets $X,Y$, recall that there is a natural quasi-isomorphism 
\[
\kappa: \mathrm{A}_{\mathrm{PL}}(X) \otimes \mathrm{A}_{\mathrm{PL}}(Y) \qi \mathrm{A}_{\mathrm{PL}}(X \times Y)
\]
of dg $u\mathcal{C}om$-algebras. See \cite{bookRHT} for instance. This gives a quasi-isomorphism
\[
\mathrm{Res}_{\varepsilon}(\kappa): \mathrm{C}_{\mathrm{PL}}(X) \otimes \mathrm{C}_{\mathrm{PL}}(Y) \qi \mathrm{C}_{\mathrm{PL}}(X \times Y)
\]
of $u\mathcal{CC}_\infty$-algebras which is natural in $X,Y$, using the fact that the restriction functor $\mathrm{Res}_{\varepsilon}$ is strong monoidal and preserves all quasi-isomorphisms.

\medskip

Using the $\infty$-quasi-isomorphisms constructed in the proof of Proposition \ref{prop: equivalence cochains CC with Apl}, we construct 

\[
\begin{tikzcd}[column sep=4pc,row sep= 0pc]
C^*_c(X) \otimes C^*_c(Y) \arrow[r,"(i_\infty)_X \otimes (i_\infty)_Y ",rightsquigarrow]
&\mathrm{C}_{\mathrm{PL}}(X) \otimes \mathrm{C}_{\mathrm{PL}}(Y) \arrow[r,"\mathrm{Res}_{\varepsilon}(\kappa)"]
&\mathrm{C}_{\mathrm{PL}}(X \times Y) \arrow[r,"(p_\infty)_{X \times Y}",rightsquigarrow]
&C^*_c(X \times Y)~,
\end{tikzcd}
\]

where the tensor product of two $\infty$-morphisms is still an $\infty$-morphism, since the $2$-colored dg operad encoding $\infty$-morphisms of $u\mathcal{CC}_\infty$-algebras is cofibrant. This gives $\infty$-quasi-morphism of $u\mathcal{CC}_\infty$-algebras 

\[
C^*_c(X) \otimes C^*_c(Y) \rightsquigarrow C^*_c(X \times Y)~.
\]
\vspace{0.25pc}

which is natural in $X$ and $Y$. Equivalently, a weak-equivalence $\mathrm{B}_\iota(C^*_c(X) \otimes C^*_c(Y)) \qi \mathrm{B}_\iota(C^*_c(X \times Y))$. Now we take $X$ and $Y$ to be finite simplicial sets (finite total number of non-degenerate simplices). By applying the linear dual functor of \cite[Theorem 2.22]{lucio2022contra}, this gives a weak-equivalence

\[
\widehat{\Omega}_\iota(C^c_*(X) \otimes C_*^c(Y)) \qi \widehat{\Omega}_\iota(C_*^c(X \times Y))
\]
\vspace{0.25pc}

of complete curved absolute $\mathcal{L}_\infty$-algebras. Now suppose $X$ and $Y$ are arbitrary simplicial sets, one can write them as the filtered colimit of finite simplicial sets 
\[
X \cong \colim_\alpha X_\alpha \quad \text{and} \quad Y \cong \colim_\beta Y_\beta~.
\]
Now we have that 
\[
\widehat{\Omega}_\iota(C^c_*(X) \otimes C_*^c(Y)) \cong \colim_{\alpha,\beta} \widehat{\Omega}_\iota(C_*^c(X_{\alpha}) \otimes C_*^c(Y_{\beta}))~, 
\]
and 
\[
\widehat{\Omega}_\iota(C_*^c(X \times Y)) \cong \colim_{\alpha,\beta} \widehat{\Omega}_\iota(C_*^c(X_{\alpha} \times Y_{\beta}))~,
\]
since both bifunctors preserve filtered colimits in each variable. Both of this colimits are in fact homotopy colimits, hence there is a weak-equivalence 
\[
\widehat{\Omega}_\iota(C^c_*(X) \otimes C_*^c(Y)) \qi \widehat{\Omega}_\iota(C_*^c(X \times Y))~,
\]
which concludes the proof.
\end{proof}

\begin{theorem}\label{thm: vrai thm mapping spaces}
Let $\mathfrak{g}$ be a complete curved absolute $\mathcal{L}_\infty$-algebra and let $X$ be a simplicial set. There is a weak-equivalence of Kan complexes

\[
\mathrm{Map}(X, \mathcal{R}(\mathfrak{g})) \simeq \mathcal{R}\left(\mathrm{hom}(C^c_*(X),\mathfrak{g})\right)~,
\]
\vspace{0.25pc}

which is natural in $X$ and in $\mathfrak{g}$, where $\mathrm{hom}(C^c_*(X),\mathfrak{g})$ denotes the convolution curved absolute $\mathcal{L}_\infty$-algebra. Furthermore, it is possible to replace the cellular chains $C^c_*(X)$ by the homology $\mathrm{H}_*(X)$ to obtain a smaller model, meaning that there is a weak-equivalence of Kan complexes

\[
\mathrm{Map}(X, \mathcal{R}(\mathfrak{g})) \simeq \mathcal{R}\left(\mathrm{hom}(\mathrm{H}_*(X),\mathfrak{g})\right)~,
\]
\vspace{0.25pc}

which is now only natural in $\mathfrak{g}$. 
\end{theorem}

\begin{proof}
There is an isomorphism

\[
\mathrm{Map}(X,\mathcal{R}(\mathfrak{g}))_\bullet \coloneqq \mathrm{Hom}_{\mathsf{sSet}}(X \times \Delta^\bullet, \mathcal{R}(\mathfrak{g})) \cong \mathrm{Hom}_{\mathsf{sSet}}( C^c_*(X \times \Delta^\bullet), \widehat{\mathrm{B}}_\iota(\mathfrak{g}))~.
\]
\vspace{0.25pc}

We can pre-compose by the $\infty$-quasi-isomorphism of Lemma \ref{prop: les chaines sont lax monoidales en infini morphismes}, giving a weak-equivalence of simplicial sets

\[
\mathrm{Hom}_{u\mathcal{CC}_\infty}(C^c_*(X \times \Delta^\bullet), \widehat{\mathrm{B}}_\iota(\mathfrak{g})) \qi \mathrm{Hom}_{u\mathcal{CC}_\infty\text{-}\mathsf{cog}}( C^c_*(X)  \otimes C^c_*(\Delta^\bullet), \widehat{\mathrm{B}}_\iota(\mathfrak{g}))~,
\]
\vspace{0.25pc}

since both $C^c_*(X \times \Delta^\bullet)$ and $C^c_*(X)  \otimes C^c_*(\Delta^\bullet)$ are Reedy cofibrant. Let's compute this last simplicial set:

\begin{align*}
\mathrm{Hom}_{u\mathcal{CC}_\infty\text{-}\mathsf{cog}}( C^c_*(X)  \otimes C^c_*(\Delta^\bullet), \widehat{\mathrm{B}}_\iota(\mathfrak{g})) 
&\cong \mathrm{Hom}_{u\mathcal{CC}_\infty\text{-}\mathsf{cog}} \left( C^c_*(\Delta^\bullet), \left\{ C^c_*(X), \widehat{\mathrm{B}}_\iota(\mathfrak{g}) \right\} \right) \\
&\cong \mathrm{Hom}_{u\mathcal{CC}_\infty\text{-}\mathsf{cog}} \left( C^c_*(\Delta^\bullet), \widehat{\mathrm{B}}_\iota \left(\mathrm{hom}( C^c_*(X), \mathfrak{g}) \right) \right) \\
& \cong \mathrm{Hom}_{\mathsf{curv}~\mathsf{abs}~\mathcal{L}_\infty\text{-}\mathsf{alg}} \left( \widehat{\Omega}_\iota(C^c_*(\Delta^\bullet)), \mathrm{hom}( C^c_*(X), \mathfrak{g}) \right) \\
& \cong \mathcal{R}\left(\mathrm{hom}( C^c_*(X), \mathfrak{g}) \right)~.
\end{align*}
\vspace{0.1pc}

Upon choosing a contraction between $C^c_*(X)$ and its homology $\mathrm{H}_*(X)$, we can use the homotopy transfer theorem to obtain a $u\mathcal{CC}_\infty$-coalgebra structure on $\mathrm{H}_*(X)$ which is weakly equivalent to $C^c_*(X)$. This implies that $\mathrm{H}_*(X)  \otimes C^c_*(\Delta^\bullet)$ and $C^c_*(X)  \otimes C^c_*(\Delta^\bullet)$ are weakly equivalent as Reedy cofibrant objects, and thus we can replace the later by the former in the above computation.
\end{proof}

\begin{Corollary}\label{Cor: True mapping spaces.}
Let $X$ be a simplicial set and let $Y$ be a finite type nilpotent simplicial set. There is a weak-equivalence of simplicial sets 

\[
\mathrm{Map}(X, Y_\mathbb{Q}) \simeq \mathcal{R}\left(\mathrm{hom}(\mathrm{H}_*(X), \mathcal{L}(Y))\right)~,
\]
\vspace{0.1pc}

where $Y_\mathbb{Q}$ denotes the $\mathbb{Q}$-localization of $Y$. If we chose $\mathcal{RL}(Y)$ as a model for the $\mathbb{Q}$-localization of $Y$, given a map $f: X \longrightarrow \mathcal{RL}(Y)$, there are isomorphism

\[
\pi_n\left(\mathrm{Map}(X, \mathcal{RL}(Y)),f\right) \cong \mathrm{H}_n^{\alpha_f}\left(\mathrm{hom}(\mathrm{H}_*(X), \mathcal{L}(Y))\right)~,
\]
\vspace{0.1pc}

of groups for for $n \geq 1$. Here $\alpha_f$ is the curved twisting morphism associated to the following composition morphism of curved absolute $\mathcal{L}_\infty$-algebras 
\[
\begin{tikzcd}
\widehat{\Omega}_\iota(\mathrm{H}_*(X))) \arrow[r,"\iota_\infty"]
&\widehat{\Omega}_\iota(C^c_*(X)) \arrow[r,"f^\sharp"]
&\mathcal{L}(Y)~,
\end{tikzcd}
\]
where $f^\sharp$ is the adjoint map to $f$ and $\iota_\infty$ the $\infty$-quasi-isomorphism given by the homotopy transfer theorem. 
\end{Corollary}

\begin{proof}
There is a weak-equivalence of simplicial sets 
\[
Y_\mathbb{Q} \simeq \mathcal{RL}(Y)~.
\]
It induces equivalences
\[
\mathrm{Map}(X, Y_\mathbb{Q}) \simeq \mathrm{Map}(X, \mathcal{RL}(Y)) \simeq \mathcal{R}\left(\mathrm{hom}(\mathrm{H}_*(X), \mathcal{L}(Y))\right)~.
\]
\vspace{0.1pc}

Finally, the second statement follows directly from Theorem \ref{thm: Berglund isomorphism}. 
\end{proof}

\begin{Remark}
Notice the following points about Corollary \ref{Cor: True mapping spaces.}.

\begin{enumerate}
\medskip

\item There is no finiteness hypothesis on the simplicial set $X$. 

\medskip

\item So far, the models for mapping spaces in the literature for mapping spaces are given in terms of complete tensor products
\[
A_X ~\widehat{\otimes}~ \mathfrak{g}_Y~,
\]
of a dg unital commutative algebra model $A_X$ of $X$ and a complete $\mathcal{L}_\infty$-algebra model $\mathfrak{g}_Y$. Our model in comparison is much smaller as it only uses the homology of $X$, see for instance \cite{Berglund15, BuijMapping, LazarevMapping}.
\end{enumerate}
\end{Remark}

\begin{Corollary}
Let $X$ be a finite type nilpotent simplicial set. The space of rational endomorphisms of $X$ is weakly equivalent to 

\[
\mathrm{End}^h_\mathbb{Q}(X) \qi \mathcal{R}\left(\mathrm{hom}(\mathrm{H}_*(X), \mathcal{L}(X))\right)~.
\]
\end{Corollary}

\begin{proof}
This is immediate from the previous corollary. 
\end{proof}

\subsection{Models for pointed rational spaces}
In this subsection, we treat the case of pointed finite type nilpotent spaces, constructing rational models for them using complete absolute $\mathcal{L}_\infty$-algebras. 

\medskip

The simplicial $u\mathcal{C}om$-algebra $\Omega_\bullet$ admits an augmentation, given by the choice of a coherent base point in $\Delta^n$ for all $n \geq 0$. The augmentation ideal $\bar{\Omega}_\bullet$ is a simplicial $\mathcal{C}om$-algebra, and using augmented version of the Dupont contraction, we can transfer a simplicial $\mathcal{CC}_\infty$-algebra structure onto the reduced cellular cochains $\tilde{C}_c^*(\Delta^\bullet)$. We can dualize this object into a cosimplicial $\mathcal{CC}_\infty$-coalgebra and thus obtain the following triangle of commuting Quillen adjunction

\[
\begin{tikzcd}[column sep=5pc,row sep=2.5pc]
&\hspace{1pc}\mathcal{CC}_\infty \textsf{-}\mathsf{coalg} \arrow[dd, shift left=1.1ex, "\widehat{\Omega}_{\iota}^\flat"{name=F}] \arrow[ld, shift left=.75ex, "\overline{\mathcal{R}}_*"{name=C}]\\
\mathsf{sSet}_*  \arrow[ru, shift left=1.5ex, "\tilde{C}^c_*(-)"{name=A}]  \arrow[rd, shift left=1ex, "\mathcal{L}_*"{name=B}] \arrow[phantom, from=A, to=C, , "\dashv" rotate=-70]
& \\
&\hspace{3pc}\mathsf{abs}~\mathcal{L}_\infty\textsf{-}\mathsf{alg}^{\mathsf{comp}}~, \arrow[uu, shift left=.75ex, "\widehat{\mathrm{B}}_{\iota}^\flat"{name=U}] \arrow[lu, shift left=.75ex, "\mathcal{R}_*"{name=D}] \arrow[phantom, from=B, to=D, , "\dashv" rotate=-110] \arrow[phantom, from=F, to=U, , "\dashv" rotate=-180]
\end{tikzcd}
\]

where the model structure on $\mathcal{CC}_\infty$-coalgebras is transferred from dg modules and the model structure on complete absolute $\mathcal{L}_\infty$-algebras is transferred along the complete bar-cobar adjunction. 

\medskip

On the one hand, the functor $\mathcal{R}_*$ gives a \textit{pointed} version of the integration functor $\mathcal{R}$ for absolute $\mathcal{L}_\infty$-algebras. On the other, the functor $\mathcal{L}_*$ gives rational models for pointed connected finite type nilpotent spaces. 

\medskip

\textbf{Comparison between the pointed and the non-pointed cases.} To turn a curved $\mathcal{L}_\infty$-algebra into a $\mathcal{L}_\infty$-algebra, one can pick a Maurer--Cartan element and twist the original algebra by this element. This algebraic procedure, called the \textit{twisting procedure}, corresponds exactly to picking a base-point.  A generalization of this twisting procedure should also work at the \textit{absolute} level, see Remark \ref{Remark: twisting procedure}. In the meantime, another way to relate the pointed and non-pointed versions is by choosing a \textit{coaugmentation} at the coalgebraic level. This turns out to be quite easy, as the type of $u\mathcal{CC}_\infty$-coalgebras that appear via the cellular chain functor $C^c_*(-)$ are \textit{stricly counital}. We also refer to \cite[Chapter 3, Sections 1 and 3]{mathez} for some of these results. 

\begin{Definition}[Strictly counital $u\mathcal{CC}_\infty$-coalgebras]
A $u\mathcal{CC}_\infty$-coalgebra $C$ is \textit{strictly counital} if the elementary decomposition maps 
\[
\Delta_{\tau}: C \longrightarrow C^{\otimes n}
\]
are the zero morphism for any $\tau$ in $\mathrm{CRT}_n^\omega$ with at least one cork, except for $\epsilon: C \longrightarrow \kk$, which is the decomposition associated to the single cork.
\end{Definition}

The category of strictly counital $u\mathcal{CC}_\infty$-coalgebras is a coreflective subcategory of the category of $u\mathcal{CC}_\infty$-coalgebras. In general, a $u\mathcal{CC}_\infty$-coalgebra is \textit{coaugmented} if the terminal morphism $\Delta_C^0: C \longrightarrow \mathscr{C}(u\mathcal{CC}_\infty)(0)$ admits a section. However, if $C$ is strictly counital, the terminal morphism $\Delta_C^0$ factors through the quasi-isomorphism $k: \kk \qi \mathscr{C}(u\mathcal{CC}_\infty)(0)$. Thus the data of a coaugmentation for a strictly counital $u\mathcal{CC}_\infty$-coalgebra $C$ is equivalent to the data of a \textit{group-like element}, that is, of a morphism $\alpha: \kk \longrightarrow C$ of $u\mathcal{CC}_\infty$-coalgebras. We will refer to coaugmented strictly counital $u\mathcal{CC}_\infty$-coalgebras as \textit{pointed} strictly counital $u\mathcal{CC}_\infty$-coalgebras from now on. Finally, notice that there is an adjunction 
\[
\begin{tikzcd}[column sep=7pc,row sep=3pc]
\mathsf{Strict}~u\mathcal{CC}_\infty\textsf{-}\mathsf{coalg}_{\bullet}    \arrow[r, shift left=.75ex, "\mathrm{Ker}(\epsilon)"{name=U}]  
&\mathcal{CC}_\infty\textsf{-}\mathsf{coalg}~, \arrow[l, shift left=1.1ex, "(-) \oplus \kk"{name=F}] 
\arrow[phantom, from=F, to=U, , "\dashv" rotate=-90]
\end{tikzcd}
\]
between the category of pointed (coaugmented) strictly counital $u\mathcal{CC}_\infty$-coalgebras and the category of $\mathcal{CC}_\infty$-coalgebras, where the functors are given by 

\begin{enumerate}
\medskip

\item the left adjoint $\mathrm{Ker}(\epsilon)$ is given by considering the $\mathcal{CC}_\infty$-coalgebra structure on the coaugmentation coideal of a coaugmented strictly counital $u\mathcal{CC}_\infty$-coalgebra;

\medskip 

\item the right adjoint $(-) \oplus \kk$ adds cofreely a strict counit to a $\mathcal{CC}_\infty$-coalgebra.

\medskip 
\end{enumerate}

Furthermore, it is straightforward to check that this adjunction is an equivalence of categories.

\begin{Proposition}\label{prop: C(X) is strictly counital}
Let $X$ be a simplicial set, then $C_*^c(X)$ is a strictly counital $u\mathcal{CC}_\infty$-coalgebra.
\end{Proposition}

\begin{proof}
By Proposition \ref{prop: formulas for the HTT}, we know that it is true for the cellular chains of the standard $n$-simplex $C_*^c(\Delta^\bullet)$. Let $X$ be a simplicial set, it can be written as the following colimit
\[
X \cong \colim_{\mathrm{E}(X)} \Delta^\bullet~,
\]
indexed by $\mathrm{E}(X)$, the category of elements of $X$. Consequently, we have that 
\[
C_*^c(X) \cong \colim_{\mathrm{E}(X)} C_*^c(\Delta^\bullet)~.
\]
Strictly counital $u\mathcal{CC}_\infty$-coalgebras form a coreflective subcategory of $u\mathcal{CC}_\infty$-coalgebras, therefore they are closed under colimits. Meaning that $C_*^c(X)$ is a strictly counital $u\mathcal{CC}_\infty$-coalgebra.
\end{proof}

\begin{Proposition}\label{Prop: commutativity of the pointing adjunctions}
The following square of Quillen adjunctions
\[
\begin{tikzcd}[column sep=5pc,row sep=5pc]
\mathsf{sSet}_* \arrow[r,"\mathcal{L}_*"{name=B},shift left=1.1ex] \arrow[d,"\mathrm{U}"{name=SD},shift left=1.1ex ]
&\mathsf{abs}~\mathcal{L}_\infty\text{-}\mathsf{alg}^{\mathsf{comp}} \arrow[d,"\mathrm{U}"{name=LDC},shift left=1.1ex ] \arrow[l,"\mathcal{R}_*"{name=C},,shift left=1.1ex]  \\
\mathsf{sSet} \arrow[r,"\mathcal{L} "{name=CC},shift left=1.1ex]  \arrow[u,"(-) \sqcup \{*\}"{name=LD},shift left=1.1ex ]
&\mathsf{curv}~\mathsf{abs}~\mathcal{L}_\infty\text{-}\mathsf{alg}^{\mathsf{comp}}~, \arrow[l,"\mathcal{R}"{name=CB},shift left=1.1ex] \arrow[u,"(-)_*"{name=TD},shift left=1.1ex] \arrow[phantom, from=SD, to=LD, , "\dashv" rotate=0] \arrow[phantom, from=C, to=B, , "\dashv" rotate=-90]\arrow[phantom, from=TD, to=LDC, , "\dashv" rotate=0] \arrow[phantom, from=CC, to=CB, , "\dashv" rotate=-90]
\end{tikzcd}
\] 

is commutative in the following sense: the left adjoints from the bottom-left to the top-right corner are naturally isomorphic.
\end{Proposition}

\begin{proof}
Let $X$ be a simplicial set. On the one hand, we have that:
\[
\mathcal{L}_*(X \sqcup \{*\}) \cong \widehat{\Omega}_\iota^\flat(\mathrm{Ker}(\epsilon)(C_*^c(X \sqcup \{*\}))) \cong  \widehat{\Omega}_\iota^\flat(\mathrm{Ker}(\epsilon)(C_*^c(X) \oplus \kk)) \cong  \widehat{\Omega}_\iota^\flat(\mathrm{Res}(C_*^c(X)))~,
\]
where $\mathrm{Res}(C_*^c(X))$ is the $\mathcal{CC}_\infty$-coalgebra obtained by restricting the $u\mathcal{CC}_\infty$-coalgebra structure on $C_*^c(X)$. And on the other hand, we have that 
\[
(\mathcal{L}(X))_* \cong \widehat{\Omega}_\iota^\flat(\mathrm{Res}(C_*^c(X)))~,
\]
by the commutative of the square in Proposition \ref{Prop: fully faithful inclusion of non-curved into curved}. 
\end{proof}

\begin{Corollary}\label{Corollary: pointed rational models}
Let $X$ be a pointed finite type nilpotent simplicial set. The unit
\[
\eta_X: X \qi \mathcal{R}_* \mathcal{L}_*(X)
\]
is a rational homotopy equivalence.
\end{Corollary}

\begin{proof}
Let $X$ be a pointed finite type nilpotent simplicial set. The unit of adjunction 
\[
\eta_X: X \longrightarrow \mathcal{R}_* \mathcal{L}_*(X)
\]
is a rational homotopy equivalence if and only if its image 
\[
\mathrm{U}(\eta_X): \mathrm{U}(X) \longrightarrow \mathrm{U}\mathcal{R}_* \mathcal{L}_*(X)
\]
via the forgetful functor from pointed simplicial sets to simplicial sets is a rational homotopy equivalence. There a natural isomorphism
\[
\mathrm{U}\mathcal{R}_* \mathcal{L}_* \cong \mathcal{R}\mathrm{U}(-)_*\mathcal{L}~,
\]
by Proposition \ref{Prop: commutativity of the pointing adjunctions}. We conclude using that the unit of the $(-)_* \dashv \mathrm{U}$ is a weak-equivalence by Proposition \ref{Prop: fully faithful inclusion of non-curved into curved} and using that the unit of the $\mathcal{L} \dashv \mathcal{R}$ adjunction is a rational homotopy equivalence by Theorem \ref{thm: modèles d'homotopie rationnel type fini}.
\end{proof}

\begin{Remark}
Notice that, like in the work of \cite{Buijandco}, there is no connectivity assumption in Corollary \ref{Corollary: pointed rational models}. 
\end{Remark}

\begin{Remark}
Let $X$ be a finite type nilpotent simplicial set and let 
\[
X \cong \underset{\alpha \in \pi_0(X)}{\coprod} X^{\alpha}
\]
be its decomposition into connected components, where $\alpha: \{*\} \longrightarrow X$ are representatives of the equivalences classes in $\pi_0(X)$. It follows from the discussion above that there are rational homotopy equivalences:
\[
X \cong \underset{\alpha \in \pi_0(X)}{\coprod} X^{\alpha} \simeq \underset{\alpha \in \pi_0(X)}{\coprod} \mathcal{R}_* \mathcal{L}_*(X^{\alpha}) \simeq \mathcal{R}\mathcal{L}(X)~.
\]
\end{Remark}

\subsection{An intrinsic homotopy theory and comparisons with other models}
The model category structures considered so far on complete (curved) absolute $\mathcal{L}_\infty$-algebras are transferred from coalgebras via the complete bar-cobar adjunctions, which are Quillen equivalences. Thus, complete (curved) absolute $\mathcal{L}_\infty$-algebras, endowed with these models structures, give another presentation of the same $\infty$-category as that of (counital) cocommutative coalgebras. In this subsection, we show that there are right Bousfield localizations of these model category structures which still encode finite type nilpotent rational homotopy types. These have the advantage of endowing complete (curved) absolute $\mathcal{L}_\infty$-algebras with an \textit{intrinsic} homotopy theory, and allow us to compare our models with other Lie type models present in the literature. 

\medskip

\textbf{The pointed connected case.} The category of complete absolute $\mathcal{L}_\infty$-algebras admits a model category structure determined by the following classes of maps

\begin{enumerate}
\medskip

\item the set of fibrations is given by degree-wise surjections;

\medskip

\item the set of weak-equivalences is given by quasi-isomorphisms;

\medskip

\item the set of cofibrations is given by maps with left lifting property with respect to acyclic fibrations.

\medskip
\end{enumerate}

Furthermore, this model category structure is a right Bousfield localization of the model category structure transferred from $\mathcal{CC}_\infty$-coalgebras via the complete bar-cobar adjunction. Both of these statements ultimately follow from \cite[Chapter 10]{grignoulejay18}, see \cite[Proof of Proposition 4.6]{lucio2022contra} for the specific arguments. This gives an adjunction

\[
\begin{tikzcd}[column sep=7pc,row sep=3pc]
\mathsf{abs}~\mathcal{L}_\infty\text{-}\mathsf{alg}^{\mathsf{comp}} \left[\mathsf{W.eq}^{-1}\right] \arrow[r, shift left=1.1ex, "\mathrm{Id}"{name=F}]      
&\mathsf{abs}~\mathcal{L}_\infty\text{-}\mathsf{alg}^{\mathsf{comp}} \left[\mathsf{Q.iso}^{-1}\right]~, \arrow[l, shift left=.75ex, "\mathbb{L}\mathrm{Id}"{name=U}]
\arrow[phantom, from=F, to=U, , "\dashv" rotate=90]
\end{tikzcd}
\]

between the $\infty$-category of complete absolute $\mathcal{L}_\infty$-algebras up to weak-equivalences and the $\infty$-category of complete absolute $\mathcal{L}_\infty$-algebras up to quasi-isomorphisms, and complete absolute $\mathcal{L}_\infty$-algebras up to quasi-isomorphisms form a coreflective sub-$\infty$-category of complete absolute $\mathcal{L}_\infty$-algebras. A fibrant-cofibrant absolute $\mathcal{L}_\infty$-algebra is in this sub-$\infty$-category if and only if it is \textit{colocal} with respect to quasi-isomorphisms.

\begin{Proposition}
Let $X$ be a pointed connected simplicial set. The complete absolute $\mathcal{L}_\infty$-algebra $\mathcal{L}_*(X)$ is a colocal object with respect to quasi-isomorphisms. 
\end{Proposition}

\begin{proof}
Let $X$ be a pointed connected simplicial set, and let $f: \mathfrak{g} \qi \mathfrak{h}$ be a quasi-isomorphism between two absolute $\mathcal{L}_\infty$-algebras. First, notice that 
\[
\mathrm{Map}_{\mathsf{abs}~\mathcal{L}_\infty\text{-}\mathsf{alg}}(\mathcal{L}_*(X),-) \simeq \mathrm{Map}_{\mathsf{sSet}_*}(X,\mathcal{R}_*(-))~,
\]
since these two functor are adjoints at the $\infty$-categorical level. Now, it suffices to notice that since $f$ is a quasi-isomorphism, $\mathcal{R}_*(f): \mathcal{R}_*(\mathfrak{g}) \longrightarrow \mathcal{R}_*(\mathfrak{h})$ induces a weak-equivalence between the connected components of the Maurer--Cartan element $0$, by Theorem \ref{thm: Berglund isomorphism}. Finally, since $X$ is pointed and connected, the mapping space out of it will only detect the connected component where the base point lives, thus post-composing with $\mathcal{R}_*(f)$ induces a weak-equivalence of mapping spaces and therefore $\mathcal{L}_*(X)$ is colocal with respect to quasi-isomorphisms. 
\end{proof}

\begin{Remark}
In other words, the homotopical image of $\mathcal{L}_*$ lands in the coreflective sub-$\infty$-category of complete absolute $\mathcal{L}_\infty$-algebras up to quasi-isomorphisms. 
\end{Remark}

Finally, let us now compare these rational models given by complete absolute $\mathcal{L}_\infty$-algebras up to quasi-isomorphism with the other Lie-type models present in the literature. The adjunction of Proposition \ref{prop: adjonction absolute non-courbe et non-absolute} is a Quillen adjunction when one considers on both sides the models structures where weak-equivalences are given by quasi-isomorphism and where fibrations are given by degree-wise surjections, since the functor $\mathrm{Res}$ obviously preserves both. It induces of the following adjunction 
\[
\begin{tikzcd}[column sep=7pc,row sep=3pc]
\mathcal{L}_\infty\text{-}\mathsf{alg}\left[\mathsf{Q.iso}^{-1}\right] \arrow[r, shift left=1.1ex, "\mathrm{Abs}"{name=F}]      
&\mathsf{abs}~\mathcal{L}_\infty\text{-}\mathsf{alg}^{\mathsf{comp}}\left[\mathsf{Q.iso}^{-1}\right]~, \arrow[l, shift left=.75ex, "\mathrm{Res}"{name=U}]
\arrow[phantom, from=F, to=U, , "\dashv" rotate=-90]
\end{tikzcd}
\]
between the $\infty$-categories of absolute $\mathcal{L}_\infty$-algebras up to quasi-isomorphism and of $\mathcal{L}_\infty$-algebras up to quasi-isomorphism. 

\begin{Proposition}
Let $X$ be a pointed connected finite type nilpotent simplicial set. The derived counit of adjunction 

\[
\epsilon_X: (\mathbb{L}\mathrm{Abs})~\mathrm{Res}~\mathcal{L}_*(X) \longrightarrow \mathcal{L}_*(X)
\]
\vspace{0.1pc}

is a quasi-isomorphism. Thus the functor $\mathrm{Res}$ is homotopically fully faithful on the homotopical essential image of pointed connected finite type nilpotent rational spaces.
\end{Proposition}

\begin{proof}
Using the triangle identities of the adjunction and the fact that $\mathrm{Res}$ is conservative, it is enough to show that $\mathrm{Res}~\mathcal{L}_*(X)$ is \textit{homotopy complete}, that is, that the derived unit of this adjunction is a quasi-isomorphism for this $\mathcal{L}_\infty$-algebra. Up to quasi-isomorphism, one can replace this quasi-free completed $\mathcal{L}_\infty$-algebra with a quasi-free completed Lie algebra, and apply the criterion of \cite[Theorem 1]{nilpotentLiemodels} to show that $\mathrm{Res}~\mathcal{L}_*(X)$ is quasi-isomorphic to the analogue non-completed quasi-free  Lie algebra, which is itself homotopy complete. An analogue statement can also be found in \cite[Theorem 5.22 and Remark 5.23]{campos2020lie}.
\end{proof}

As a conclusion, let us mentioned that for pointed connected finite type nilpotent rational spaces, the $\mathcal{L}_\infty$-algebra $\mathrm{Res}~\mathcal{L}_*(X)$ is quasi-isomorphic to more classical models like Quillen models of \cite{Quillen69}, the Neisendorfer model of \cite{Neisendorfer} and to more recent models like those of Buij--Félix--Murillo--Tanré in \cite{Buijandco} and to those of Robert-Nicoud--Vallette in \cite{robertnicoud2020higher}. 

\medskip

\textbf{The general case.} We transfer the Kan--Quillen model structure on simplicial sets to the category of complete curved absolute $\mathcal{L}_\infty$-algebras, thus obtaining another model category structure on complete curved absolute $\mathcal{L}_\infty$-algebras. The identity functor can be promoted to a Quillen adjunction, and complete curved absolute $\mathcal{L}_\infty$-algebras with this model structure turn out to be a coreflective sub-$\infty$-category of the $\infty$-category of $u\mathcal{CC}_\infty$-coalgebras localized at quasi-isomorphisms.

\begin{Definition}[$\alpha$-quasi-isomorphisms]
Let $f: \mathfrak{g} \longrightarrow \mathfrak{h}$ be a morphism of complete curved $\mathcal{L}_\infty$-algebras. It is an $\alpha$-\textit{quasi-isomorphisms} if 

\begin{enumerate}
\medskip

\item it induces a bijection 
\[
\mathcal{MC}(f)/\sim_{\mathrm{gauge}}: \mathcal{MC}(\mathfrak{g})/\sim_{\mathrm{gauge}}~ \longrightarrow \mathcal{MC}(\mathfrak{h})/\sim_{\mathrm{gauge}}~,
\]

at the set of Maurer--Cartan elements up to gauge equivalences;

\medskip

\item and if, for every Maurer--Cartan element $\alpha$ in $\mathfrak{g}$, the induced map between the $\alpha$-homology groups
\[
f^\alpha: \mathrm{H}^{\alpha}_*(\mathfrak{g}) \longrightarrow \mathrm{H}^{f(\alpha)}_*(\mathfrak{h})~,
\]
is an isomorphism. 
\end{enumerate}
\end{Definition}

\begin{Proposition}
Let $f: \mathfrak{g} \qi \mathfrak{h}$ be a weak equivalence of complete curved $\mathcal{L}_\infty$-algebras. Then it is a $\alpha$-quasi-isomorphism.
\end{Proposition}

\begin{proof}
Since the functor $\mathcal{R}$ sends weak equivalences to weak homotopy equivalences, $f$ induces a bijection between the Maurer--Cartan elements up to gauge equivalences and is an $\alpha$-quasi-isomorphism in degrees $\geq 1$ by Theorem \ref{thm: Berglund isomorphism}. 

\medskip

Let us deal with degrees $\leq 0$: it will follow from considering spheres in negative degrees in Lemma \ref{lemma: structure de cogèbre sur la sphere}. Indeed, we can consider the same underlying chain complex as in the lemma, where the class $[n]$ can now be in any degree. The same formulas endow it with a $u\mathcal{CC}_\infty$-coalgebra structure as it is straightforward to check that, by degree reasons, these are the only possible operations when $n \leq 0$. Finally, the bijection of Theorem \ref{thm: Berglund isomorphism} generalizes to 
\[
\mathrm{Hom}_{u\mathcal{CC}_\infty\text{-}\mathsf{coalg}}\left(C_*^c(\mathbb{S}^n), \widehat{\mathrm{B}}_\iota(\mathfrak{g}) \right)/\sim_{\mathrm{hmt}} ~ \cong \mathrm{H}^{\alpha}_n(\mathfrak{g})~,
\]
for all $n$ in $\mathbb{Z}$. Since the left hand side is invariant under weak equivalences, so is the right hand side. 
\end{proof}

Using the right transfer theorem, we endow complete curved absolute $\mathcal{L}_\infty$-algebras with a model structure transferred from simplicial sets, where weak equivalences are given by $\alpha$-quasi-isomorphisms in degrees $\geq 1$.

\begin{theorem}\label{thm: intrinsic model structure}
There exists a model category structure on the category of complete curved absolute $\mathcal{L}_\infty$-algebras, determined by the following classes of maps
\begin{enumerate}
\medskip

\item the class of weak-equivalences is given by $\alpha$-quasi-isomorphisms in degrees $\geq 1$;

\medskip

\item the class of fibrations is given by maps $f$ which are sent to fibrations under the functor $\mathcal{R}$, and they include all degree-wise surjections;

\medskip

\item the class of cofibrations is given by left lifting property with respect to acyclic fibrations.
\end{enumerate}
\end{theorem}

\begin{proof}
The idea is to apply the right transfer theorem using the adjunction $\mathcal{L} \dashv \mathcal{R}$, where we consider the Kan--Quillen model structure on simplicial sets. The class of maps which are sent to weak-equivalences can be easily identified with $\alpha$-quasi-isomorphisms, using Theorems \ref{thm: characterisation pi zero et jauges} and \ref{thm: Berglund isomorphism}.

\medskip

We will rely on the "old" model category structure on complete curved absolute $\mathcal{L}_\infty$-algebras, transferred from $u\mathcal{CC}_\infty$-coalgebras, to show that it satisfies the right transfer hypothesis. First, notice that any "old" fibration and any "old" weak-equivalence is in the new classes of fibrations and weak-equivalences determined by the adjunction $\mathcal{L} \dashv \mathcal{R}$. Hence every complete curved absolute $\mathcal{L}_\infty$-algebra is fibrant with respect to the new fibrations, and the diagonal morphism $\mathfrak{g} \longrightarrow \mathfrak{g} \oplus \mathfrak{g}$ can be factored into a weak-equivalence followed by a fibration, simply using the old factorization. This concludes the proof.
\end{proof}

\begin{Remark}
This model category structure is very closely related to those constructed in \cite[Section 8.1]{Buijandco} and in \cite[Section 6.4]{robertnicoud2020higher}. Nevertheless, here fibrations do not admit the same easy description as morphisms which are surjective in degrees $\geq 0$. Indeed, consider the following counter-example: let $\mathfrak{g}_\emptyset$ be the model for the empty set constructed in Example \ref{Example: two point set and empty set}. Now let $\mathfrak{g}$ be any curved absolute $\mathcal{L}_\infty$-algebra. The canonical map $\mathfrak{g}_\emptyset \rightarrowtail \mathfrak{g}$ is sent to the empty bundle $\emptyset \longrightarrow \mathcal{R}(\mathfrak{g})$, which is a Kan fibration. However, the canonical map $\mathfrak{g}_\emptyset \rightarrowtail \mathfrak{g}$ is clearly not surjective in degrees $\geq 0$. 
\end{Remark}

\begin{Remark}
There should be a model structure on complete curved absolute $\mathcal{L}_\infty$-algebras where the class of weak equivalences is given by \textit{all} $\alpha$-quasi-isomorphisms and not just by $\alpha$-quasi-isomorphisms in degrees $\geq 1$, akin to the model structure on the complete absolute $\mathcal{L}_\infty$-algebras where the class of weak equivalences is given by \textit{all} quasi-isomorphisms. 
\end{Remark}

Since the classes of fibration and weak-equivalences in the transferred model structure from $u\mathcal{CC}_\infty$-coalgebras are included in the classes of fibrations and weak-equivalences in the model structure of Theorem \ref{thm: intrinsic model structure}, the identity functor from the first to the second becomes a right Quillen functor. This results in the following adjunction 

\[
\begin{tikzcd}[column sep=7pc,row sep=3pc]
\mathsf{curv}~\mathsf{abs}~\mathcal{L}_\infty\text{-}\mathsf{alg}^{\mathsf{comp}} \left[\mathsf{W.eq}^{-1}\right] \arrow[r, shift left=1.1ex, "\mathrm{Id}"{name=F}]      
&\mathsf{curv}~\mathsf{abs}~\mathcal{L}_\infty\text{-}\mathsf{alg}^{\mathsf{comp}} \left[\alpha\text{-}\mathsf{q.iso}^{-1}_{\geq 1}\right]~, \arrow[l, shift left=.75ex, "\mathbb{L}\mathrm{Id}"{name=U}]
\arrow[phantom, from=F, to=U, , "\dashv" rotate=90]
\end{tikzcd}
\]

between the $\infty$-category of complete curved absolute $\mathcal{L}_\infty$-algebras up to weak-equivalences and the $\infty$-category of complete curved absolute $\mathcal{L}_\infty$-algebras up to $\alpha$-quasi-isomorphisms in degrees $\geq 1$. 

\begin{Proposition}
The functor $\mathbb{L}\mathrm{Id}$ is homotopically fully faithful.
\end{Proposition}

\begin{proof}
It is straightforward to see that the derived unit of the adjunction is an isomorphism.
\end{proof}

\begin{Remark}
Let $X$ be a simplicial set. Its image $\mathcal{L}(X)$ will land on the coreflective sub-$\infty$-category of complete curved absolute $\mathcal{L}_\infty$-algebras up to $\alpha$-quasi-isomorphisms.
\end{Remark}

\section{Deformation theory}
In this section, we consider different applications of the theory developed so far. First, we show that convolution (curved in our case) $\mathcal{L}_\infty$-algebras constructed by D. Robert-Nicoud and F. Wierstra in \cite{RobertNicoudWierstra17} carry naturally an "absolute" structure without the need of an underlying filtration. These convolution algebras encode $\infty$-morphisms between algebras over a (non-augmented) dg operad as their Maurer--Cartan elements and all the homotopies between those. Secondly, for any derived affine stack, we construct a curved absolute $\mathcal{L}_\infty$-algebra which is "geometrical model". More concretely, we shows that this geometrical model recovers the formal neighborhood of any of the $\mathbb{L}$-points of our derived affine stack, where $\mathbb{L}$ is a finite algebraic extension of the base field $\kk$. These results should lead to a "non-pointed" version of formal moduli problems.

\subsection{Convolution curved absolute $\mathcal{L}_\infty$-algebras and $\infty$-morphisms of algebras}\label{Subsection: convolution} Let $\mathcal{P}$ denote a dg operad and let $\mathcal{C}$ denote a conilpotent curved cooperad. 

\begin{Proposition}\label{prop: curved twisting morphisms}
There is a bijection 

\[
\mathrm{Tw}(\mathcal{C},\mathcal{P}) \cong \mathrm{Hom}_{\mathsf{curv}~\mathsf{ab}~\mathsf{pOp}}\left(\widehat{\Omega}^{\mathrm{s.c}} u\mathcal{C}om^*, \mathcal{H}om(\mathcal{C},\mathcal{P})\right)~,
\]

between morphisms of curved absolute partial operads $\phi_\alpha : \widehat{\Omega}^{\mathrm{s.c}} u\mathcal{C}om^* \longrightarrow \mathcal{H}om(\mathcal{C},\mathcal{P})$ and curved twisting morphisms $\alpha: \mathcal{C} \longrightarrow \mathcal{P}$ such that $\alpha(1)$ is the zero morphism.
\end{Proposition}

\begin{proof}
First, notice that the convolution curved partial operad of \cite[Subsection 6.2]{lucio2022curved} is in fact a curved absolute partial operad in the sense of \cite[Subsection 10.2]{lucio2022curved}. Indeed, one can check that infinite sums of compositions in it have a well-defined image in it since any such sum is \textit{locally finite} because of the conilpotency of $\mathcal{C}$. There is an isomorphism 

\[
\mathcal{H}om\left(u\mathcal{C}om^*,\mathcal{H}om(\mathcal{C},\mathcal{P})\right) \cong \mathcal{H}om(\mathcal{C},\mathcal{P})
\]
\vspace{0.2pc}

of curved absolute partial operads since $u\mathcal{C}om^*(n) = \kk$ for all $n \geq 0$. Therefore,

\[
\mathrm{Hom}_{\mathsf{curv}~\mathsf{ab}~\mathsf{pOp}}\left(\widehat{\Omega} u\mathcal{C}om^*, \mathcal{H}om(\mathcal{C},\mathcal{P})\right) \cong \mathrm{Tw}\left(u\mathcal{C}om^*,\mathcal{H}om(\mathcal{C},\mathcal{P})\right) \cong \mathrm{Tw}(\mathcal{C},\mathcal{P})~,
\]
\vspace{0.2pc}

using the curved operadic bar-cobar adjunction of \cite[Section 6]{lucio2022curved}. Now one can check that $\widehat{\Omega}^{\mathrm{s.c}} u\mathcal{C}om^*$ represents exactly those curved twisting morphisms that have a trivial arity one component.
\end{proof}

\begin{Remark}
For an analogue statement concerning convolution $\mathcal{L}_\infty$-algebra structures, see the constructions in \cite{RobertNicoudWierstra17}.
\end{Remark}

Let $(A,\allowbreak \gamma_A,d_A)$ be a dg $\mathcal{P}$-algebra and let $(C,\Delta_C,d_C)$ be a curved $\mathcal{C}$-coalgebra. The graded module of graded morphisms $\mathrm{hom}(C,A)$ is naturally a pdg module endowed with the pre-differential $\partial(f) \coloneqq d_A \circ f - f \circ d_C$.

\begin{Proposition}
Let $(A,\gamma_A,d_A)$ be a dg $\mathcal{P}$-algebra and let $(C,\Delta_C,d_C)$ be a curved $\mathcal{C}$-coalgebra. The pdg module 
\[
\left(\mathrm{hom}(C,A),\partial \right)
\]
can be endowed with a curved $\mathcal{L}_\infty$-algebra structure is given by 
\[
\begin{tikzcd}[column sep=2.5pc,row sep=0.5pc]
\gamma_{\mathrm{hom}(C,A)}: \displaystyle \bigoplus_{n \geq 0} \widehat{\Omega} u\mathcal{C}om^*(n) \otimes_{\mathbb{S}_n} \mathrm{hom}(C,A)^{\otimes n} \arrow[r]
&\mathrm{hom}(C,A) \\
c_n (f_1 \otimes \cdots \otimes f_n) \arrow[r,mapsto]
&\displaystyle \gamma_A \circ \left[\phi_\alpha(c_n) \otimes (f_1 \otimes \cdots \otimes f_n) \right] \circ \Delta_C~,
\end{tikzcd}
\]
where $c_n$ denotes the $n$-corolla.
\end{Proposition}

\begin{proof}
This follows directly from Proposition \ref{prop: curved twisting morphisms} by considering the composition 
\[
\begin{tikzcd}[column sep=4pc,row sep=0.5pc]
\widehat{\Omega} u\mathcal{C}om^* \arrow[r," \phi_\alpha"] 
&\mathcal{H}om(\mathcal{C},\mathcal{P}) \arrow[r,"\lambda_{\mathrm{hom}(C,A)} "]
&\mathrm{End}_{\mathrm{hom}(C,A)}~,
\end{tikzcd}
\]
of morphisms of curved operads, where $\lambda_{\mathrm{hom}(C,A)}$ is the natural curved algebra structure of $\mathrm{hom}(C,A)$ over $\mathcal{H}om(\mathcal{C},\mathcal{P})$.
\end{proof}

\begin{Remark}[Local nilpotency]
One can notice that the Maurer--Cartan equation is always well-defined without the need of imposing any filtration on the space $\mathrm{hom}(C,A)$. Indeed, let $f: C \longrightarrow A$ be a degree zero map, then the sum
\[
\gamma_A \circ \left[\sum_{\substack{n \geq 0 \\ n \neq 1}} \phi_\alpha(c_n) \otimes f^{\otimes n} \right] \circ \Delta_C(c) + \partial(f)(c)
\]
is well defined for any element $c$ in $C$, since there are only a finite number of non-zero terms in $\Delta_C(c)$. Filtrations introduced to make this type of sum converge like in \cite{DolgushevHoffnungRogers14} are redundant. 
\end{Remark}

\begin{Proposition}
The extension of $\gamma_{\mathrm{hom}(C,A)}$ given by 
\[
\begin{tikzcd}[column sep=1pc,row sep=-0.5pc]
\displaystyle \prod_{n \geq 0} \widehat{\Omega} u\mathcal{C}om^*(n)~\widehat{\otimes}_{\mathbb{S}_n}~\mathrm{hom}(C,A)^{\otimes n} \arrow[r]
&\mathrm{hom}(C,A) \\
\displaystyle \sum_{\substack{n\geq 0 \\ \omega \geq 0}} \sum_{\tau \in \mathrm{CRT}_\omega^n } \lambda_\tau \tau(f_1 \otimes \cdots \otimes f_n) \arrow[r,mapsto]
&\displaystyle \gamma_A \circ \left[\sum_{\substack{n \geq 0 \\ \omega \geq 0}} \sum_{\tau \in \mathrm{CRT}_\omega^n } \phi_\alpha(\tau) \otimes (f_1 \otimes \cdots \otimes f_n) \right] \circ \Delta_C~,
\end{tikzcd}
\]
defines a structure of curved absolute $\mathcal{L}_\infty$-algebra, denoted by $\mathrm{hom}(C,A)$.
\end{Proposition}

\begin{proof}
One can check by hand that this formula satisfies conditions \ref{pdg condition}, \ref{associativity condition}, and \ref{curved condition}.
\end{proof}

\begin{Remark}
This convolution curved absolute $\mathcal{L}_\infty$-algebra structure coincides with the one constructed in \cite{grignou2022mapping} between a dg $\Omega \mathcal{C}$-algebra and a curved $\mathcal{C}$-coalgebra, using the Hopf comodule structure of the cofibrant dg operad $\Omega \mathcal{C}$.
\end{Remark}

\begin{Proposition}
Let $A$ and $B$ be two dg $\mathcal{P}$-algebras. The simplicial set $\mathcal{R}(\mathrm{hom}(\mathrm{B}_\alpha A, B))$ has as $0$-simplices the set of $\infty_\alpha$-morphisms between $A$ and $B$.
\end{Proposition}

\begin{proof}
Maurer--Cartan elements correspond by definition to curved twisting morphisms $\mathrm{B}_\alpha A$ between and $B$, which are in bijection with morphisms of curved $\mathcal{C}$-coalgebras between $\mathrm{B}_\alpha A$ and $\mathrm{B}_\alpha B$.
\end{proof}

Therefore $1$-simplices in $\mathcal{R}(\mathrm{hom}(\mathrm{B}_\alpha A, B))$ induce a notion of \textit{homotopies} between these $\infty_\alpha$-morphism, and similarly for $n$-simplices which correspond to higher homotopies. Using these methods to simplicially enrich the category of dg $\mathcal{P}$-algebras and thus obtain a model for the $\infty$-category of dg $\mathcal{P}$-algebras will be the subject of a future work. 

\begin{Remark}[The dual case of $\infty$-morphism of coalgebras]
In \cite{grignou2022mapping}, Brice Le Grignou constructs a convolution curved absolute $\mathcal{L}_\infty$-algebra from a dg $\Omega\mathcal{C}$-coalgebra and a curved $\mathcal{C}$-algebra. Integrating these convolution algebras provides us with the right set of $\infty_\iota$-morphisms between two coalgebras. This method also provides us with a way in which one can simplicially enrich the category of dg $\Omega\mathcal{C}$-coalgebras with $\infty_\iota$-morphisms.
\end{Remark}

\begin{Remark}[Extension to dg properads]
Let $\Omega \mathcal{C}$ be a properad. Consider the notion of $\infty$-morphism for gebras over $\Omega \mathcal{C}$ as defined in  \cite{hoffbeck2019properadic}. The convolution construction that controls them is again an absolute $\mathcal{L}_\infty$-algebra, as it is also \textit{locally nilpotent}. Thus the integration theory of absolute $\mathcal{L}_\infty$-algebras could be applied to  settle a suitable simplicial enrichment for the category of $\Omega \mathcal{C}$-gebras with $\infty$-morphisms.
\end{Remark}

\newcommand\htimes{\stackrel{\mathclap{\normalfont\mbox{\tiny{h}}}}{\times}}

\subsection{The formal geometry of Maurer--Cartan spaces}
Let $A$ be a dg $u\mathcal{C}om$-algebra, viewed equivalently as a \textit{derived affine stack}. Its \textit{functor of points} is given by 

\[
\begin{tikzcd}[column sep=1pc,row sep=0.5pc]
\mathrm{Spec}(A)(-): \mathsf{dg}~u\mathcal{C}om\text{-}\mathsf{alg}_{\geq 0} \arrow[r]
&\mathsf{sSet} \\
B \arrow[r,mapsto]
&\mathrm{Spec}(A)(B)_\bullet \coloneqq \mathbb{R}\mathrm{Hom}_{\mathsf{dg}~u\mathcal{C}om\text{-}\mathsf{alg}}(A, B \otimes \Omega_\bullet)~,
\end{tikzcd}
\]
where $\Omega_\bullet$ is again the simplicial Sullivan algebra and where $B$ is a dg $u\mathcal{C}om$-algebra concentrated in non-negative homological degrees, that is, a \textit{derived affine scheme}. The simplicial set $\mathrm{Spec}(A)(B)_\bullet$ is called the $B$-points of $A$.

\medskip

The first goal of this subsection is to construct a functor from dg $u\mathcal{C}om$-algebras to curved absolute $\mathcal{L}_\infty$-algebras that recovers the "formal geometry" of $A$. Let $x: A \longrightarrow \kk$ be a $\kk$-point of $A$, it provides an $\kk$-\textit{augmentation} of $A$. The standard way of looking at \textit{infinitesimal thickenings} of $x$ inside $A$ has been to test $A$ against $\kk$-augmented dg Artinian algebras. 

\begin{Definition}[$\kk$-augmented dg Artinian algebra]
Let $R$ be a dg $u\mathcal{C}om$-algebra with homology concentrated in non-negative homological degrees. It is a $\kk$\textit{-augmented dg Artinian algebra} if it satisfies the following conditions.

\medskip

\begin{enumerate}
\item Its homology is degree-wise finite dimensional and concentrated in a finite number of degrees.

\medskip

\item There is an unique augmentation morphism $p: R \longrightarrow \kk$ (hence $R$ is local).

\medskip

\item Let $\overline{R}$ be the non-unital dg commutative algebra given by the kernel of $p$. Then $\overline{R}$ is nilpotent.
\end{enumerate}
\end{Definition}

This notion allows us to define the notion of a $\kk$-pointed formal moduli problem.

\begin{Definition}[$\kk$-pointed formal moduli problem]
Let 
\[
F: \mathsf{dg}~\mathsf{Art}\text{-}\mathsf{alg}^{\kk\text{-}\mathsf{aug}}_{\geq 0} \longrightarrow \mathsf{sSet}
\]

be a functor. It is a $\kk$-\textit{pointed formal moduli problem} if it satisfies the following conditions

\medskip

\begin{enumerate}
\item We have that $F(\kk) \simeq \{*\}$.

\medskip

\item The functor $F$ sends quasi-isomorphisms to weak-equivalences.

\medskip

\item The functor $F$ preserves homotopy pullbacks of the of dg Artinian algebras $X,Y$ and $Z$ 
\[
\begin{tikzcd}
X \htimes_Z Y \arrow[r] \arrow[d]
&Y \arrow[d,"\pi_1"] \\
X \arrow[r,"\pi_2 ",swap]
&Z \\
\end{tikzcd}
\]

such that $\pi_1$ and $\pi_2$ are surjections on the zeroth homology groups $\mathrm{H}_0$. 
\end{enumerate}
\end{Definition}

\begin{Example}
Any $\kk$-augmented derived affine stack $A$ defines such a pointed deformation problem by considering 

\[
\mathrm{Spec}^*(A)(-)_\bullet = \mathbb{R}\mathrm{Hom}_{\mathsf{dg}~u\mathcal{C}om\text{-}\mathsf{alg}^{\kk\text{-}\mathsf{aug}}}(A, - \otimes \Omega_\bullet)~,
\]

where morphisms of $\kk$-augmented dg commutative algebras are required to preserve the augmentation. This functor preserves homotopy limits and sends quasi-isomorphisms to weak homotopy equivalences, hence it defines a pointed formal moduli problem. 
\end{Example}

The $\infty$-category of $\kk$-pointed formal moduli problems is equivalent to the $\infty$-category of dg Lie algebras over $\kk$, or for that matter, to the $\infty$-category of $\mathcal{L}_\infty$-algebras, both of them constructed by consider quasi-isomorphisms as weak-equivalences. This was shown independently by J. Lurie in \cite{Lurie11} by working directly with $\infty$-categories and by J.P Pridham \cite{Pridham10} using model categories. 

\medskip

\textbf{Non-pointed version.} Heuristically speaking, curved absolute $\mathcal{L}_\infty$-algebras should correspond to a \textit{non-pointed} version of formal moduli problems. A first idea into what this notion might be is the following: in a pointed formal moduli problem, the point around which one considers the infinitesimal thickening is already specified in advance. In the non-pointed version, one should remove the $\kk$-augmented hypothesis. 

\medskip

As a first approach to non-pointed deformation problems, we build a curved absolute $\mathcal{L}_\infty$-algebra $\mathfrak{g}_A$ from a derived affine stack $A$. We consider the following commutative square given by \cite[Theorem 2.22]{lucio2022contra}:

\[
\begin{tikzcd}[column sep=5pc,row sep=5pc]
\left(\mathsf{dg}~u\mathcal{C}om\text{-}\mathsf{alg}\right)^{\mathsf{op}} \arrow[r,"\mathrm{B}_\pi^{\mathsf{op}}"{name=B},shift left=1.1ex] \arrow[d,"(-)^\circ "{name=SD},shift left=1.1ex ]
&\left(\mathsf{curv}~\mathcal{L}_\infty\text{-}\mathsf{coalg}^{\mathsf{conil}}\right)^{\mathsf{op}} \arrow[d,"(-)^*"{name=LDC},shift left=1.1ex ] \arrow[l,"\Omega_\pi^{\mathsf{op}}"{name=C},,shift left=1.1ex]  \\
\mathsf{dg}~u\mathcal{C}om\text{-}\mathsf{coalg} \arrow[r,"\widehat{\Omega}_\pi "{name=CC},shift left=1.1ex]  \arrow[u,"(-)^*"{name=LD},shift left=1.1ex ]
&\mathsf{curv}~\mathsf{abs}~\mathcal{L}_\infty\text{-}\mathsf{alg}^{\mathsf{comp}}~, \arrow[l,"\widehat{\mathrm{B}}_\pi"{name=CB},shift left=1.1ex] \arrow[u,"(-)^\vee"{name=TD},shift left=1.1ex] \arrow[phantom, from=SD, to=LD, , "\dashv" rotate=0] \arrow[phantom, from=C, to=B, , "\dashv" rotate=-90]\arrow[phantom, from=TD, to=LDC, , "\dashv" rotate=0] \arrow[phantom, from=CC, to=CB, , "\dashv" rotate=-90]
\end{tikzcd}
\] 

which relates the category of dg $u\mathcal{C}om$-algebras to the category of curved absolute $\mathcal{L}_\infty$-algebras.

\begin{Remark}\label{Rmk: pas de top model}
The category of dg $u\mathcal{C}om$-coalgebras does not admit a model category structure transferred from dg modules, where weak-equivalences are quasi-isomorphisms and where cofibrations are degree-wise monomorphisms. Therefore the left adjunction on this square is not a Quillen adjunction in this case, see \cite[Section 9]{grignoulejay18} for more details.
\end{Remark}

\begin{Definition}[Geometrical model]
Let $A$ be a dg $u\mathcal{C}om$-algebra. Its \textit{geometrical model} $\mathfrak{g}_A$ is the complete curved absolute $\mathcal{L}_\infty$-algebra given by $(\mathrm{B}_\pi A)^*$. It defines a functor 
\[
\mathfrak{g}_{(-)} \coloneqq (\mathrm{B}_\pi(-))^*: \left(\mathsf{dg}~u\mathcal{C}om\text{-}\mathsf{alg}\right)^{\mathsf{op}} \longrightarrow \mathsf{curv}~\mathsf{abs}~\mathcal{L}_\infty\text{-}\mathsf{alg}^{\mathsf{comp}}~.
\]
\end{Definition}

\begin{Remark}
This approach is close to the one of \cite{calaqueformal}, where the authors consider the linear dual of the bar construction relative to $\kappa$ of a dg $\mathcal{P}$-algebra. Nevertheless, they view it as a dg $\mathcal{P}^{\hspace{1pt}!}$-algebra instead of viewing this linear dual as dg algebra over the cooperad $\mathcal{P}^{\hspace{1pt}\ac}$. They use the linear dual of the bar construction in order to show an equivalence between the $\infty$-category of pointed formal moduli problems defined using dg Artinian $\mathcal{P}$-algebras and the $\infty$-category of dg $\mathcal{P}^{\hspace{1pt}\ac}$-algebras, under certain hypothesis on the operad $\mathcal{P}$.
\end{Remark}

\begin{Proposition}
The geometrical model functor $\mathfrak{g}_{(-)}$ sends quasi-isomorphisms to weak-equivalences. 
\end{Proposition}

\begin{proof}
Let $f: A \qi B$ be a quasi-isomorphism of dg $u\mathcal{C}om$-algebras, then $\mathrm{B}_\pi(f): \mathrm{B}_\pi A \qi \mathrm{B}_\pi B$ is a weak-equivalence of conilpotent curved $\mathcal{L}_\infty$-coalgebras, since every dg $u\mathcal{C}om$-algebra is fibrant. The linear dual $(-)^*$ preserves all weak-equivalences since it is a right Quillen functor and since every conilpotent curved $\mathcal{L}_\infty$-coalgebra is cofibrant (fibrant in the opposite category). 
\end{proof}

\begin{Definition}[dg Artinian algebra]
Let $A$ be a dg $u\mathcal{C}om$-algebra in non-negative homological degrees. It is a \textit{dg Artinian algebra} if its homology is degree-wise finite dimensional and bounded above. 
\end{Definition}

\begin{Example}
Consider a dg Artinian algebra $A$ concentrated in degree zero. Since $A$ is a finite dimensional $\kk$-algebra, it can be written as a product 
\[
A \cong \prod_{i = 0}^n A_i~, 
\]
where $A_i$ are local Artinian algebras. Thus one can consider the following examples

\medskip

\begin{enumerate}
\item the algebra $\kk^n$, which gives a collection of $n$ points, 

\medskip

\item any finite field extension $\mathbb{L}$ of $\kk$,

\medskip

\item classical thickenings like $\kk[t]/(t^n)$. 
\end{enumerate}
\end{Example}

\begin{lemma}\label{lemma: equivalence dg Artinien avec les uCC}
There is a an equivalence of $\infty$-categories between the opposite $\infty$-category of $u\mathcal{CC}_\infty$-coalgebras with total finite dimensional homology concentrated in non-positive degrees and the $\infty$-category of dg Artinian algebras.
\end{lemma}

\begin{proof}
There are equivalences of $\infty$-categories 

\hspace{-2pc}
\begin{tikzcd}[column sep=5pc,row sep=2.5pc]
\left(u\mathcal{CC}_\infty\text{-}\mathsf{coalg}_{\leq 0}^{\mathsf{f.d}}\right)^{\mathsf{op}}\left[\mathsf{Q.iso}^{-1}\right] \arrow[r, shift left=1.5ex, "(-)^*"{name=B}]
&u\mathcal{CC}_\infty\text{-}\mathsf{alg}_{\geq 0}^{\mathsf{f.d}}\left[\mathsf{Q.iso}^{-1}\right] \arrow[r, shift left=1.5ex, "\mathbb{L}\mathrm{Ind}_\varepsilon"{name=A}] \arrow[l, shift left=.75ex, "\mathbb{R}(-)^\circ"{name=D}]
&\mathsf{dg}~\mathcal{C}om\text{-}\mathsf{alg}_{\geq 0}^{\mathsf{f.d}}\left[\mathsf{Q.iso}^{-1}\right]~, \arrow[l, shift left=.75ex, "\mathrm{Res}_\varepsilon"{name=C}] \arrow[phantom, from=A, to=C, , "\dashv" rotate=-90] \arrow[phantom, from=B, to=D, , "\dashv" rotate=90]
\end{tikzcd}

where the first equivalence is given by \cite[Theorem 2.25]{lucio2022contra} and the second equivalence is induced by the quasi-isomorphism of dg operads $\varepsilon: \Omega \mathrm{B}u\mathcal{C}om \qi u\mathcal{C}om$, which respects the degrees of the homology.
\end{proof}

\begin{theorem}\label{thm: modèle géométrique sur les points Artiniens}
Let $A$ be a dg $u\mathcal{C}om$-algebra. Let $B$ be a dg Artinian algebra. There is a weak homotopy equivalence of simplicial sets

\[
\mathrm{Spec}(A)(B) \simeq \mathcal{R}(\mathrm{hom}(\mathbb{R}(\mathrm{Res}_\varepsilon B)^\circ, \mathfrak{g}_A))~,
\]
\vspace{0.1pc}

which is natural in $B$ and in $A$, where $\mathbb{R}(-)^\circ$ is the derived Sweedler dual functor. A model, non-functorial in $B$, for this simplicial set is given by $\mathcal{R}(\mathrm{hom}((\mathrm{H}_*B)^*, \mathfrak{g}_A))$, where the homology $(\mathrm{H}_*B)^*$ is endowed with a transferred $u\mathcal{CC}_\infty$-algebra structure via the homotopy transfer theorem. 
\end{theorem}

\begin{proof}
We start by computing the derived mapping space with the classical bar-cobar resolution
\[
\mathrm{Spec}(A)(B) \simeq \mathrm{Hom}_{\mathsf{dg}~u\mathcal{C}om\text{-}\mathsf{alg}}(\Omega_\pi \mathrm{B}_\pi A, B \otimes \Omega_\bullet)~,
\]

since $\Omega_\pi \mathrm{B}_\pi A$ is a cofibrant resolution of $A$. Using the same methods as in the proof of Proposition \ref{prop: equivalence cochains CC with Apl}, we show that

\[
\mathrm{Hom}_{\mathsf{dg}~u\mathcal{C}om\text{-}\mathsf{alg}}(\Omega_\pi \mathrm{B}_\pi A, B \otimes \Omega_\bullet) \simeq \mathrm{Hom}_{u\mathcal{CC}_\infty\text{-}\mathsf{alg}}(\Omega_\iota \mathrm{B}_\iota \mathrm{Res}_\epsilon A, \mathrm{Res}_\epsilon(B) \otimes C^*_c(\Delta^\bullet))~.
\]
\vspace{0.1pc}

We consider now the following quasi-isomorphism of $u\mathcal{CC}_\infty$-algebras 

\[
\begin{tikzcd}[column sep=4pc,row sep=0pc]
\mathbb{R}\mathrm{Res}_\epsilon(B) \otimes C^*_c(\Delta^\bullet) \arrow[r,"\mathbb{R}\eta_{\mathrm{Res}_\epsilon(B)} \otimes \mathrm{id}"]
&((\mathbb{R}\mathrm{Res}_\varepsilon B)^\circ)^* \otimes C^*_c(\Delta^\bullet) \arrow[r,"\Gamma"]
&\left((\mathbb{R}\mathrm{Res}_\varepsilon B)^\circ \otimes C_*^c(\Delta^\bullet))\right)^*~,
\end{tikzcd}
\]

where $\mathbb{R}\eta_{\mathrm{Res}_\epsilon(B)}$ is the derived unit of the $(-)^* \dashv (-)^\circ$ adjunction, which is a quasi-isomorphism since $B$ has degree-wise finite dimensional homology, and where $\Gamma$ is the lax monoidal structure of the linear dual functor $(-)^*$ with respect to the tensor product, which is a quasi-isomorphism since both objects have degree-wise finite dimensional homology. This quasi-isomorphism gives a direct natural weak-equivalence of simplicial sets 

\[
\begin{tikzcd}[column sep=0pc,row sep=2pc]
\mathrm{Hom}_{u\mathcal{CC}_\infty\text{-}\mathsf{alg}}(\Omega_\iota \mathrm{B}_\iota \mathrm{Res}_\epsilon A, \mathbb{R}\mathrm{Res}_\epsilon(B) \otimes C^*_c(\Delta^\bullet)) \arrow[d,"\simeq"] \\
\mathrm{Hom}_{u\mathcal{CC}_\infty\text{-}\mathsf{alg}}(\Omega_\iota \mathrm{B}_\iota \mathrm{Res}_\epsilon A, \left(\mathbb{R}(\mathrm{Res}_\varepsilon B)^\circ \otimes C_*^c(\Delta^\bullet))\right)^*)
\end{tikzcd}
\]

There is an isomorphism of simplicial sets 

\[
\begin{tikzcd}[column sep=0pc,row sep=2pc]
\mathrm{Hom}_{u\mathcal{CC}_\infty\text{-}\mathsf{alg}}(\Omega_\iota \mathrm{B}_\iota \mathrm{Res}_\varepsilon A, \left(\mathbb{R}(\mathrm{Res}_\varepsilon B)^\circ \otimes C_*^c(\Delta^\bullet)\right)^*) \arrow[d,"\cong"] \\
\mathrm{Hom}_{u\mathcal{CC}_\infty\text{-}\mathsf{coalg}}(\mathbb{R}(\mathrm{Res}_\varepsilon B)^\circ \otimes C_*^c(\Delta^\bullet),\left(\Omega_\iota \mathrm{B}_\iota \mathrm{Res}_\varepsilon A\right)^\circ)
\end{tikzcd}
\]

induced by the adjunction $(-)^\circ \dashv (-)^*$. We compute that 

\[
\left(\Omega_\iota \mathrm{B}_\iota \mathrm{Res}_\varepsilon A\right)^\circ \cong \widehat{\mathrm{B}}_\iota (\mathrm{B}_\iota \mathrm{Res}_\varepsilon A)^* \cong \widehat{\mathrm{B}}_\iota (\mathrm{B}_\pi A)^*~.
\]
\vspace{0.1pc}

Finally we get

\begin{align*}
\mathrm{Hom}_{u\mathcal{CC}_\infty\text{-}\mathsf{coalg}}(\mathbb{R}(\mathrm{Res}_\varepsilon B)^\circ \otimes C_*^c(\Delta^\bullet),\widehat{\mathrm{B}}_\iota (\mathrm{B}_\pi A)^*) \\
\cong \mathrm{Hom}_{u\mathcal{CC}_\infty\text{-}\mathsf{coalg}}(C_*^c(\Delta^\bullet),\{\mathbb{R}(\mathrm{Res}_\varepsilon B)^\circ,\widehat{\mathrm{B}}_\iota (\mathrm{B}_\pi A)^*\}) \\
\cong \mathrm{Hom}_{u\mathcal{CC}_\infty\text{-}\mathsf{coalg}}(C_*^c(\Delta^\bullet),\widehat{\mathrm{B}}_\iota(\mathrm{hom}(\mathbb{R}(\mathrm{Res}_\varepsilon B)^\circ,(\mathrm{B}_\pi A)^*)) \\
\cong \mathrm{Hom}_{\mathsf{curv}~\mathsf{abs}~\mathcal{L}_\infty\text{-}\mathsf{alg}}(\widehat{\Omega}_\iota(C_*^c(\Delta^\bullet)),\mathrm{hom}(\mathbb{R}(\mathrm{Res}_\varepsilon B)^\circ,(\mathrm{B}_\pi A)^*) \\
\cong \mathcal{R}(\mathrm{hom}(\mathbb{R}(\mathrm{Res}_\varepsilon B)^\circ, \mathfrak{g}_A))~.
\end{align*}

Finally, using the homotopy transfer theorem one can compute that a model for $\mathbb{R}(\mathrm{Res}_\varepsilon B)^\circ$ is given by 
\[
\left(\Omega_\iota \mathrm{B}_\iota H_*B\right)^\circ \cong \widehat{\mathrm{B}}_\iota (\mathrm{B}_\iota H_*B)^* \cong \widehat{\mathrm{B}}_\iota (\mathrm{B}_\iota H_*B)^* \simeq \widehat{\mathrm{B}}_\iota \widehat{\Omega}_\iota(H_*B)^* \simeq H_*B~,
\]
using that $H_*B$ is degree-wise finite dimensional and bounded below. 
\end{proof}

\begin{Example}
Let $A$ be a derived affine scheme, that is, a dg $u\mathcal{C}om$-algebra concentrated in non-negative degrees. There is an isomorphism of constant simplicial sets
\[
\mathrm{Spec}(A)(\kk) \cong \mathcal{MC}(\mathfrak{g}_A) \cong \mathcal{R}(\mathfrak{g}_A)~,
\]
since $\mathfrak{g}_A$ is concentrated in non-positive degrees. For $A$ a derived affine stack, we have that
\[
\mathrm{Spec}(A)(\kk) \simeq \mathcal{R}(\mathfrak{g}_A)~.
\]
Geometrically speaking, these isomorphisms allow us to recover all the $\kk$-points of $A$ as well as the formal neighborhood \textit{of any of these points} from its geometrical model $\mathfrak{g}_A$. 

\medskip

Furthermore, given a finite field extension $\mathbb{L}$ of $\kk$, we have  
\[
\mathrm{Spec}(A)(\mathbb{L}) \cong \mathcal{MC}(\mathrm{hom}(\mathbb{L}^*,\mathfrak{g}_A) \cong \mathcal{MC}(\mathbb{L} \otimes \mathfrak{g}_A)~,
\]
so we can base-change $\mathfrak{g}_A$ to any finite extension $\mathbb{L}$, and recover the $\mathbb{L}$-points of $A$ as well as their formal neighborhoods. Finally, any of these combinations can be done for a finite number $n$ of points in $A$. 
\end{Example}

\begin{Example}[The complex numbers as an $\mathbb{R}$-algebra]
Let us compute a particular example. The field of complex numbers $\mathbb{C}$ is given, as an $\mathbb{R}$-algebra, by the following presentation
\[
\mathbb{C} \cong \frac{\mathbb{R}[x]}{(x^2 + 1)}~.
\]
This is a non-augmented $\mathbb{R}$-algebra, thus its Koszul dual will be a conilpotent curved $\mathcal{L}_\infty$-coalgebra, which can be dualized to a curved absolute $\mathcal{L}_\infty$-algebra. Using the constant-linear-quadratic Koszul duality for algebras over a binary operad developed in \cite{najibcurved}, its geometrical model admits a quite simple description:
\[
\mathfrak{g}_{\mathbb{C}} \simeq \mathbb{R}.y \oplus \mathbb{R}.\frac{1}{2}[y,y]~, 
\]
where $y$ is in degree $0$, the bracket of $y$ with itself $[y,y]$ is in degree $-1$, and other brackets are zero except for the curvature $l_0: \kk \longrightarrow \mathfrak{g}_{\mathbb{C}}$, which sends $1$ to $-[y,y]$. (Remember we work with the \textit{shifted} convention for the degree of the brackets). Notice that higher iterations of the bracket with itself vanish because of the Jacobi relation.

\medskip

The curved absolute $\mathcal{L}_\infty$-algebra $\mathfrak{g}_{\mathbb{C}}$ does not admit any Maurer--Cartan elements: if $\lambda.y$ was a Maurer--Cartan element, it would have to satisfy:
\[
\frac{1}{2}[\lambda.y,\lambda.y] = \frac{\lambda^2}{2}.[y,y] = - \frac{1}{2}[y,y] \quad \text{thus} \quad \lambda^2 = -1~, 
\]
which has no solution for $\lambda$ in $\mathbb{R}$. This reflects the fact that there are no $\mathbb{R}$-algebra morphisms from $\mathbb{C}$ to $\mathbb{R}$. However $\mathfrak{g}_{\mathbb{C}} \otimes \mathbb{C}$ has a two canonical Maurer--Cartan elements given by $i.y$ and $-i.y$, which reflects the fact that there are two $\mathbb{R}$-algebra endomorphisms of $\mathbb{C}$.
\end{Example}

On the other hand, given a curved absolute $\mathcal{L}_\infty$-algebra $\mathfrak{g}$, one can construct a functor from dg Artinian algebras to simplicial sets. 

\begin{Definition}[Deformation functor]
Let $\mathfrak{g}$ be a curved absolute $\mathcal{L}_\infty$-algebra. Its \textit{deformation functor}
\[
\mathrm{Def}_{\mathfrak{g}}: \mathsf{dg}~\mathsf{Art}\text{-}\mathsf{alg}_{\geq 0} \longrightarrow \mathsf{sSet}
\]
is given by 
\[
\mathrm{Def}_{\mathfrak{g}}(B) \coloneqq \mathcal{R}(\mathrm{hom}(\mathrm{Res}_\varepsilon (-)^\circ, \mathfrak{g}))~.
\]
\end{Definition}

\begin{Remark}
This construction is analogous to the following classical construction: given a dg Lie algebra $\mathfrak{g}$ and a dg Artinian algebra $A$, one can consider the simplicial set given by
\[
\mathrm{MC}(\mathfrak{g} \otimes \overline{A})_\bullet~,
\]
where $\overline{A}$ is the augmentation ideal of $A$ and $\mathrm{MC}(-)_\bullet$ is the integration functor constructed in \cite{Hinich01}. This construction has been the classical way to associate to a dg Lie/$\mathcal{L}_\infty$-algebra a "deformation problem". 
\end{Remark}

\begin{Proposition}
Let $f: \mathfrak{h} \qi \mathfrak{g}$ be a weak-equivalence of complete curved absolute $\mathcal{L}_\infty$-algebras. It induces a natural weak-equivalence of simplicial sets $\mathrm{Def}_{f} :  \mathrm{Def}_{\mathfrak{h}} \qi \mathrm{Def}_{\mathfrak{g}}$ between the associated deformation functors.
\end{Proposition}

\begin{proof}
Let $f: \mathfrak{h} \qi \mathfrak{g}$ is a weak-equivalence of complete curved absolute $\mathcal{L}_\infty$-algebras. It induces a natural weak-equivalence of curved absolute $\mathcal{L}_\infty$-algebras $f_* : \mathrm{hom}(-,\mathfrak{h}) \qi \mathrm{hom}(-,\mathfrak{g})$. Indeed, the map $f$ is a weak-equivalence if and only if
\[
\widehat{\mathrm{B}}_\iota(f): \widehat{\mathrm{B}}_\iota \mathfrak{h}  \qi \widehat{\mathrm{B}}_\iota \mathfrak{g}
\]
is a quasi-isomorphism of $u\mathcal{CC}_\infty$-coalgebras. The following composition is a quasi-isomorphism
\[
\widehat{\mathrm{B}}_\iota(\mathrm{hom}(-,\mathfrak{h})) \cong \{-, \widehat{\mathrm{B}}_\iota \mathfrak{h} \} \qi \{-, \widehat{\mathrm{B}}_\iota \mathfrak{g} \}  \cong  \widehat{\mathrm{B}}_\iota(\mathrm{hom}(-,\mathfrak{g}))~,
\]
since $\{-,-\}$ is a right Quillen functor by Proposition \ref{prop: monoidal model category} and thus preserves weak-equivalences between fibrant objects, guaranteeing that $\{-,\widehat{\mathrm{B}}_\iota(f)\}$ is a quasi-isomorphism. One then concludes using the fact that $\mathcal{R}$ sends weak-equivalences to weak homotopy equivalences of simplicial sets by Theorem \ref{thm: propriétés de l'intégration}. 
\end{proof}

Let us explore other properties of the deformation functor.

\begin{lemma}\label{lemma: Def g preserves}
Let $\mathfrak{g}$ be a complete curved absolute $\mathcal{L}_\infty$-algebra. The functor 
\[
\mathcal{R}(\mathrm{hom}(-,\mathfrak{g})): \left(u\mathcal{CC}_\infty\text{-}\mathsf{coalg}\right)^{\mathsf{op}} \longrightarrow \mathsf{sSet}
\]
sends any homotopy colimit of $u\mathcal{CC}_\infty$-coalgebras to a homotopy limit of simplicial sets. Furthermore, it sends quasi-isomorphisms to weak-equivalences.
\end{lemma}

\begin{proof}
Let $\mathrm{hocolim}~ C_{\alpha}$ be a homotopy colimit of $u\mathcal{CC}_\infty$-coalgebras, we consider the following equivalences

\begin{align*}
\mathcal{R}(\mathrm{hom}(\mathrm{hocolim}~ C_{\alpha},\mathfrak{g})) &\simeq \mathrm{Hom}_{u\mathcal{CC}_\infty\text{-}\mathsf{coalg}}(C_*^c(\Delta^\bullet),  \{ \mathrm{hocolim}~ C_{\alpha},\widehat{\mathrm{B}}_\iota \mathfrak{g} \}) \\
&\simeq \mathrm{holim}~ \mathrm{Hom}_{u\mathcal{CC}_\infty\text{-}\mathsf{coalg}}(C_*^c(\Delta^\bullet),  \{C_{\alpha},\widehat{\mathrm{B}}_\iota \mathfrak{g} \}) \\
&\simeq \mathrm{holim}~ \mathcal{R}(\mathrm{hom}(C_{\alpha},\mathfrak{g}))~.
\end{align*}

\end{proof}

\begin{Proposition}
Let $\mathfrak{g}$ be a complete curved absolute $\mathcal{L}_\infty$-algebra. 

\medskip

\begin{enumerate}
\item The deformation functor $\mathrm{Def}_{\mathfrak{g}}$ sends quasi-isomorphisms to weak-equivalences of simplicial sets.

\medskip

\item The deformation functor $\mathrm{Def}_{\mathfrak{g}}$ preserves homotopy limits.
\end{enumerate}
\end{Proposition}

\begin{proof}
This follows directly from Lemma \ref{lemma: equivalence dg Artinien avec les uCC} and Lemma \ref{lemma: Def g preserves}.
\end{proof}

\begin{Remark}
Developing the theory of \textit{non-pointed formal moduli problems} will be the subject of a future work. For a brief heuristic exposition, see \cite[Chapter 3, Section 4, Conclusion]{mathez}.
\end{Remark}

\section*{Appendix}\label{Section: Appendix B}
In order to work with simpler algebraic structures, we recall the \textit{semi-augmented} bar construction introduced in \cite[Section 3.3]{HirshMilles12}, which provides us with smaller but non-functorial cofibrant resolutions for non-augmented dg operads. We introduce its dual construction for dg counital partial cooperads in the sense of \cite{lucio2022curved}.

\begin{Definition}[Semi-augmented dg operad]
A \textit{semi-augmented} dg operad $(\PP,\gamma,\eta,d_\PP,\epsilon)$ is the data of a dg operad $(\PP,\gamma,\eta,d_\PP)$ together with a degree $0$ morphism
\[
\epsilon: \PP \longrightarrow \I~,
\]
such that $\epsilon \circ \eta = \mathrm{id}_\I~.$ 
\end{Definition}

Let $\overline{\PP}$ be the graded $\mathbb{S}$-module given by $\overline{\PP} \coloneqq \mathrm{Ker}(\epsilon)$. 

\begin{Definition}[Semi-augmented bar construction]
Let $(\PP,\gamma,\eta,d_\PP,\epsilon)$ be a semi-augmented dg operad. Its \textit{semi-augmented bar construction}, denoted by $\mathrm{B}^{\mathrm{s.a}}\PP$, is given by
\[
\mathrm{B}^{\mathrm{s.a}}\PP \coloneqq \left(\mathscr{T}^c(s\overline{\PP}),d_{\mathrm{bar}} \coloneqq d_1 + d_2 ,\Theta_{\mathrm{bar}}\right)~.
\] 
Here $\mathscr{T}^c(s\overline{\PP})$ denotes the cofree conilpotent coaugmented cooperad generated by the suspension of the graded $\mathbb{S}$-module $\overline{\PP}$. The pre-differential $d_{\mathrm{bar}}$ is given by the sum of two terms. The first term $d_1$ is given by the unique coderivation extending the map
\[
\begin{tikzcd}[column sep=4pc,row sep=0.5pc]
\mathscr{T}^c(s\overline{\PP}) \arrow[r,twoheadrightarrow]
&s\overline{\PP} \arrow[r,"sd_{\overline{\PP}}"]
&s\overline{\PP}~.
\end{tikzcd}
\]
The second term $d_2$ is given by the unique coderivation extending the map
\[
\begin{tikzcd}[column sep=4pc,row sep=0.5pc]
\mathscr{T}^c(s\overline{\PP}) \arrow[r,twoheadrightarrow]
&s^2(\overline{\PP} \circ_{(1)} \overline{\PP}) \arrow[r,"s^{-1}\overline{\gamma}_{(1)}"]
&s\overline{\PP}~,
\end{tikzcd}
\]
where $\overline{\gamma}_{(1)}$ is given by the composition 
\[
\begin{tikzcd}[column sep=4pc,row sep=0.5pc]
\overline{\PP} \circ_{(1)} \overline{\PP} \arrow[r,"\gamma_{(1)}"]
&s\PP \arrow[r,twoheadrightarrow]
& s\overline{\PP}~.
\end{tikzcd}
\]
The curvature $\Theta_{\mathrm{bar}}$ is given by 
\[
\begin{tikzcd}[column sep=3.5pc,row sep=0.5pc]
\Theta_{\mathrm{bar}}: \mathscr{T}^c(\overline{\PP}) \arrow[r,twoheadrightarrow]
&s^2(\overline{\PP} \circ_{(1)} \overline{\PP}) \arrow[r,"s^{-2}\gamma_{(1)}"]
&s\PP \arrow[r,"\epsilon"]
&\I~.
\end{tikzcd}
\]
The semi-augmented bar construction $\mathrm{B}^{\mathrm{s.a}}\PP$ forms a conilpotent curved cooperad.
\end{Definition}

It provides smaller but non-functorial resolutions. Indeed, it is only functorial with respect to morphisms of \textit{semi-augmented} dg operads.

\begin{theorem}[{\cite[Theorem 3.4.4]{HirshMilles12}}]
Let $(\PP,\gamma,\eta,d_\PP,\epsilon)$ be a semi-augmented dg operad. There is a quasi-isomorphism of dg operads
\[
\varepsilon_{\mathcal{P}}: \Omega \mathrm{B}^{\mathrm{s.a}}\PP \qi \PP~.
\]
\end{theorem}

Furthermore, it coincides in certain cases with the Boardmann-Vogt construction. 

\begin{theorem}\label{thm: Boardman-Vogt and cellular chains}
Let $\PP$ an operad in the category of cellular topological spaces. There is an isomorphism of dg operads
\[
\Omega \mathrm{B}^{\mathrm{s.a}} C_*^c(\PP,\kk) \cong C_*^c( \mathrm{W} \PP,\kk)~,
\]
\vspace{0.5pc}

where $C_*^c(-,\kk)$ denotes the cellular chain functor and where $\mathrm{W}(-)$ denotes the Boardmann-Vogt construction. 
\end{theorem}

\begin{proof}
This result is shown for reduced dg operads in \cite{BergerMoerdijk06}. See \cite{grignou2022mapping} for the extension to semi-augmented dg operads using this resolution. 
\end{proof}

\begin{Example}
Let $\mathrm{uCom}$ be the operad in the category of sets defined by 
\[
\mathrm{uCom}(n) \coloneqq \{*\}~,
\]
for all $n \geq 0$, together with the obvious action of $\mathbb{S}_n$ and the obvious operad structure. Then
\[
\Omega \mathrm{B}^{\mathrm{s.a}} \ucom \cong \Omega \mathrm{B}^{\mathrm{s.a}} C_*^c(\mathrm{uCom},\kk) \cong  C_*^c( \mathrm{W} \mathrm{uCom},\kk)~,
\]
where $\mathrm{uCom}$ is viewed as an operad in the category of cellular topological spaces by endowing it with the discrete topology. 
\end{Example}

\medskip

We introduce the dual version for semi-coaugmented dg counital partial cooperads. For the definition of complete curved absolute (partial) operads, we refer to \cite[Appendix]{lucio2022curved}.

\begin{Definition}[Semi-coaugmented dg counital partial cooperad]
A \textit{semi-coaugmented} dg counital partial cooperad $(\C,\{\Delta_i\},\epsilon,d_\C,\eta)$ amounts to the data of a dg counital partial cooperad $(\C,\{\Delta_i\},\epsilon,d_\C)$ together with a degree $0$ morphism
\[
\eta: \I \longrightarrow \C~,
\]
such that $\epsilon \circ \eta = \mathrm{id}_\I~.$ 
\end{Definition}

Let $\overline{\C}$ be the graded $\mathbb{S}$-module given by $\overline{\C} \coloneqq \mathrm{Ker}(\epsilon)$. 

\begin{Definition}[Semi-coaugmented complete cobar construction]
Let $(\C,\{\Delta_i\},\epsilon,d_\C,\eta)$ be a semi-coaugmented dg counital partial cooperad. Its \textit{semi-coaugmented complete cobar construction}, denoted by $\widehat{\Omega}^{\mathrm{s.c}}\C$, is given by
\[
\widehat{\Omega}^{\mathrm{s.c}}\C \coloneqq \left(\mathscr{T}^\wedge(s^{-1}\overline{\C}),d_{\mathrm{cobar}} \coloneqq d_1 -d_2 ,\Theta_{\mathrm{cobar}}\right)~.
\] 
Here $\mathscr{T}^\wedge(s^{-1}\overline{\C})$ denotes the completed tree monad applied to the desuspension of the graded $\mathbb{S}$-module $\overline{\C}$. The pre-differential $d_{\mathrm{cobar}}$ is given by the difference of two terms. The first term $d_1$ is given by the unique derivation extending the map
\[
\begin{tikzcd}[column sep=4pc,row sep=0.5pc]
s^{-1}\overline{\C} \arrow[r,"sd_{\overline{\C}}"]
&s^{-1}\overline{\C} \arrow[r,hookrightarrow]
&\mathscr{T}^\wedge(s^{-1}\overline{\C})~.
\end{tikzcd}
\]
The second term $d_2$ is given by the unique derivation extending the map
\[
\begin{tikzcd}[column sep=4pc,row sep=0.5pc]
s^{-1}\overline{\C} \arrow[r,"s\overline{\Delta}_{(1)}"]
&s^{-2}(\overline{\C} \circ_{(1)} \overline{\C}) \arrow[r,hookrightarrow]
&\mathscr{T}^\wedge(s^{-1}\overline{\C})
\end{tikzcd}
\]
where $\overline{\Delta}_{(1)}$ is given by the composition
\[
\begin{tikzcd}[column sep=4pc,row sep=0.5pc]
\overline{\C} \arrow[r,"\Delta_{(1)}"]
&\C \circ_{(1)} \C \arrow[r,twoheadrightarrow]
&\overline{\C} \circ_{(1)} \overline{\C}~.
\end{tikzcd}
\]
The curvature $\Theta_{\mathrm{cobar}}$ is given by 
\[
\begin{tikzcd}[column sep=3.5pc,row sep=0.5pc]
\Theta_{\mathrm{cobar}}: \I \arrow[r,"\epsilon"]
&\C \arrow[r,"s^{-2}\overline{\Delta}_{(1)}"]
&s^{-2}(\overline{\C} \circ_{(1)} \overline{\C}) \arrow[r,hookrightarrow]
&\mathscr{T}^\wedge(s^{-1}\overline{\C})~.
\end{tikzcd}
\]
The semi-coaugmented complete cobar construction $\widehat{\Omega}^{\mathrm{s.c}}\C$ forms a complete curved augmented absolute operad. 
\end{Definition}

\begin{Proposition}\label{lemma: Bucom dual lin}
Let $(\PP,\gamma,\eta,d_\PP,\epsilon)$ be a semi-augmented dg operad which is arity and degree-wise finite dimensional. There is an isomorphism of complete curved augmented absolute operads
\[
\left(\mathrm{B}^{\mathrm{s.a}}\PP \right)^* \cong \widehat{\Omega}^{\mathrm{s.c}}\PP^*~.
\]
\end{Proposition}

\begin{proof}
We first view $\PP$ as a semi-augmented dg unital partial operad. Since it is arity and degree-wise finite dimensional, its linear dual is a semi-coaugmented dg counital partial cooperad. The result follows from a straightforward computation.
\end{proof}

Let us first make explicit what a curved algebra over $\widehat{\Omega}^{\mathrm{s.c}}\ucom$ is. Let us recall the definition of a \textit{classical} curved $\mathcal{L}_\infty$-algebra. 

\begin{Definition}[Curved $\mathcal{L}_{\infty}$-algebra]\label{def: classical curved linfty alg}
Let $\mathfrak{g}$ be a graded module. A \textit{curved} $\mathcal{L}_{\infty}$-\textit{algebra} structure $\{l_n\}_{n \geq 0}$ on $\mathfrak{g}$ is the data of a family of symmetric morphisms $l_n: \mathfrak{g}^{\odot n} \longrightarrow \mathfrak{g}$ of degree $-1$ such that, for all $n \geq 0$, the following equation holds:
\[
\sum_{p+q= n+1} \sum_{\sigma \in \mathrm{Sh}^{-1}(p,q)}(l_{p} \circ_1 l_q)^{\sigma} = 0~,
\]
where $\mathrm{Sh}^{-1}(p,q)$ denotes the inverse of the $(p,q)$-shuffles. The morphism $l_0 : \kk \longrightarrow \mathfrak{g}$ is equivalent to an element $\vartheta$ in $\mathfrak{g}_{-1}$ which is called the \textit{curvature} of $\mathfrak{g}$. 
\end{Definition}

Recall that any complete curved augmented absolute operad structure is in particular a curved operad in the sense of \cite{lucio2022curved}, hence one can consider curved algebras over the underlying curved operad. 

\begin{Proposition}\label{Prop: Curved algebras over OmegauCom}
Let $(V,d_V)$ be a pdg module. The data of a curved $\widehat{\Omega}^{\mathrm{s.c}}\ucom$-algebra structure on $V$ is equivalent to a family of symmetric operations 
\[
\left\{ l_n: V^{\odot n} \longrightarrow V \right\}
\]
of degree $-1$ for $n \neq 1$. When setting $l_1 \coloneqq d_V$, the family $\{l_n\}_{n \geq 0}$ forms a curved $\mathcal{L}_\infty$-algebra structure in the classical sense.
\end{Proposition}

\begin{proof}
This follows from a straightforward computation.
\end{proof}

\begin{Remark}[Mixed curved $\mathcal{L}_\infty$-algebras] 
One may wonder what type of algebraic structure would have appeared had we chosen to use the bar-cobar constructions that appear in \cite{lucio2022curved} instead of their "semi-(co)augmented" counterparts. The type of algebraic structure that appears is not substantially different. A \textit{mixed curved} $\mathcal{L}_\infty$-algebra $(\mathfrak{g},\{l_n\}_{n \geq 0},d_\mathfrak{g})$ amounts to the data of a pdg module $(\mathfrak{g},d_\mathfrak{g})$ together with a family $\{l_n\}_{n \geq 0}$ of symmetric morphisms $l_n: \mathfrak{g}^{\odot n} \longrightarrow \mathfrak{g}$ of degree $-1$ such that 
\[
\sum_{p+q= n+1} \sum_{\sigma \in \mathrm{Sh}^{-1}(p,q)}(l_{p} \circ_1 l_q)^{\sigma} = -\partial(l_n)~,
\]
for all $n \geq 0$. In particular, for $n=1$, the relation satisfied is
\[
d_\mathfrak{g}^2 - \partial(l_1) = l_1 \circ_1 l_1 + l_2 \circ_1 l_0~.
\]
It is immediate that $d_\mathfrak{g} - l_1$ together with the other structure maps does define a curved $\mathcal{L}_\infty$-algebra structure. Furthermore, this gives a forgetful functor
\[
\begin{tikzcd}[column sep=4.5pc,row sep=0.3pc]
\mathsf{mix} ~ \mathsf{curv}\text{-}s\mathcal{L}_\infty \text{-}\mathsf{alg} \arrow[r,"\mathrm{blend}"]
&\mathsf{curv}\text{-}\mathcal{L}_\infty \text{-}\mathsf{alg}~.
\end{tikzcd}
\]
This type of structures can be found in \cite{calaque2021lie}. Considering curved algebras over the curved cooperad $\mathrm{B}\ucom$ would result in "mixed curved absolute" $\mathcal{L}_\infty$-algebras.
\end{Remark}

\bibliographystyle{alpha}
\bibliography{bibe}
\end{document}